\documentclass[11pt]{sat}

\newif\ifsattoc\sattoctrue
\newread\testfl\immediate\openin\testfl=\jobname.toc
    \ifeof\testfl\sattocfalse\fi\closein\testfl

\usepackage{float}
\usepackage[breaklinks]{hyperref}

\usepackage[lite,backrefs,alphabetic]{amsrefs}

\usepackage{amsmath,amssymb,latexsym,amsthm,enumerate,graphicx,tabularx,upgreek,mathabx}

\usepackage{ltablex}
\usepackage{paralist}

\renewcommand{\figurename}{Table}

\def\norm#1#2{\ensuremath{\left\|#1\right\|_{#2}}}

\newcommand{\ds}{\displaystyle}
\newcommand{\C}{C}

\renewcommand{\P}{\mathbb{P}}
\newcommand{\B}{\mathbb{B}}
\newcommand{\N}{\mathbb{N}}
\newcommand{\R}{\mathbb{R}}
\newcommand{\Y}{\mathbb{Y}}
\newcommand{\W}{\mathbb{W}}
\newcommand{\ie}{{\em i.e.,} }
\newcommand{\eg}{{\em e.g.,}{} }
\newcommand{\st}{\; : \; }
\newcommand{\hboxx}{ \mbox}
\newcommand{\LL}{\mathbb{L}}
\newcommand{\Lp}{\LL_p}
\newcommand{\ineq}[1]{(\ref{#1})}
\newcommand{\thm}[1]{Theorem~\ref{#1}}
\newcommand{\prop}[1]{Statement~\ref{#1}}
\newcommand{\props}[2]{Statements~\ref{#1} and \ref{#2}}
\newcommand{\newfig}[1]{Table~\ref{#1}}
\newcommand{\newfigs}[2]{Tables~\ref{#1} and \ref{#2}}
\newcommand{\op}[1]{Open Problem~\ref{#1}}
\newcommand{\cor}[1]{Corollary~\ref{#1}}

\newcommand{\NN}{\mathcal{N}}
\newcommand{\oversim}[1]{\stackrel{#1}{\sim}}

\newcommand{\esssup}{\mathop{\rm ess \; sup}}
\newcommand{\ww}{{\mathrm{w}}}
\newcommand{\fra}{f_{\mu,\lambda}}
\newcommand{\gml}{g_{\mu,\lambda}}

\newcommand{\circplusnew}{\circplus}

\newcounter{subfig}

\newcommand{\samefigure}{
\renewcommand{\thefigure}{\arabic{figure}\alph{subfig}}
\setcounter{subfig}{1} \addtocounter{figure}{-1} }
\newcommand{\esamefigure}{
\renewcommand{\thefigure}{\arabic{figure}}
}

\newtheorem{theorem}{Theorem}
\newtheorem{corollary}[theorem]{Corollary}
\newtheorem{problem}{Open Problem}

\theoremstyle{definition}
\newtheorem{proposition}{Statement}

\theoremstyle{remark}
\newtheorem*{remark}{Remark}

\newcommand\updots{\rotatebox[origin=lb]{-45}{\scalebox{1.3}{$\vdots$}}}

\def\versiondate{\today}

\title{Uniform and Pointwise Shape Preserving Approximation by Algebraic Polynomials}
\def\shorttitle{Shape Preserving Approximation $($SPA$)$}

\author{K. A. Kopotun\footnote{The first and third authors were supported in part by NSERC of Canada.},
D. Leviatan\footnote{Part of this work was done while the second and
fourth authors visited the University of Manitoba.}, A.
Prymak\footnotemark[1] {} and I. A. Shevchuk\footnotemark[2]}
\def\shortauthor{K. A. Kopotun, D. Leviatan, A. Prymak and I. A. Shevchuk}

\def\abstracttext
{We survey developments, over the last thirty years, in the theory
of Shape Preserving Approximation (SPA) by algebraic polynomials on
a finite interval. In this article, ``shape'' refers to (finitely
many changes of) monotonicity, convexity, or $q$-monotonicity of a
function (for definition, see Section 4). It is rather well known
that it is possible to approximate a function by algebraic
polynomials that preserve its shape (\ie the Weierstrass
approximation theorem is valid for  SPA). At the same time, the
degree of SPA is much worse than the degree of best unconstrained
approximation in some cases, and it is ``about the same'' in others.
Numerous results quantifying this difference in degrees of SPA and
unconstrained approximation have been obtained in recent years, and
the main purpose of this article is to provide a ``bird's-eye view''
on this area, and discuss various approaches used.

In particular, we present results on the validity and invalidity of
uniform and pointwise estimates in terms of various moduli of
smoothness. We compare various constrained and unconstrained
approximation spaces as well as orders of unconstrained and shape
preserving approximation of particular functions, etc.
There are quite a few interesting phenomena and several open
questions. }

\def\MSCnumbers{41A10, 41A17, 41A25, 41A29} % see http://www.ams.org/msc/

\def\keywords{Shape Preserving Approximation (SPA), degree of approximation by algebraic
polynomials, Jackson-Stechkin and Nikolskii type estimates, moduli
of smoothness,
Ditzian-Totik weighted moduli of smoothness, monotonicity, convexity and $q$-monotonicity}

\def\startpagenumber{24}
\def\volumenumber{6} % = current year - 2011
\def\year{2011}

\def\dword#1{{\bf #1}}\def\eword#1{{\it #1}}
\def\dd{\,{\rm d}}  % for integration (making them mathop places them less well)
\def\ee{{\rm e}}  % for the base of the natural log
\newcommand{\beginddoc}{
\maketitle
\begin{abstract}
\abstracttext \vskip1pt MSC: \MSCnumbers
\ifx\keywords\empty\else\vskip1pt Keywords: \keywords\fi
\end{abstract}
\insert\footins{\scriptsize
\medskip
\baselineskip 8pt \leftline{Surveys in Approximation Theory}
\leftline{Volume \volumenumber, \year.
pp.~\thepage--\pageref{endpage}.} \leftline{\copyright\ \year\
Surveys in Approximation Theory.} \leftline{ISSN 1555-578X}
\leftline{All rights of reproduction in any form reserved.}
\smallskip
\par\allowbreak}
\ifsattoc\else\tableofcontents\fi}

\renewcommand\rightmark{\ifodd\thepage{\it \hfill\shorttitle\hfill}\else {\it \hfill\shortauthor\hfill}\fi}
\def\endddoc{\label{endpage}\end{document}}
\date{{\small \versiondate}}
\begin{document}
\beginddoc
\ifsattoc
\bigskip
%%%%%%%%%%%%%%%%%%%%%%%%%%%  toc material
\def\toczer{0}\def\tochalf{.5}\def\tocone{1}
\def\tocindent{0}
\def\ection{section}\def\ubsection{subsection}
\def\numberline#1{\hskip\tocindent truecm{} #1\hskip1em}
\newread\testfl
\def\inputifthere#1{\immediate\openin\testfl=#1
    \ifeof\testfl\message{(#1 does not yet exist)}
    \else\input#1\fi\closein\testfl}
\countdef\counter=255
\def\diamondleaders{\global\advance\counter by 1
  \ifodd\counter \kern-10pt \fi
  \leaders\hbox to 15pt{\ifodd\counter \kern13pt \else\kern3pt \fi
  \hss.\hss}\hfill}
\newdimen\lextent
\newtoks\writestuff
\medskip
\begingroup
\small
\def\contentsline#1#2#3#4{
\def\argu{#1}
\ifx\argu\ection\let\tocindent\toczer\else
\ifx\argu\ubsection\let\tocindent\tochalf\else\let\tocindent\tocone\fi\fi
\setbox1=\hbox{#2}\ifnum\wd1>\lextent\lextent\wd1\fi}
\lextent0pt\inputifthere{\jobname.toc}\advance\lextent by 2em\relax
\def\contentsline#1#2#3#4{
\def\argu{#1}
\ifx\argu\ection\let\tocindent\toczer\else
\ifx\argu\ubsection\let\tocindent\tochalf\else\let\tocindent\tocone\fi\fi
\writestuff={#2}
\centerline{\hbox to \lextent{\rm\the\writestuff%
\ifx\empty#3\else\diamondleaders{}
\hfil\hbox to 2 em\fi{\hss#3}}}}
\inputifthere{\jobname.toc}\endgroup
\immediate\openout\testfl=\jobname.toc % to empty \jobname.toc in order to
\immediate\closeout\testfl             % get \tableofcontents to initiate
\renewcommand{\contentsname}{}         % regeneration the toc file without
\tableofcontents\bigskip\bigskip               % printing a ToC.
\fi

\section{Scope of the survey}
This is a comprehensive survey   of  \eword{uniform} and
\eword{pointwise} estimates of  \eword{polynomial} shape-preserving
approximation (\dword{SPA}) of the following types:
\begin{inparaenum}[(i)]
\item  monotone,
\item convex,
\item $q$-monotone,
\item como\-notone,
\item coconvex.
\end{inparaenum}

The following types of SPA are not covered:
\begin{inparaenum}[(i)]
\item co-$q$-monotone (no positive results are known for $q \geq 3$),
\item estimates of SPA in the $\Lp$ (quasi)norms,
\item positive, copositive, onesided and intertwining approximation,
\item nearly shape preserving approximation,
\item pointwise estimates of SPA with interpolation at the
endpoints,
\item SPA by splines with fixed and free knots,
\item SPA by rational functions,
\item simultaneous SPA,
\item SPA of periodic functions.
\end{inparaenum}

\renewcommand\listfigurename{\section{List of Tables}}

\begingroup\small
\listoffigures
\endgroup

\section{Motivation} \label{motivation}

The purpose of this section is to provide some motivation to the
problems/questions discussed in this survey. Because convex
functions are ``nicer'' than monotone ones (for example, a convex
function defined on an open interval is necessarily locally
absolutely continuous on that interval), we pick ``convexity'' as
our example of ``shape'' in the examples appearing in this section.
Readers who are not familiar with the notations, definitions and/or
symbols that we use here can consult Section~\ref{S1} ``Definitions,
Notations and Glossary of Symbols''.

With $E_n(f)$ denoting the degree of approximation of $f$ by
polynomials of degree $<n$, it was shown by  Bernstein in 1912 that
\[
    E_n(|x|) \sim  n^{-1}, \quad n\in\N   .
\]
(The upper estimate was obtained in 1908 by de la Vall\'ee Poussin.)
Since the polynomial of best approximation to $f$ necessarily
interpolates $f$ at ``many'' points it is clear that despite the
fact that $|x|$ is convex on $[-1,1]$, its polynomial of best
approximation of degree $<n$ is not going to be convex for $n\geq 5$.

A natural question is then if one can approximate $|x|$ by a convex
polynomial with the error of approximation still bounded by $c
n^{-1}$ (the answer is ``yes'' -- we omit almost all references in
this section for the reader's convenience, since  the results
mentioned here are either   simple or   particular cases of general
theorems discussed in detail later in this survey).
Is this the case for other functions having the same approximation
order? In other words, if a convex function $f$ is such that $E_n(f)
=  O\left( n^{-1}\right)$, does this imply that there exist convex
polynomials providing the same rate of approximation, \ie
$E_n^{(2)}(f) =  O\left( n^{-1}\right)$ ?  (The answer is ``yes''.)
Is this the case for other rates of convergence of errors to zero?
For example, if a convex function $f$ is such that $E_n(f) = O\left(
n^{-\alpha}\right)$ for some $\alpha > 0$, does this imply that
$E_n^{(2)}(f) =  O\left( n^{-\alpha}\right)$?  (The answer is
``yes''.)
Does this mean that the rates of approximation of convex functions
by all polynomials  and by convex polynomials   from $\P_n$ are
``about the same''? In other words, because it is clear that $E_n(f)
\leq  E_n^{(2)}(f)$ (since the set of all convex polynomials from
$\P_n$ is certainly a proper subset of $\P_n$), is it possible that
the converse is also true, \ie there exists a constant $c(f)$ such
that $E_n^{(2)}(f) \leq c(f) E_n(f)$ for all $n$? (The answer is
``no''.)

So, approximation by convex polynomials cannot simply be reduced to
unconstrained approximation.

Going back to  a convex function $f$ such that $E_n(f) =  O\left(
n^{-\alpha}\right)$, we know that this implies that $E_n^{(2)}(f) =
O\left( n^{-\alpha}\right)$, but these are asymptotic rates that
only deal with the behavior of the quantities $E_n(f)$ and
$E_n^{(2)}(f)$ as $n\to\infty$. Thus, we conclude that
$E_n^{(2)}(f)\le cn^{-\alpha}$, $n\ge m$, for some large $m$. What
can be said about $E_n^{(2)}(f)$ for small values of $n$? Clearly,
we can say something about it even when we do not have asymptotic
information on $E_n(f)$. Indeed,
\[
E_1 (f)=E_1^{(2)}(f)=(\max_{x\in [-1,1]}f(x)-\min_{x\in
[-1,1]}f(x))/2 \le\|f\|
\]
implies that, for any $m\in\N$,
\[
E_n^{(2)}(f)\le E_1^{(2)}(f)\le(\|f\|  m^\alpha)n^{-\alpha},\quad
1\le n\le m.
\]
Hence, for a convex function $f$ such that $E_n(f)=O(n^{-\alpha})$,
we have
\begin{equation} \label{newfirst}
E_n(f)=O\left(n^{-\alpha}\right) \quad \Longrightarrow \quad
E_n^{(2)}(f)\leq c(f)n^{-\alpha}, \quad n\in\N .
\end{equation}

It is   easy to see that this  implication  is no longer valid if we
require that the constant $c$ in \ineq{newfirst} remain independent
of $f$ (for example, for  $f_\gamma(x) = \gamma x^{2}$, it is
obvious that $E_n(f_\gamma) =  O\left(n^{-\alpha}\right)$ but $E_{1}
(f_\gamma) = E_{1}^{(2)}(f_\gamma) = |\gamma|/2 = \|f_\gamma\|/2
\to\infty$ as $\gamma\to\infty$).

We should therefore word this question more carefully. Namely,
suppose that $f$ is a convex function   such that $E_n(f) \leq
n^{-\alpha}$ for $n\in\N$, does this imply that $E_n^{(2)}(f) \leq
c(\alpha) n^{-\alpha}$ for all $n\in\N$? (The answer is ``yes''.)
What if we only know that $E_n(f) \leq n^{-\alpha}$ for $n\geq
2011$, can we say that $E_n^{(2)}(f) \leq c(\alpha) n^{-\alpha}$ for
$n\geq 2011$? (The answer is ``no''.) So, the fact that $E_n(f) \leq
n^{-\alpha}$ for $n\geq\NN$  sometimes implies that $E_n^{(2)}(f)
\leq c(\alpha) n^{-\alpha}$ for  $n\geq\NN$, and sometimes does not.
Do we know precisely for which $\NN\in\N$ this implication is valid,
and for which it is invalid in general? (The answer is ``yes,
we do know: if $\alpha \leq 4$, this implication is valid for
$\NN\leq 4$, and it is in general invalid in the stated form for
$\NN\geq 5$; if $\alpha > 4$, then this implication is valid for all
$\NN\in\N$''.)

In general, given a function $f\in\C^r[-1,1]$, in order to obtain an estimate of the form $E_n(f) \leq c
n^{-\alpha}$ one can try to apply a Jackson-Stechkin type uniform
estimate
\begin{equation}\label{tmp1}
E_n (f)\le c(k,r) n^{-r}\omega_k \bigl(f^{(r)},n^{-1}\bigr), \quad
n\ge k+r  .
\end{equation}
For example, since $|x|\in\C[-1,1]$ and $\omega (|x|, \delta) =
\min\{1, \delta\}$, using \ineq{tmp1} with $r=0$ and $k=1$, one
immediately obtains the above mentioned de la Vall\'ee Poussin
estimate $E_n (|x|) \leq c n^{-1}$, $n\in\N$. Hence,  it is
desirable to have  analogs of Jackson-Stechkin type estimates for
approximation by convex polynomials. So, does one have an estimate
$E_n^{(2)} (f)\le  c \omega \left(f, n^{-1} \right)$, $n\in\N$,
which would then immediately imply that $E_n^{(2)} (|x|) \leq c
n^{-1}$, $n\in\N$? (The answer is ``yes''.)

One can naturally try to use the same process for functions having
greater smoothness (and, hence, having better approximation rates).
For example, it is known (see
Timan \cite{T1994}*{p. 416},
Bernstein \cite{collected2}*{pp. 262--272, pp. 402--404} and Ibragimov
\cite{ibr}) that, for the function $\fra(x) = x^{\mu-1}
|x|^{1+\lambda}$, $\mu\in\N$, $|\lambda|<1$,
 \[
 \lim_{n\to \infty} n^{\mu+\lambda} E_n(\fra) = c(\mu,\lambda) > 0 .
 \]
In particular, this implies that, if $\alpha>0$ is such that
$\alpha/2\not\in\N$ (otherwise, the function becomes a polynomial),
then
\[
E_n(|x|^\alpha) \oversim{\alpha}n^{-\alpha},\quad n\in\N.
\]
An upper estimate can be obtained directly from \ineq{tmp1} by
setting $r=0$, $k  = \lceil\alpha\rceil$, and using the fact that
\[
\omega_k(|x|^\alpha, \delta) \leq c(\alpha)  \min\{1,
\delta^\alpha\}.
\]
Now, since $|x|^{2011}$ is convex, we would get the estimate
$E_n^{(2)}(|x|^{2011}) \leq c  n^{-2011}$ if the estimate $E_n^{(2)}
(f)\le  c \omega_{2011} \left(f, n^{-1} \right)$ were true for all
convex functions. Unfortunately, this is not the case, and it is
known that even the estimate $E_n^{(2)} (f)\le  O\left( \omega_{5}
(f, n^{-1})\right)$  is invalid for some convex functions (and the
estimate $E_n^{(2)} (f)\le  c \omega_{4} \left(f, n^{-1} \right)$
cannot be valid for all convex functions $f$ and all
  $n\geq \NN$ with the constant $\NN$ independent of $f$). In this
particular case, one can overcome this obstacle  by using the
estimate $E_n^{(2)}(f)\le c(k)n^{-2}\omega_{k}\left(f'', n^{-1}
\right)$ which is true for all $k\in\N$ and all convex functions $f$
from $\C^2[-1,1]$, and which implies the required result.

Can this approach always be used? In other words, is it enough to
use Jackson-Stechkin type estimates in order to get the ``right''
order of polynomial approximation? Again, unfortunately, the answer
is ``no''. It has been known for a long time that the
Jackson-Stechkin type estimates, while producing the ``right''
orders of approximation for some functions, produce rather weak
estimates in some special cases. For example, Ibragimov
\cite{ibr1946} verified Bernstein's conjecture (see
\cite{ber1937}*{p. 91}) that, for $g(x)=(1-x)\ln(1-x)$ (we will later call
this function $g_{1,1}$),
\[
  E_n(g) \sim n^{-2} \, .
\]
Because $g\not\in\C^1[-1,1]$, one can only use \ineq{tmp1} with
$r=0$, and since, for any $k\geq 2$, $\omega_k(g,n^{-1}) \oversim{k}
n^{-1}$,
one can only get from \ineq{tmp1} the rather weak estimate $E_n(g)
\leq c n^{-1}$ . In particular, this means that, given a function
$f\in\C[-1,1]$ such that   $E_n(f) = O(n^{-2})$, in general we
cannot use Jackson-Stechkin type estimates in order to prove
$E_n^{(2)}(f) = O(n^{-2})$, and more precise estimates (yielding
constructive characterization of function classes) are required.
This brings us to the Ditzian-Totik type estimates (see \cites{DT,
Dit2007}). In particular, since the following Ditzian-Totik type
estimate is valid
\begin{equation}\label{tmp2}
E_n (f)\le c(k) \omega_k^\varphi \left(f, n^{-1} \right), \quad n\ge
k   ,
\end{equation}
and since
\[
\omega_2^\varphi \left(g, n^{-1} \right) \sim n^{-2} ,
\]
one immediately obtains   $E_n(g) \leq c n^{-2}$. Now, since $g$ is
convex on $[-1,1]$, how does one go about proving the estimate
$E_n^{(2)}(g) \leq c n^{-2}$ (which, as we mentioned above, has to
be valid)? This estimate would immediately follow  if the estimate
$E_n^{(2)} (f)\le c  \omega_2^\varphi \left(f, n^{-1} \right)$ were
true for all convex $f$, but is it true? (The answer is ``yes''.)

Suppose now that instead of the function $g$ we consider a similar
function having higher smoothness (and, hence, a better rate of
polynomial approximation). For example,
Ibragimov \cite{ibr1946}*{Theorem VII} showed that, for $\gml(x)  = (1-x)^\mu \ln^\lambda
(1-x)$, $\nu$, $\lambda\in\N$,
\begin{equation}\label{ibra}
E_n(\gml) \oversim{\lambda, \mu} \frac{(\ln n)^{\lambda
-1}}{n^{2\mu}},\quad n\ge2.
\end{equation}
In particular, consider $h(x) = 10x^2-(1-x)^3 \ln(1-x)$ which is
convex on $[-1,1]$. Since
\[
\omega_6^\varphi (h, n^{-1}) \sim n^{-6}
\]
the estimate \ineq{tmp2} implies that $E_n(h) \leq c n^{-6}$. Now, a
similar estimate for the rate of convex approximation would
immediately follow if the estimate $E_n^{(2)} (f)\le c
\omega_6^\varphi \left(f, n^{-1} \right)$ were valid for all convex
$f$, but is it valid? (The answer is ``no''.) In the above example
for  the Jackson-Stechkin type estimates we managed to overcome a
similar difficulty by using an estimate involving the modulus of a
derivative of $f$. What if we try a similar approach now? We know
that
 \begin{equation}\label{badder}
\omega_6^\varphi \left(f, \delta \right) \leq c \delta^2
\omega_4^\varphi (f'',\delta)
\end{equation}
and the function $h$ is clearly in $\C^2[-1,1]$. Hence, if we show
that $\omega_4^\varphi (h'',n^{-1}) \leq c n^{-4}$, and if the
estimate $E_n^{(2)} (f)\le c  n^{-2} \omega_4^\varphi \left(f'',
n^{-1} \right)$ is valid for all convex functions $f\in\C^2[-1,1]$,
then we will have proved what we want. Now, since $h''(x) = 15+5x-6
(1-x)\ln(1-x)$, we have
\[
\omega_4^\varphi (h'',n^{-1}) = c \omega_4^\varphi (g_{1,1}, n^{-1})
\leq c n^{-2} ,
\]
and the last inequality cannot be improved (this can be verified
directly, but it also follows immediately from the fact that if we
could replace $n^{-2}$ by $o(n^{-2})$, then \ineq{ibra} would not be
valid). Hence, this approach would only give us a rather weak
estimate (different by the factor of $n^2$ from the optimal). The
reason why this approach fails is that the inequality \ineq{badder}
is very imprecise, and one needs to work with generalized
Ditzian-Totik moduli involving derivatives of functions. In this
particular case, we can obtain the needed result taking into account
the fact that
\[
\omega_6^\varphi \left(f,\delta \right) \leq c\delta^5
\omega_{1,5}^\varphi (f^{(5)},\delta) ,
\]
and using the estimate
\[
E_n^{(2)} (f)\le c  n^{-5} \omega_{1,5}^\varphi \bigl(f^{(5)},
n^{-1}\bigr) ,
\]
and the inequality $\omega_{1,5}^\varphi \left(h^{(5)},\delta\right)
\leq c\delta$.

The conclusion that we can reach at this time is that one needs to
work with generalized Ditzian-Totik moduli of smoothness
$\omega^\varphi_{k,r}(f^{(r)},\delta)$ instead of the regular
(ordinary or Ditzian-Totik) moduli in order to get exact uniform
estimates. Of course, it is now well known that one can also obtain
exact estimates of approximation by algebraic polynomials (yielding
constructive characterization of classes of functions) in terms of
the ordinary moduli of smoothness, but only if pointwise estimates
are used. It is well known that for a function $f\in \C^r[-1,1]$ and
each $n\ge k+r$ there is a polynomial $P_n\in\P_n$, such that for
every $x\in[-1,1]$,
\[
|f(x)-P_n(x)|\le
c(k,r)\rho_n^r(x)\omega_k\bigl(f^{(r)},\rho_n(x)\bigr),
\]
where $\rho_n(x):=n^{-2} +n^{-1} \sqrt{1-x^2}$. (This is a so-called
Nikolskii type pointwise estimate.)  This immediately implies
\ineq{tmp1}, and is clearly more precise than \ineq{tmp1}. Hence,
the same questions that we discussed above can also be asked for
pointwise estimates of rates of approximation by convex polynomials.

So far, we only discussed functions which are convex on the entire
interval $[-1,1]$. What if a function has one or more inflection
points? Can we approximate it by polynomials having the same
``shape''? Again, all of the above questions naturally arise. For
example, consider the above mentioned function $g_{2,1}(x) = (1-x)^2
\ln (1-x)$. This function is convex on $[-1, y_1]$ and concave on
$[y_1,1]$, where $y_1 = 1-\ee^{-3/2}$. We already know (see
\ineq{ibra}) that $E_n(g_{2,1}) \leq c n^{-4}$. Do we have the same
rate of coconvex approximation, \ie is $E_n^{(2)}(g_{2,1}, Y_1) \leq
c n^{-4}$ valid? (The answer is ``yes''.) Is it always the case that
if a function $f$ changes its convexity at one point of $[-1,1]$ and
has the rate of its unconstrained approximation bounded by $n^{-4}$,
then its rate of coconvex approximation is   bounded by $c n^{-4}$?
The answer is ``yes, but the constant $c$ has to depend on $y_1$,
and may become large as $y_1$ gets closer to $\pm 1$''. This seems
like a very obvious and intuitively clear conclusion (the closer the
inflection point is to the endpoints, the harder it is to control
the approximating polynomial preserving the ``shape'' of the
function). However, this is where intuition fails. We now know that
given a function $f$ that has one inflection point in $[-1,1]$ and
whose rate of unconstrained approximation is bounded by
$n^{-\alpha}$ with $\alpha\neq 4$, its rate of coconvex
approximation is bounded by $c(\alpha) n^{-\alpha}$. In other words,
in this sense, the case $\alpha=4$ is totally different from all
other cases.

The main purpose of this survey is to summarize the state of the art
of this area of research as of the summer of 2011.

\section{Definitions, Notations and Glossary of Symbols}\label{S1}
Most symbols used throughout this paper are listed in the table
below (however, for the reader's convenience, we also define some of
them the first time they are used in this survey).

\setlength{\extrarowheight}{5.0pt}

 \begin{tabularx}{\textwidth}{ |l|X| }

  \hline
$\N$ & $\{1, 2, 3, \dots\}$\\
$\N_0$ & $\N \cup \{0\}$ \\
 $\C(S)$ & space of  continuous functions on $S$ \\
 $\C^r(S)$ & space of $r$-times continuously differentiable functions on $S$, $r\in\N$ \\
$ \norm{f}{\C(S)}$ & $\max_{x\in S} |f(x)| $\\
  $\norm{\cdot}{}$ & $ \norm{\cdot}{\C[-1,1]}$\\
  $\norm{f}{\LL_\infty(S)}$ & $\esssup_{x\in S} |f(x)|$\\
 $\Delta^{(1)}$ & $   \left\{ f\in \C[-1,1] \st \mbox{\rm $f$ is nondecreasing on $[-1,1]$} \right\}$\\
$\Delta^{(2)}$ & $   \left\{ f\in \C[-1,1] \st \mbox{\rm $f$ is convex on $[-1,1]$} \right\}$\\
$\Delta^{(q)}$ & $   \left\{ f\in \C[-1,1]\cap \C^{q-2}(-1,1) \st
\mbox{\rm $f^{(q-2)}$
is convex on $(-1,1)$} \right\}$, $q\geq 3$; $f$ is $q$-monotone\\
$\P_n$ & space of algebraic polynomials of degree $\leq n-1$\\
$ E_n(f) $ & $ \inf_{P_n\in\P_n}\|f-P_n\|$
(degree of best unconstrained approximation) \\
$ E_n^{(q)}(f)$ & $ \inf_{P_n\in\P_n\cap\Delta^{(q)}}\|f-P_n\| $
(degree of $q$-monotone approximation)\\
$Y_s$ &  collection $\{y_i\}_{i=1}^{s}$ of $s\in\N$
points $ -1=: y_{s+1} <y_s< \dots <y_1< y_0 := 1 $\\
$Y_0$ &  $\emptyset$\\
$\Delta^{(1)}(Y_s)$ & set of all functions $f\in \C[-1,1]$ that
change  monotonicity at the points $Y_s$, and are non-decreasing on $[y_1, 1]$ \\
$\Delta^{(2)}(Y_s)$ & set of all functions $f\in \C[-1,1]$ that
change  convexity at the points $Y_s$, and are convex on $[y_1, 1]$ \\
$\Delta^{(q)}(Y_s)$ & set of all functions $f\in\C[-1,1]\cap
\C^{q-2}(-1,1)$ such that $f^{(q-2)}$ changes convexity at the
points $Y_s$, and is convex on $[y_1, 1]$,  $q\geq 3$ \\
$\Delta^{(q)}(Y_0)$ &  $\Delta^{(q)}$,  $q\geq 1$ \\
$E_n^{(q)}(f,Y_s)$ & $\inf_{P_n\in{\Delta^{(q)}(Y_s)}\cap\P_n}
\|f-P_n\|$ (degree of co-$q$-monotone approximation)\\
$\Delta^k_h(f,x)$ & $\sum_{i=0}^{k}{k\choose i}(-1)^{k-i}
f(x-kh/2+ih)$, if $|x\pm kh/2|< 1$, and
$0$, otherwise ($k$th symmetric difference) \\
$\omega_k(f,t)$ & $ \sup_{0\le h\le t}\|\Delta ^k_h(f,\cdot)\|$
($k$th modulus of smoothness)\\
$\omega(f,t)$ & same as $\omega_1(f,t)$ (ordinary modulus
of continuity of $f$)\\
$\omega_0(f,t)$ & $  \|f\|_{\LL_\infty[-1,1]}$ \\
$\varphi(x)$ & $   \sqrt{1-x^2}$\\
$\rho_n(x)$ & $n^{-1}\sqrt{1-x^2} + n^{-2}$\\
$\omega^\varphi_{k}(f,t) $ & $ \sup_{0\le h\le t}
\|\Delta_{h\varphi(\cdot)}^k(f,\cdot)\|$ (Ditzian-Totik (D-T) $k$th
modulus of smoothness)\\ $\W^r$ &  space of functions $f$ defined on
$[-1,1]$ such that $f^{(r-1)}$ is absolutely continuous
and $f^{(r)}\in\LL_\infty[-1,1]$, $r\geq 1$ \\
$\B^r$ &  space of functions $f$ defined on $[-1,1]$ such that
$f^{(r-1)}$ is locally absolutely continuous in $(-1,1)$ and
$\varphi^r f^{(r)}\in \LL_\infty[-1,1]$, $r\geq 1$\\
$\C^0_\varphi$ & $ \C[-1,1]$\\
$\C^r_\varphi$ & $\left\{f\in\C^{r}(-1,1)\st\lim_{x\to\pm1}
\varphi^r(x)f^{(r)}(x)=0 \right\}$, $r\geq 1$\\
$K(x,\mu)$ & $ \varphi(|x|+\mu\varphi(x))$\\
$\omega^\varphi_{k,r}(f^{(r)},t) $ & $ \sup_{0\le h\le
t}\sup_{x:|x|+ kh \varphi(x)/2<1}
K^r(x,kh/2)|\Delta_{h\varphi(x)}^k(f^{(r)},x)|
$ (D-T generalized modulus of smoothness) \\
$\omega^\varphi_{0,r}(f^{(r)},t)$ & $\|\varphi^r  f^{(r)}  \|_{\LL_\infty[-1,1]}$  \\
$c$  & absolute positive constants that can be different even if they appear on the same line\\
$c(\cdot)$   &  positive constants that depend on the parameters appearing inside the parentheses and nothing else \\
$a_n \sim b_n$ & there exists an absolute positive constant $c$ such
that
$c^{-1} a_n \leq b_n \leq c a_n$ for all $n\in\N$\\
$a_n \oversim{\alpha_1, \alpha_2, \dots} b_n$ & there exists a
positive constant $c=c(\alpha_1, \alpha_2, \dots)$ such that $c^{-1}
a_n \leq b_n \leq c a_n$ for all $n\in\N$\\
$a_n(f) = O(b_n)$   & there exist $c(f)\in \R$ and $m\in\N$ such
that
$|a_n(f)| \leq c(f) |b_n|$, for  $n\geq m$\\
\hline
\end{tabularx}

\begin{remark}
Throughout this survey, the letters ``$q$'', ``$s$'', ``$r$'' and
``$k$'' always stand for nonnegative integers. The letter ``$q$'' is
always used to describe the shape of a function (\eg
$\Delta^{(q)}$), the letter ``$s$'' always stands for the number of
changes of monotonicity, convexity or $q$-monotonicity (\eg $\Y_s$,
$Y_s$, $\Delta^{(1)} (Y_s)$),  the letter ``$r$'' always refers to
the $r$-th derivative (\eg $\W^r$, $\B^r$, $f^{(r)}$,
$\omega^\varphi_{k,r}(f^{(r)},\delta)$), and the letter ``$k$'' is
used to describe the order of the appropriate moduli of smoothness
or similar quantities (\eg $\omega_k$, $\omega_{k}^\varphi$,
$\Phi^k$).
\end{remark}

The following are the types of estimates of errors of SPA and
unconstrained approximation that we discuss in this survey.
\newpage
\begin{tabularx}{\textwidth}{ |l|l|l| }
\hline Type of estimate & Estimate is given in terms of & Notation
\\  \hline Jackson-Stechkin (uniform) &  $
n^{-r}\omega_k\left(f^{(r)},1/n\right)$ &
$\delta_n(x) := 1/n$, $\ww_k :=\omega_k$ \\
Ditzian-Totik (uniform) & $ n^{-r}
\omega^\varphi_{k,r}\bigl(f^{(r)},1/n\bigr)$ & $\delta_n(x) := 1/n$,
$\ww_k :=\omega^\varphi_{k,r}$ \\
Nikolskii (pointwise) &
$\underbrace{\rho_n^r(x)\omega_k\bigl(f^{(r)},\rho_n(x)\bigr)}_{\displaystyle
\delta_n^r(x) \ww_k\bigl(f^{(r)},\delta_n(x)\bigr)}$ & $\delta_n(x)
:= \rho_n(x)$,
$\ww_k :=\omega_k$ \\
\hline
  \end{tabularx}

With the above notation,  it is convenient  to refer to all of the
above types of estimates at once. Namely, given $q$, $s$, $r$ and
$k$ as above, for collections $Y_s$ and functions $f\in
\Delta^{(q)}(Y_s)$ (assumed to have appropriate smoothness so that
$\ww_k(f^{(r)},t)$ is defined and finite), we formally write an
estimate for the error of (co)-$q$-monotone polynomial approximation
of $f$ (which may or may not be valid for $P_n\in\P_n\cap\Delta^{(q)}(Y_s)$):
\begin{equation} \label{eq0}
\left|f(x)-P_n(x)\right|\leq
 c(k,r,q)\delta_n^r(x)\ww_k\bigl(f^{(r)},\delta_n(x)\bigr),\quad n\ge\NN,
\end{equation}
and distinguish the following cases for the quadruple $(k,r,q,s)$.

\begin{description}
\item[Case ``$+$'' (``strongly positive case"):]
Inequality \ineq{eq0} holds with $\NN$ depending only on $k$, $r$,
$q$ and $s$.

\item[Case ``$\oplus$" (``weakly positive case"):] Inequality \ineq{eq0}
holds with a constant $\NN$ that, in addition, depends on the set
$Y_s$ (\ie location of points where $f$  changes its monotonicity,
convexity or $q$-monotonicity), and does not hold in general with
$\NN$ independent of $Y_s$. (Note that this case is only applicable
if $s\geq 1$ since, in the case $s=0$, $Y_0 = \emptyset$.)
\item[Case ``$\ominus$" (``weakly negative case"):] Inequality  \ineq{eq0}
holds with a constant $\NN$ that depends on the function $f$, and
does not hold in general with $\NN$ independent of $f$.
\item[Case ``$-$" (``strongly negative case"):]  Inequality  \ineq{eq0}
does not hold in general even with a constant $\NN$ that depends on
the function $f$.
This means that there exists a function $f\in \Delta^{(q)}(Y_s) \cap
\C^r[-1,1]$ such that, for every sequence of polynomials
$P_n\in\P_n\cap\Delta^{(q)}(Y_s)$,
\[
\limsup_{n\to\infty} \left\|\frac{f-P_n}{\delta_n^r
\ww_k\left(f^{(r)},\delta_n\right)}\right\|=\infty.
\]
\end{description}

We will also have the case ``$\oslash$'', which we do not discuss
here because, so far, it appears only when pointwise estimates for
coconvex approximation are considered (see Section~\ref{subS11-3},
page~\pageref{slashref}).

\begin{remark}
In all cases ``$+$" discussed in this survey, inequality \ineq{eq0}
holds with $\NN=k+r$. (This is the best possible case since
\ineq{eq0} certainly cannot hold if $\NN < k+r$ because
$\ww_k(f^{(r)},t)$ vanishes if  $f\in\P_{k+r}$.)
\end{remark}

Along with different direct estimates on the error of SPA, we will
also consider relations between degrees of best unconstrained
approximation and SPA, in particular, the so-called
$\alpha$-relations. Namely, given $q$ and $s$ as above, and
$\alpha>0$, for collections $Y_s$ and functions
$f\in\Delta^{(q)}(Y_s)$, we investigate the validity of the
implication
\begin{equation} \label{eq0-alpha}
n^\alpha E_n(f)\le1,  \;\;  n\ge \NN   \quad  \Longrightarrow \quad
n^\alpha E_n^{(q)}(f,Y_s)\le c(s,\alpha),\;\; n\ge \NN^*,
\end{equation}
distinguishing the following cases. (Note that, while we use the
same symbols as for direct estimates, their meaning for
$\alpha$-relations is different.)
\begin{description}
\item[Case ``$+$'' (``strongly positive case"):] Implication~\eqref{eq0-alpha}
is valid with $\NN^*$ depending on $\NN$, $\alpha$, $s$.
\item[Case ``$\oplus$" (``weakly positive case"):]
Implication~\eqref{eq0-alpha} is valid with $\NN^*$ depending on
$\NN$, $\alpha$, $s$ and $Y_s$, and is not valid in general with
$\NN$ independent of $Y_s$.
\item[Case ``$\ominus$" (``weakly negative case"):] Implication~\eqref{eq0-alpha}
is valid with  $\NN^*$ that depends on the function $f$ (as well as
on $\NN$, $\alpha$, $s$ and $Y_s$), and is invalid in general with
$\NN$ independent of $f$.
\item[Case ``$-$" (``strongly negative
case"):]  Implication~\eqref{eq0-alpha} is invalid  in general even
with a constant $\NN^*$ that depends on the function $f$.
\end{description}

\begin{remark}
We  emphasize that, in all cases ``$+$'' of this type discussed in
this survey, inequality \ineq{eq0-alpha} holds with $\NN^*=\NN$.
\end{remark}

We finally mention that all ``Statements'' throughout this survey
are valid for some values of parameters and invalid for some other
values. The above cases/symbols ``$+$'', ``$\oplus$", ``$\ominus$"
and ``$-$" will be used to describe various cases of their validity
(expressions of type ``the quadruple $(k,r,q, s)$ is strongly
positive in Statement~$X$'', ``Statement~$X$ is strongly positive
for the quadruple $(k,r,q, s)$", etc. all have the same meaning and
will be used interchangeably).

We use the symbol ``$?$'' to indicate that we do not know at all
which of the symbols ``$+$'', ``$\oplus$" (if applicable),
``$\ominus$" or ``$-$" should be put in its place (\ie if a problem
is completely open). The symbol ``$?^*$'' is used if we know
something about this case, but the problem is still not completely
resolved.

\section{Historical background}\label{S2}

Perhaps the first investigation of SPA was done in 1873 by
Chebyshev, who constructed an algebraic polynomial having the
minimum uniform norm on $[-1,1]$ among all nondecreasing polynomials
of the form $\epsilon x^n + a_1 x^{n-1} + \dots + a_n$ with
$\epsilon = 1$ or $-1$ (see  \cite{Ch1873} or \cite{C1955}). Namely,
Chebyshev \cite{Ch1873} showed  that
\begin{eqnarray*}
\lefteqn{
\inf \left\{ \norm{P_n}{} \st P_n(x) = x^n + Q_{n}(x)\,, \;
Q_n\in\P_n \; \mbox{\rm and } P_n\in\Delta^{(1)} \right\} }\\
& & \mbox{}  = \left\{
\begin{array}{ll} \ds
 2\left( \frac{m!}{(2m-1)!!}\right)^2  \,, & \mbox{\rm if $n=2m$}, \\
\ds \left( \frac{m!}{(2m-1)!!} \right)^2 \,, & \mbox{\rm if
$n=2m+1$},
\end{array}
\right.
\end{eqnarray*}
and
\begin{eqnarray*}
\lefteqn{ \inf \left\{ \norm{P_n}{} \st P_n(x) = x^n + Q_{n}(x)\,, \; Q_n\in\P_n  \; \mbox{\rm and }
(-P_n)\in\Delta^{(1)} \right\} }\\
& & \mbox{}   = \left\{
\begin{array}{ll} \ds
2\left( \frac{m!}{(2m-1)!!}\right)^2   \,, & \mbox{\rm if $n=2m$}, \\
\ds \left(1+ \frac 1m \right)\left( \frac{m!}{(2m-1)!!} \right)^2
\,, & \mbox{\rm if $n=2m+1$}.
\end{array}
\right.
\end{eqnarray*}

In 1927, Bernstein \cite{B1927} (see also \cite{collected1}*{pp. 339--349}) obtained several analogous results for
multiply monotone polynomials (smooth functions are called
``multiply monotone of order $\mu$'' if their first $\mu$
derivatives are nonnegative; sometimes these functions are referred
to as ``absolutely monotone of order $\mu$'').

It is rather well known by now that, if $f\in\Delta^{(q)}$, then its
Bernstein polynomial  (introduced by Bernstein in 1912, see
\cite{Bern1912} and \cite{collected1}*{pp. 105-106})
$$
B_n(f,x)=\frac1{2^n}\sum_{j=0}^n\binom
njf\left(\frac{n-2j}n\right)(1+x)^{n-j}(1-x)^j
$$
is also in $\Delta^{(q)}$. It is not clear who was the first to
notice this shape preserving property of Bernstein polynomials.
Popoviciu knew about it as early as in 1934 (see \cite{p1934}),
but it is not clear if Bernstein himself was aware of this property
before then. Since Bernstein polynomials associated with $f\in
\C[-1,1]$ uniformly approximate $f$, it follows that the Weierstrass
approximation theorem is valid for SPA, that is, for every
$f\in\Delta^{(q)}$,
$$
E_n^{(q)}(f)\to0,\quad n\to\infty.
$$

In 1965, Shisha\hboxx{\cite{S1965}} proved that, for $f\in\C^r \cap
\Delta^q$, $1\leq q\leq r$,
\begin{equation} \label{eqshisha}
E_n^{(q)}(f) \leq c (q, r) \frac{1}{n^{r-q}} \omega(f^{(r)}, 1/n) .
\end{equation}
The proof was based on the rather obvious observation that, for
$f\in\C^q[-1,1]\cap\Delta^q$, we have
\begin{equation}
 E_n^{(q)}(f)  \leq  c(q) E_{n-q} (f^{(q)})  . \label{eq2}
\end{equation}
Indeed, let $Q_{n-q}\in\P_{n-q}$ be such that $E_{n-q} (f^{(q)})  =
\|f^{(q)}-Q_{n-q}\|$, and $P_n\in\P_{n}$ be such that
$P_n^{(\nu)}(0) = f^{(\nu)}(0)$, $0\leq \nu\leq q-1$, and
$P_n^{(q)}(x) := Q_{n-q}(x)+ E_{n-q} (f^{(q)})$. Hence,
\[
\left\| f^{(q)} -P_n^{(q)} \right\| \leq 2 E_{n-q} (f^{(q)}) ,
\]
 and
$P_n^{(q)}(x) \geq f^{(q)}(x) \geq 0$, \ie $P_n\in\Delta^q$.
Finally,
\begin{eqnarray*}
 E_{n}^{(q)}(f)  & \leq &  \|f-P_n\|  = \left\| \int_{0}^x \int_{0}^{t_1}
\cdots\int_{0}^{t_{q-1}} \left( f^{(q)}(t_q)-P_n^{(q)}(t_q)
\right)\dd t_q
\cdots \dd t_1 \right\| \\
 & \leq &   \frac{1}{q!}  \, \left\| f^{(q)} -P_n^{(q)} \right\|
\leq  \frac{2}{q!} E_{n-q} (f^{(q)})  .
\end{eqnarray*}
While inequality \ineq{eqshisha} differs from ``the optimal estimate
of this type''  by the factor of $n^q$, it was perhaps the first
attempt to  obtain a nontrivial direct estimate for SPA and brought
attention to this area.

The first major developments in this area appeared in papers by
Lorentz and Zeller and by DeVore. Lorentz and
Zeller\hboxx{\cite{LZ1969}} constructed, for each $q\ge1$, a
function $f\in\Delta^{(q)}\cap \C^{q}[-1,1]$, such that
$$
\limsup_{n\to\infty}\frac{E_n^{(q)}(f)}{E_n(f)}=\infty.
$$
This means  that questions on SPA do not trivially reduce to those
on unconstrained approximation. Yet, Lorentz and
Zeller\hboxx{\cite{LZ1968}} (for $r=0$), Lorentz\hboxx{\cite{L1972}}
($r=1$) and DeVore\hboxx{\cite{DV1977}}, proved the exact analogue
of the Jackson type estimate, namely that, for each function
$f\in\Delta^{(1)}\cap \C^{r}[-1,1]$,
$$
E_n^{(1)}(f)\le \frac
{c(r)}{n^r}\omega\bigl(f^{(r)},n^{-1}\bigr),\quad n\ge r.
$$
Furthermore, for each function $f\in\Delta^{(1)}$,
DeVore\hboxx{\cite{DV1976}} proved the estimate for the second
modulus of smoothness,
$$
E_n^{(1)}(f)\le c\omega_2\left(f,n^{-1}\right),\quad n\ge 2.
$$
It is impossible to replace $\omega_2$ by $\omega_3$ in the above
estimate due to a negative result by Shvedov
 (see \cite{S1981}), who proved that, for each $A>0$ and $n\in\N$, there is a function $f=f_{n,A}\in\Delta^{(q)}$
such that
\begin{equation}\label{q+2}
E_n^{(q)}(f)\ge A\omega_{q+2}\left(f,1\right).
\end{equation}
Newman \cite{N1979} obtained the first ``optimal'' estimate for
comonotone approximation (earlier results on comonotone
approximation are due to Newman, Passow and Raymon \cite{NPR1972},
Passow, Raymon and Roulier \cite{PRR1974}, Passow and Raymon
\cite{PR1974},
 and Iliev \cite{I1978}).
He showed that, if $f\in\Delta^{(1)}(Y_s)$, then
$$
E_n^{(1)}(f,Y_s)\le c(s)\omega \left(f, n^{-1} \right),\quad n\ge1.
$$
 Shvedov \cite{S1981-co}   proved that, if
$f\in\Delta^{(1)}(Y_s)$, then
$$
E_n^{(1)}(f,Y_s)\le c(s)\omega_2\left(f, n^{-1} \right),\quad n\ge
\NN,
$$
where $\NN=\NN(Y_s)$, and that this estimate is no longer valid with
$\NN$ independent of $Y_s$.

We also mention the following pre-1980 papers which are somewhat
related to the topics discussed in this survey: Roulier
\cites{R1968, R1971,R1973,R1975,R1976}, Lorentz and Zeller
\cite{LZ1970}, Lim \cite{lim}, R. Lorentz \cite{LR1971}, Zeller
\cite{Z1973}, DeVore \cite{DV1974}, Popov and Sendov \cite{PS1974},
Gehner \cite{gehner}, Kimchi and Leviatan \cite{KL1976}, Passow and
Roulier \cite{PRo1976}, Ishisaki \cite{I1977}, Iliev
\cite{I1978-serdica}, and Myers and Raymon \cite{MR1978}.

\section{Unconstrained polynomial approximation}\label{S3}

In this section, we remind the reader of the direct and inverse
theorems for polynomial approximation. In particular, we emphasize
again that, in order to get matching direct and inverse estimates,
one should use (i) pointwise estimates involving the usual moduli of
smoothness $\omega_k\left(f^{(r)},\rho_n(x)\right)$, or (ii) uniform
estimates in terms of generalized Ditzian-Totik moduli of smoothness
$\omega^\varphi_{k,r}(f^{(r)},n^{-1})$. Jackson-Stechkin type
estimates (see \cor{jscor} below) are NOT optimal in this sense.

\subsection{Nikolskii type pointwise estimates}

In 1946, Nikolskii \cite{N1946} showed that, for any function $f$
such that $\omega(f,\delta)\leq \delta$, it is possible to construct
a sequence of polynomials $P_n\in\P_n$ such that, for all $x\in
[-1,1]$,
\[
|f(x)-P_n(x)|\leq\frac{\pi}{2}\cdot\frac{\sqrt{1-x^2}}{n}+|x|
O\left(\frac{\ln n}{n^2}\right) .
\]
Towards the end of the 1960's, the constructive theory of
approximation of functions by algebraic polynomials was completed.

Timan \cite{T1951}  (for $k=1$), Dzjadyk \cite{Dz1958} and,
independently, Freud \cite{F1959}  (for $k=2$), and Brudnyi
\cite{Br1963}  (for $k\ge 3$) proved the following direct theorem
for approximation by algebraic polynomials involving Nikolskii type
pointwise estimates.

\begin{theorem}[Direct theorem: pointwise estimates] \label{directth}
Let $k\in\N$ and $r\in\N_0$. If $f\in \C^r[-1,1]$, then for each
$n\ge k+r$ there is a polynomial $P_n\in\P_n$, such that for every
$x\in[-1,1]$,
$$
|f(x)-P_n(x)|\le
c(k,r)\rho_n^r(x)\omega_k\bigl(f^{(r)},\rho_n(x)\bigr).
$$
\end{theorem}

This pointwise estimate readily implies the classical
Jackson-Stechkin type estimate.

\begin{corollary}[Jackson-Stechkin type estimates] \label{jscor}
  If $f\in \C^r[-1,1]$, then
$$
E_n(f)\le \frac{c(k,r)}{n^r}\omega_k\bigl(f^{(r)},n^{-1}\bigr),
\quad n\ge k+r.
$$
\end{corollary}

Trigub \cite{Trigub62} ($k=1$) and Gopengauz \cite{gop} ($k\geq 1$)
proved the generalization of \thm{directth} for the simultaneous
approximation of a function and its derivatives by a polynomial and
its corresponding derivatives. (For $k=1$ and $1/n$ instead of
$\rho_n(x)$, this result was proved by Gelfond \cite{gel}.)

\begin{theorem}[Simultaneous approximation of a function and its derivatives:
pointwise estimates]  \label{trigthm} \mbox{} Let $k\in\N$ and
$r\in\N_0$. If $f\in \C^r[-1,1]$, then for each $n\ge k+r$ there is
a polynomial $P_n\in\P_n$ satisfying, for $0\leq \nu\leq r$ and
$x\in[-1,1]$,
$$
|f^{(\nu)}(x)-P_n^{(\nu)}(x)|\le
c(k,r)\rho_n^{r-\nu}(x)\omega_k\bigl(f^{(r)},\rho_n(x)\bigr).
$$
\end{theorem}

The following stronger result on simultaneous polynomial
approximation, \thm{kthm}, first appeared in \cite{kop-sim}. We also
note that, while \thm{kthm} was not stated in \cite{M1992}, it can
be proved similarly to \cite{M1992}*{Theorem 15.3} using \cite{M1992}*{Lemmas 15.3 and 4.2$'$}.

\begin{theorem}[Simultaneous approximation of a function and its derivatives:
pointwise estimates] \label{kthm} \mbox{} Let $k\in\N$ and
$r\in\N_0$. If $f\in \C^r[-1,1]$, then for each $n\ge k+r$ there is
a polynomial $P_n\in\P_n$ satisfying, for $0\leq \nu\leq r$ and
$x\in[-1,1]$,
$$
|f^{(\nu)}(x)-P_n^{(\nu)}(x)|\le c(k,r)
\omega_{k+r-\nu}\bigl(f^{(\nu)},\rho_n(x)\bigr).
$$
\end{theorem}

\begin{remark} One of the consequences of \thm{kthm} is that if $k,q\in\N$,  $f\in\C^q[-1,1]$  and $f^{(q)}$ is
{\em strictly} positive on $[-1,1]$, then, for sufficiently large
$n$ (depending on $k$, $q$ and $f$), there exists a polynomial
$P_n\in\P_n$  with a positive $q$-th derivative on $[-1,1]$ (\ie
$P_n\in\Delta^{(q)}$) such that
\[
|f (x)-P_n (x)|\le c(k,q) \omega_k\bigl(f ,\rho_n(x)\bigr).
\]
\end{remark}

We call a  function $\upphi$ a $k$-majorant (and write
$\upphi\in\Phi^k$) if it satisfies the following conditions:
\begin{enumerate}[(i)]
\item $\upphi \in\C[0,\infty)$,
$\upphi(0)=0$,
\item $\upphi$ is nondecreasing on $(0,\infty)$,
\item $x^{-k} \upphi(x)$ is nonincreasing on $(0,\infty)$.
\end{enumerate}

\begin{theorem}[Inverse theorem for pointwise estimates: Dzjadyk \cite{D1956},
Timan \cite{T1957}, Lebed' \cite{L1957}, Brudnyi \cite{B1959}]
\mbox{} \label{inverse} Let $k\in\N$, $r\in\N_0$, $\upphi\in\Phi^k$,
and let $f$ be a given function. If  for every $n\geq k+r$ there
exists $P_n\in\P_n$ such that
\[
|f(x)-P_n(x)| \leq \rho_n^r(x)\upphi\left(\rho_n(x)\right), \quad
x\in [-1,1] ,
\]
then
\[
\omega_k(f^{(r)},\delta)\leq c(k,r)\left(\int_0^\delta r u^{-1}
\upphi(u)\dd u+\delta^k \int_\delta^1 u^{-k-1}\upphi(u)\dd
u\right),\quad 0\leq \delta \leq 1/2 .
\]
In particular, if $\displaystyle \int_0^1ru^{-1}\upphi(u)\dd u <
\infty$, then $f\in\C^r[-1,1]$.
\end{theorem}

Theorems~\ref{directth} and \ref{inverse} (with $\upphi(u) :=
u^\alpha$) imply the following result.

\begin{corollary}[Constructive characterization: pointwise estimates]
Let $k\in\N$, $r\in\N_0$, $0<\alpha <k$, and let $f$ be a given
function. Then, for every $n\geq k+r$ there exist
$P_n\in\P_n$  such that
\[
|f(x)-P_n(x)| \leq c(k,r,\alpha) \rho_n^{r+\alpha}(x) , \quad x\in
[-1,1] ,
\]
if and only if  $f\in\C^r[-1,1]$ and
\[
\omega_k(f^{(r)},\delta)\leq c(k,r,\alpha) \delta^\alpha,\quad 0\leq
\delta \leq 1/2 .
\]
\end{corollary}

\subsection{Ditzian-Totik type estimates}

As discussed above, Jackson-Stechkin type estimates are rather weak
in the sense that they do not provide a constructive
characterization of classes of functions having prescribed order of
approximation by algebraic polynomials. One has to measure
smoothness taking into account the distance from the endpoints of
$[-1,1]$ in order to get this characterization.

\begin{theorem}[Direct theorem: Ditzian and Totik \cite{DT}] \label{dt-thm}
If $f\in \C[-1,1]$ and $k\in\N$, then for each $n\ge k$,
\[
E_n(f) \leq c(k)\omega_k^\varphi \left(f, n^{-1}\right) ,
\]
where $\omega_k^\varphi \left(f, \delta \right) = \sup_{0\le h\le
\delta} \|\Delta_{h\varphi(\cdot)}^k(f,\cdot)\|$  and
$\varphi(x):=\sqrt{1-x^2}$.
\end{theorem}

The  generalized Ditzian-Totik moduli
$\omega^\varphi_{k,r}(f^{(r)},\delta)$  discussed below in this
section were introduced  in \cite{M1992}.

For $r\ge1$, we say that $f\in \B^r$ if $f^{(r-1)}$ is locally
absolutely continuous in $(-1,1)$ and $\varphi^r f^{(r)}\in
\LL_\infty[-1,1]$. Also, let $\C^0_\varphi:=\C[-1,1]$ and
\[
\C^r_\varphi = \left\{ f\in \C^{r}(-1,1) \st
\lim_{x\to\pm1}\varphi^r(x)f^{(r)}(x)=0 \right\}, \quad r\geq 1 .
\]
For $f\in \C^r_\varphi$, $r\geq 0$, we denote
$$
\omega^\varphi_{k,r}(f^{(r)},\delta):=\sup_{0\le h\le
\delta}\sup_{x:|x|+\frac{kh}2\varphi(x)<1} K^r
\left(x,\frac{kh}2\right)|\Delta_{h\varphi(x)}^k(f^{(r)},x)|,
$$
where $K(x,\mu):=\varphi(|x|+\mu\varphi(x))$.

Note that if $r=0$, then
$$
\omega^\varphi_{k,0}(f,\delta)\equiv\omega_k^\varphi(f,\delta),
$$
where $\omega_k^\varphi(f,\delta)$ is the $k$th Ditzian-Totik
modulus of smoothness defined above.

Clearly $\C^r_\varphi\subset\B^r$, while it is known (see, \eg\,
\cite{s-book}*{Chapter 3.10}) that if $f\in\B^r$, then
$f\in\C^l_\varphi$ for all $0\leq l<r$, and
\begin{equation} \label{uniform2}
\omega_{r-l,l}^\varphi(f^{(l)},\delta)\leq c
\delta^{r-l}\left\|\varphi^r f^{(r)}\right\|_{\LL_\infty[-1,1]}
,\quad \delta>0.
 \end{equation}
Note that if $f\in\C^r_\varphi$, then the following inequality holds
for all $0\leq l\leq r$ and $k\geq 1$ (see \cite{s-book}*{Chapter 3.10}):
\begin{equation} \label{uniform1}
\omega_{k+r-l,l}^\varphi (f^{(l)},\delta) \,\leq\, c \,
\delta^{r-l}\omega_{k,r}^\varphi(f^{(r)},\delta),\quad \delta>0.
\end{equation}

By virtue of \ineq{uniform1}, \thm{dt-thm} immediately implies the
following.

\begin{corollary}[Direct theorem] \label{dirdt}
If $k\in\N$, $r\in\N_0$ and $f\in\C^r_\varphi$, then for each $n\ge
k+r$,
\[
E_n(f) \leq c(k,r)n^{-r}\omega_{k,r}^\varphi\bigl(f^{(r)},
n^{-1}\bigr).
\]
\end{corollary}

A matching inverse theorem is the following generalization of
\cite{DT}*{Theorem 7.2.4} in the case $p=\infty$ (see  \cite{M1992}
or \cite{KLS2010}*{Theorem 3.2}).

\begin{theorem}[Inverse theorem for uniform estimates] \label{Theorem 3.199}
Let $k\in\N$, $r\in\N_0$, $\NN\in\N$, and  let $\upphi:
[0,\infty)\rightarrow[0,\infty)$ be a nondecreasing function such
that $\upphi(0+)=0$ and
$$
\int_0^1\frac{r\upphi(u)}{u^{r+1}} \dd u<+\infty.
$$
If
$$
E_n(f)\le\upphi\left(n^{-1}\right),\quad \mbox{for all} \quad n\ge
\NN,
$$
then $f\in\C_{\varphi}^r$, and, for any $0\leq \delta \leq 1/2$,
\[
\omega_{k,r}^\varphi(f^{(r)},\delta) \le
c(k,r)\int_0^\delta\frac{r\upphi(u)}{u^{r+1}}\dd u+c(k,r)\delta^k
\int_\delta^1\frac{\upphi(u)}{u^{k+r+1}}\dd u
 + c(k,r,\NN)\delta^kE_{k+r}(f).
\]
If, in addition, $\NN\le k+r$, then the following Bari--Stechkin
type estimate holds:
$$
\omega_{k,r}^\varphi(f^{(r)},\delta)\le
c(k,r)\int_0^\delta\frac{r\upphi(u)}{u^{r+1}}\dd
u+c(k,r)\delta^k\int_\delta^1 \frac{\upphi(u)}{u^{k+r+1}}\dd u,
\quad \delta\in\left[0,1/2\right].
$$
\end{theorem}

If $\upphi(u) := u^{\alpha}$ we get the following corollary (see
\cite{KLS2010}*{Theorem 3.3}).

\begin{corollary} \label{cor199}
Let $k\in\N$, $r\in\N_0$  and $\alpha >0$, be such that $r<\alpha<
k+r$.
If
$$
n^\alpha E_n(f)\le 1,\quad \mbox{for all} \quad n \ge \NN,
$$
where $\NN \geq k+r$, then $f\in \C_{\varphi}^r$ and
$$
\omega_{k,r}^\varphi(f^{(r)},\delta)\le c(\alpha,k,r)
\delta^{\alpha-r} + c(\NN,k, r) \delta^{k} E_{k+r}(f) \,.
$$
In particular, if $\NN=k+r$, then
$$
\omega_{k,r}^\varphi(f^{(r)},\delta)\le c(\alpha,k,r)
\delta^{\alpha-r}\,.
$$
\end{corollary}

Corollaries~\ref{dirdt} and \ref{cor199} imply the following result.

\begin{corollary}[Constructive characterization: uniform estimates]
Let $k\in\N$, $r\in\N_0$  and $\alpha >0$, be such that $r<\alpha<
k+r$. Then,
$$
 E_n(f)\le c(k,r,\alpha) n^{-\alpha},\quad \mbox{for all} \quad n \ge k+r,
$$
if and only if $f\in \C_{\varphi}^r$ and
$$
\omega_{k,r}^\varphi(f^{(r)},\delta)\le c(k,r,\alpha)
\delta^{\alpha-r}\,.
$$
\end{corollary}

\section{Jackson-Stechkin type   estimates for $q$-monotone approximation, $q\geq 1$}\label{S4}

\subsection{$q$-monotone approximation of functions from $\C[-1,1]$ ($r=0$)} \label{sub7.1}

In this section, we discuss the validity  of the following statement
(note that this is the case $r=0$ of the more general \prop{prop1}
discussed in Section~\ref{sub7.2}).

\begin{proposition} \label{prop0}
Let $q\in\N$, $k\in\N$ and $\NN \in\N$. If $f\in \Delta^{(q)}\cap \C
[-1,1]$, then
\begin{equation} \label{eq44}
E_n^{(q)}(f)\le c(k, q) \omega_k\left(f , n^{-1}\right), \quad n\ge
\NN .
\end{equation}
\end{proposition}

\subsubsection{Case ``$+$''}

  Case ``$+$'' is the case when \ineq{eq44} holds with $\NN=k+r$.

For $(k=1,q=1)$, the estimate \ineq{eq44} with $\NN=1$ was proved by
Lorentz and Zeller\hboxx{\cite{LZ1968}}. Beatson \cite{Be1978}
proved \ineq{eq44} for $k=1$ and all $q$. For $(k=2,q=1)$,
\ineq{eq44} with $\NN=2$ was established by DeVore \cite{DV1976}.
Shevdov \cite{S1980} later extended it to $k=2$ and all $q\ge1$.

So, for the first and second moduli of smoothness, we have exactly
the same estimate as in the unconstrained approximation. Thus the
cases $(k=1,q\in\N)$ and  $(k=2,q\in\N)$ are of type ``$+$''. Only
two other cases of type ``$+$'' are known. Namely, Hu, Leviatan and
Yu \cite{HLY1994}  and, independently, Kopotun \cite{K1994} proved
\ineq{eq44} in the case $(k=3,q=2)$, with $\NN=3$, and Bondarenko
\cite{Bo2002}  proved \ineq{eq44} for $(k=3, q=3)$, with $\NN=3$.

\subsubsection{Case ``$-$''}

Wu and Zhou \cite{WZ1992}  (for $k \geq q+3$ and $q\ge1$), and
Bondarenko and Prymak \cite{BP2004} (for $k\ge3$ and $q\ge 4$),
proved that there is a function $f\in\Delta^{(q)}$, such that
\begin{equation} \label{eq5}
\limsup_{n\to\infty}\frac{E_n^{(q)}(f)}{\omega_k\left(f,1/n\right)}=\infty.
\end{equation}
In other words, in the cases $(k \geq q+3,  q\geq 1)$ and $(k\ge3,
q\ge 4)$, estimate \ineq{eq44} is not valid in general even if $\NN$
is allowed to depend on $f$.

\subsubsection{Case ``$\ominus$''}

Shvedov's negative result (see \cite{S1981}) implies that, in the
cases $(k=3, q=1)$,   $(k=4, q=2)$, and $(k=5, q=3)$, for each $A>0$
and $n\in\N$ there is a function $f=f_{n,A}\in\Delta^{(q)}$ such
that
\begin{equation} \label{shvneg}
E_n^{(q)}(f)\ge A\omega_k\left(f,1\right).
\end{equation}
At the same time,  positive results for the first two of these cases
were proved in \cite{LS1998}   and
 \cite{LS2003}. Namely, it was shown that, in the cases
$(k=3, q=1)$ and $(k=4, q=2)$, if $f\in\Delta^{(q)}$, then
\ineq{eq44} holds with a constant $\NN$ that depends on the function $f$.
We emphasize that \ineq{shvneg} implies that, in these cases,
\ineq{eq44} cannot be valid with the constant $\NN$  independent of
$f$.

We refer to all such cases
as ``$\ominus$''. Namely, we denote by ``$\ominus$'' the cases when
\ineq{eq44} holds with $\NN$ that depends on $f$ (this means that
there is no $f$ for which \ineq{eq5} is valid) and does not hold in
general with $\NN$ independent of $f$.

All of the above cases are conveniently summarized in \newfig{fig1}.

\begin{figure}[H]
$$
\begin{matrix}
 q &\vdots  &\vdots&\vdots&\vdots&\vdots&\vdots&\vdots&\updots
      \cr 5&+&+  &-  &- &- &-& -&\cdots&
      \cr 4&+ &+   &-  &- &- &-& -&\cdots&
      \cr 3 &+&+ &+  &?  &?^* &-&-&\cdots&
      \cr 2 &+& +&+&\ominus &-  &-&-& \cdots&
      \cr 1&+ &+&\ominus&- &- &-&-&\cdots&
      \cr &1&2 &3&4&5&6&7&k&
\end{matrix}
$$
\caption{$q$-monotone approx. of functions from $\C[-1,1]$
($r=0$), validity of $E_n^{(q)}(f)\le c(k, q) \omega_k\left(f ,
n^{-1}\right)$, $n\ge \NN$ }\label{fig1}
\end{figure}

\begin{remark}
 It follows from \ineq{shvneg}
that  ``$?^*$'' in \newfig{fig1} cannot be replaced by ``$+$".
\end{remark}

We therefore have two problems related to \ineq{eq44}.
\begin{problem} \label{op1}
Does there exist a function $f\in \C[-1,1]$ with a convex derivative
on $(-1,1)$, that is, $f\in\Delta^{(3)}$, such that, for each
sequence $\{P_n\}_{n=1}^\infty \subset \P_n$ of algebraic
polynomials satisfying
$$
P_n^{(3)}(x)\ge0,
$$
we have
$$
\limsup_{n\to\infty}\frac{\|f-P_n\|_{\C[-1,1]}}{\omega_5\left(f,1/n
\right)}=\infty ?
$$
In other words, is the case $(k=5, q=3)$  strongly negative
$(``-")$?
\end{problem}

\begin{problem} What can be said if $\omega_5$ in {\rm\op{op1}} is replaced by $\omega_4$?
\end{problem}

\subsection{$q$-monotone approximation of functions from $\C^r[-1,1]$ and $\W^r$, $r\geq 0$} \label{sub7.2}

In this section, we discuss the validity  of the following
\prop{prop1}, and ask for which triples $(k,r,q)$ this statement is
valid and for which it is invalid (note that the case $r=0$ was
already considered in \prop{prop0} in Section~\ref{sub7.1}).

\begin{proposition} \label{prop1}
If $k\in\N$, $r\in\N_0$, $q\in\N$, $\NN\in\N$ and
$f\in\Delta^{(q)}\cap \C^r[-1,1]$, then
\begin{equation} \label{eq4}
E_n^{(q)}(f)\le\frac{c(k,r,q)}{n^r} \omega_k\bigl(f^{(r)},
n^{-1}\bigr), \quad n\ge \NN.
\end{equation}
\end{proposition}

Recall that an estimate similar to \ineq{eq4} is valid with
$\NN=k+r$ in the unconstrained case. Also we remind the reader that
we say that, for the triple $(k,r,q)$,  \prop{prop1} is
\begin{itemize}
\item ``strongly positive''
(``$+$'')  if \ineq{eq4} holds with $\NN=k+r$,
\item   ``weakly negative'' (``$\ominus$'') if \ineq{eq4} holds with
$\NN=\NN(f)$ and is not valid with $\NN$ independent of $f$,
\item ``strongly negative'' (``$-$'')  if \ineq{eq4} is
not valid at all, that is, there is a function $f\in\Delta^{(q)}\cap
\C^r[-1,1]$ such that
\begin{equation} \label{eq6}
\limsup_{n\to\infty}\frac{n^rE_n^{(q)}(f)}{\omega_{k}\left(f^{(r)},
n^{-1}\right)}
=\infty.
\end{equation}
\end{itemize}

For completeness, we also consider the case $k=0$ in \ineq{eq4}
requiring that $f$ belong to   $\W^r$, the space of $(r-1)$ times
continuously differentiable functions $f$ in $[-1,1]$ such that
$f^{(r-1)}$ is absolutely continuous in $(-1,1)$ and
$\|f^{(r)}\|_{\LL_\infty[-1,1]} < \infty$.

\begin{proposition}[$k=0$] \label{prop2}
If  $r\in\N$, $q\in\N$, $\NN\in\N$, and $f\in\Delta^{(q)}\cap \W^r$,
then
$$
E_n^{(q)}(f)\le
\frac{c(r,q)}{n^r}\|f^{(r)}\|_{\LL_\infty[-1,1]},\quad n\ge \NN.
$$
\end{proposition}

The following are the ``truth tables''  for \props{prop1}{prop2}.

\begin{figure}[H]
$$
\begin{matrix}
 r &\vdots  &\vdots&\vdots&\vdots&\vdots&\vdots&\updots
      \cr 3 &+&+ &+  &+  &+ &+ &\cdots&
      \cr 2 &+&+ &+  &+  &+ &+ &\cdots&
      \cr 1 &+& +&+ &+ &+  &+ &\cdots&
      \cr 0& &+&+&\ominus &- &-&\cdots&
      \cr &0&1 &2&3&4&5&k&
\end{matrix}
$$
\caption{ Monotone approx. ($q=1$), validity of $E_n^{(1)}(f)\le
c(k,r) n^{-r} \omega_k\left(f^{(r)}, n^{-1}\right)$, $n\ge \NN$ }
\label{fig2}
\end{figure}

These results appeared in the papers by Lorentz and Zeller
\cite{LZ1968}, Lorentz \cite{L1972}, DeVore \cite{DV1976}  and
\cite{DV1977}, Shvedov \cite{S1979}, Wu and Zhou \cite{WZ1992},
Shevchuk \cite{S1992} (see also \cite{S-doklady}), and Leviatan and
Shevchuk \cite{LS1998}.

It is convenient to summarize the above references in another table.
Note that, in the case ``$\ominus$'', we first put a reference to
the negative result, then to the positive one.

\samefigure
\begin{figure}[H]
{\scriptsize
$$
\begin{matrix}
 r &\vdots  &\vdots&\vdots&\vdots&\vdots&\vdots&\updots
      \cr 3 & \hboxx{\cite{DV1977}}&\hboxx{\cite{DV1977}}&\hboxx{\cite{S1992}}
      &\hboxx{\cite{S1992}}&\hboxx{\cite{S1992}}&\hboxx{\cite{S1992}}&\cdots&
      \cr 2 &\hboxx{\cite{L1972}}&\hboxx{\cite{DV1977}}&\hboxx{\cite{S1992}}
      &\hboxx{\cite{S1992}}&\hboxx{\cite{S1992}}&\hboxx{\cite{S1992}}&\cdots&
      \cr 1 &\hboxx{\cite{LZ1968}}&\hboxx{\cite{L1972}}&\hboxx{\cite{S1992}}
      &\hboxx{\cite{S1992}}&\hboxx{\cite{S1992}}&\hboxx{\cite{S1992}}&\cdots&
      \cr 0& &\hboxx{\cite{LZ1968}}&\hboxx{\cite{DV1976}}&\hboxx{\cite{S1979}},
      \hboxx{\cite{LS1998}}&\hboxx{\cite{WZ1992}}&\hboxx{\cite{WZ1992}}&\cdots&
      \cr &0&1 &2&3&4&5&k&
\end{matrix}
$$
}
\caption{References for \newfig{fig2}} \label{fig2a}
\end{figure}
\esamefigure

\begin{figure}[H]
$$
\begin{matrix}
 r &\vdots  &\vdots&\vdots&\vdots&\vdots&\vdots&\vdots&\updots
      \cr 3 &+&+   &+  &+ &+ &+& +& \cdots&
      \cr 2 &+&+ &+  &+  &+ &+ & +&\cdots&
      \cr 1 &+& +&+ &\ominus &-  & - &  -&\cdots&
      \cr 0& &+&+&+ &\ominus &-& - &\cdots&
      \cr &0&1 &2&3&4&5&6&k&
\end{matrix}
$$
\caption{ Convex approx. ($q=2$), validity of $E_n^{(2)}(f)\le
c(k,r) n^{-r} \omega_k\left(f^{(r)}, n^{-1}\right)$, $n\ge \NN$
}\label{fig3}
\end{figure}

These results appeared in the papers by  Beatson \cite{Be1978},
Shvedov \cite{S1981}, Wu and Zhou \cite{WZ1992}, Mania (see in
\cite{M1992}*{Theorems 17.2 and 16.1}), Hu, Leviatan and Yu
\cite{HLY1994},  Kopotun \cite{K1994}, Nissim and Yushchenko
\cite{NY2003}, and Leviatan and Shevchuk \cite{LS2003}.

\samefigure
\begin{figure}[H]
{\scriptsize
$$
\begin{matrix}
 r &\vdots  &\vdots&\vdots&\vdots&\vdots&\vdots&\vdots&\updots
      \cr 3 &\hboxx{\cite{M1992}}&\hboxx{\cite{M1992}}&\hboxx{\cite{M1992}}&\hboxx{\cite{M1992}}
      &\hboxx{\cite{M1992}}&\hboxx{\cite{M1992}}& \hboxx{\cite{M1992}}& \cdots&
      \cr 2 &\hboxx{\cite{S1981}}&\hboxx{\cite{M1992}}&\hboxx{\cite{M1992}}&\hboxx{\cite{M1992}}
      &\hboxx{\cite{M1992}}&\hboxx{\cite{M1992}}&\hboxx{\cite{M1992}}&\cdots&
      \cr 1 &\hboxx{\cite{Be1978}}& \hboxx{\cite{S1981}}& \mbox{\cites{HLY1994, K1994} }
      &\hboxx{\cite{M1992}},\hboxx{\cite{LS2003}}&\hboxx{\cite{NY2003}}&\hboxx{\cite{NY2003}}
      &\hboxx{\cite{NY2003}}&\cdots&
      \cr 0& &\hboxx{\cite{Be1978}}&\hboxx{\cite{S1981}}&  \mbox{\cites{HLY1994, K1994}}
      &\hboxx{\cite{S1981}},\hboxx{\cite{LS2003}}&\hboxx{\cite{WZ1992}}&\hboxx{\cite{WZ1992}}&\cdots&
      \cr &0&1 &2&3&4&5&6&k&
\end{matrix}
$$
}
\caption{References for \newfig{fig3}}\label{fig3a}
\end{figure}
\esamefigure

\begin{remark}
  The cases $(k=3,r=0,q=2)$ and $(k=2,r=1,q=2)$
(both of type ``$+$''), were proved in \hboxx{\cite{HLY1994}} and
\hboxx{\cite{K1994}} simultaneously and independently.
\end{remark}

\begin{figure}[H]
$$
\begin{matrix}
 r &\vdots  &\vdots&\vdots&\vdots&\vdots&\updots
      \cr 4 &-&-  &-  &- &- & \cdots&
      \cr 3 &- &-   &-  &- &- & \cdots&
      \cr 2 &+&- &-  &-  &- &\cdots&
      \cr 1 &+& +&-&- &-  & \cdots&
      \cr 0& &+&+&- &- &\cdots&
      \cr &0&1 &2&3&4&k&
\end{matrix}
$$
\caption{ $q$-monotone approx. ($q\geq 4$), validity of
$E_n^{(q)}(f)\le c(k,r) n^{-r} \omega_k\left(f^{(r)},
n^{-1}\right)$, $n\ge \NN$ }\label{fig4}
\end{figure}

The results appeared in the papers by  Beatson \cite{Be1978},
Shvedov \cite{S1980}, and Bondarenko and
Prymak\hboxx{\cite{BP2004}}.

\samefigure
\begin{figure}[H]
{\scriptsize
$$
\begin{matrix}
 r &\vdots  &\vdots&\vdots&\vdots&\vdots&\updots
      \cr 4 &\hboxx{\cite{BP2004}}&\hboxx{\cite{BP2004}}
      &\hboxx{\cite{BP2004}}&\hboxx{\cite{BP2004}}&\hboxx{\cite{BP2004}}& \cdots&
      \cr 3 &\hboxx{\cite{BP2004}}&\hboxx{\cite{BP2004}}&\hboxx{\cite{BP2004}}
      &\hboxx{\cite{BP2004}}&\hboxx{\cite{BP2004}}& \cdots&
      \cr 2
      &\hboxx{\cites{Be1978, S1980}}
      &\hboxx{\cite{BP2004}}&\hboxx{\cite{BP2004}}
      &\hboxx{\cite{BP2004}}&\hboxx{\cite{BP2004}}&\cdots&
      \cr 1 &\hboxx{\cite{Be1978}} &\hboxx{\cites{Be1978, S1980}}
      &\hboxx{\cite{BP2004}}&\hboxx{\cite{BP2004}}
      &\hboxx{\cite{BP2004}}& \cdots&
      \cr 0& &\hboxx{\cite{Be1978}} & \hboxx{\cite{S1980}}
      &\hboxx{\cite{BP2004}}&\hboxx{\cite{BP2004}}&\cdots&
      \cr &0&1 &2&3&4&k&
\end{matrix}
$$
}
\caption{References for \newfig{fig4}} \label{fig4a}
\end{figure}
\esamefigure

\begin{remark} It is worth mentioning that the
breakthrough in the surprising negative results is due to Konovalov
and Leviatan\hboxx{\cite{KL2003}} who proved that the negative
assertion of ``$\ominus$" is valid for all ``$-$" entries in Table
\ref{fig4}, breaking the Shvedov \cite{S1981} pattern (compare with
\ineq{q+2}). In fact, the negative results in \cite{KL2003} are of
width type and are valid for any increasing sequence of
$n$-dimensional linear subspaces of $\C[-1,1]$.
\end{remark}

Finally, in the case of $3$-monotone approximation, the following is
known.

\begin{figure}[H]
$$
\begin{matrix}
r &\vdots  &\vdots&\vdots&\vdots&\vdots&\vdots&\vdots&\updots
      \cr 5 &?&?   &?  &? &? &?& ?& \cdots&
      \cr 4 &?&?   &?  &? &? &?& ?& \cdots&
      \cr 3 &+& ? &? &?&?&?& ?& \cdots&
      \cr 2 &+& + &? &?^*&-&-& -& \cdots&
      \cr 1 &+& +&+ & ?& ?^* & - &  -& \cdots&
      \cr 0& &+&+&+ & ?&?^*& - &\cdots&
      \cr &0&1 &2&3&4&5&6&k&
\end{matrix}
$$
\caption{ $3$-monotone approx. ($q=3$), validity of $E_n^{(3)}(f)\le
c(k,r) n^{-r} \omega_k\left(f^{(r)}, n^{-1}\right)$, $n\ge \NN$
}\label{fig5}
\end{figure}

\begin{remark}
 It follows from Shvedov\hboxx{\cite{S1981}}
 and Mania (see \cite{M1992}) that ``$?^*$'' in the cases
$(k=5-r,0\leq r\leq 2)$  in \newfig{fig5} cannot be replaced by
``$+$''.
\end{remark}

These results appeared in the papers by  Beatson \cite{Be1978},
Shvedov \cites{S1980, S1981}, Bondarenko \cite{Bo2002}, Mania (see
 \cite{M1992}*{Theorem 16.1}), Nissim and
Yushchenko \cite{NY2003}, and Wu and Zhou \cite{WZ1992}.

\samefigure
\begin{figure}[H]
{\scriptsize
$$
\begin{matrix}
r &\vdots  &\vdots&\vdots&\vdots&\vdots&\vdots&\vdots&\updots
      \cr 5 &?&?   &?  &? &? &?& ?& \cdots&
      \cr 4 &?&?   &?  &? &? &?& ?& \cdots&
      \cr 3 &\hboxx{\cite{Bo2002}}& ? &? &?&?&?& ?& \cdots&
      \cr 2 &\hboxx{\cites{Be1978, S1980}}&\hboxx{\cite{Bo2002}}&? &\hboxx{\cite{M1992}}
      &\hboxx{\cite{NY2003}}&\hboxx{\cite{NY2003}}&\hboxx{\cite{NY2003}}& \cdots&
      \cr 1 &\hboxx{\cite{Be1978}}&\hboxx{\cites{Be1978, S1980}}&\hboxx{\cite{Bo2002}}& ?
      &\hboxx{\cite{M1992}}&\hboxx{\cite{NY2003}}&\hboxx{\cite{NY2003}}& \cdots&
      \cr 0& &\hboxx{\cite{Be1978}}&\hboxx{\cite{S1980}}&\hboxx{\cite{Bo2002}}& ?&\hboxx{\cite{S1981}}
      &\hboxx{\cite{WZ1992}}&\cdots&
      \cr &0&1 &2&3&4&5&6&k&
\end{matrix}
$$
}
\caption{References for \newfig{fig5}} \label{fig5a}
\end{figure}
\esamefigure

Obviously, we would like to replace the question marks in \newfig{fig5}
with definitive answers. We emphasize that it seems to be a hard
problem to replace the question mark in \eword{any} place with a
definitive symbol.

\begin{problem}
Determine which symbol $($among ``$+$'', ``$\ominus$'' and
``$-$''$)$ replaces ``$?$'' or ``$?^*$'' in {\rm\newfig{fig5}} in  any
of the places.
\end{problem}

\section{Nikolskii type pointwise estimates for $q$-monotone approximation, $q\geq 1$}\label{S5}

In this section, we discuss the validity of the following statement
on pointwise estimates for $q$-monotone polynomial approximation.

\begin{proposition} \label{prop4}
If $k\in\N$, $r\in\N_0$, $q\in\N$, $\NN\in\N$, and
$f\in\Delta^{(q)}\cap \C^r[-1,1]$, then there is a sequence
$\{P_n\}_{n=1}^\infty$ of polynomials $P_n\in\Delta^{(q)}\cap \P_n$
such that, for every $n\ge \NN$  and each $x\in[-1,1]$, we have
\begin{equation}
|f(x)-P_n(x)|\le
c(k,r,q)\rho_n^r(x)\omega_k\left(f^{(r)},\rho_n(x)\right). \label{eq7}
\end{equation}
\end{proposition}

Recall that the case ``$+$'' means that \prop{prop4}  is valid with
$\NN=k+r$. On the other hand, note that the case ``$-$'' here means
  (compare with \ineq{eq6}) that there exists a function $f\in\Delta^{(q)}\cap \C^r[-1,1]$ such that,
 for every sequence of polynomials
$P_n\in\P_n\cap\Delta^{(q)}$,
\begin{equation}
\limsup_{n\to\infty}E^{(q)}_{n,r,k}(f):=\limsup_{n\to\infty}
\left\|\frac{f-P_n}{\rho_n^r
\omega_k\left(f^{(r)},\rho_n\right)}\right\|=\infty. \label{eq8}
\end{equation}

For $k=0$, similarly to what was done in Section~\ref{sub7.2}, we
have the following modification of \prop{prop4}.

\begin{proposition}[$k=0$] \label{prop5}
If  $r\in\N$, $q\in\N$, $\NN\in\N$, and $f\in\Delta^{(q)}\cap
\W^{r}$, then there is a sequence $\{P_n\}_{n=1}^\infty$ of
polynomials $P_n\in\Delta^{(q)}\cap \P_n$  such that, for every
$n\ge \NN$  and each $x\in[-1,1]$, we have
\begin{equation}
|f(x)-P_n(x)|\le
c(r,q)\rho_n^r(x)\|f^{(r)}\|_{\LL_\infty[-1,1]}. \label{eq9}
\end{equation}
\end{proposition}

Note that ``$+$'' now means that \ineq{eq9} is true with $\NN=r$,
and   ``$-$''    means that, for every sequence of polynomials
$P_n\in\P_n\cap\Delta^{(q)}$,
\begin{equation}
\limsup_{n\to\infty}E^{(q)}_{n,r,0}(f):=\limsup_{n\to\infty}
\left\|\frac{f-P_n}{\rho_n^r}\right\|=\infty. \label{eq10}
\end{equation}

Since
\begin{equation} \label{newJuly9-2}
\rho_n^r(x)\omega_k\left(f^{(r)},\rho_n(x)\right) \leq c(k,r) n^{-r}
\omega_{k}\left(f^{(r)},n^{-1} \right) ,
\end{equation}
all positive estimates in all tables in this section imply positive
estimates in the corresponding places (in the corresponding tables)
in Section~\ref{sub7.2} (and all negative results in tables in
Section~\ref{sub7.2} imply negative results in the corresponding
places in the corresponding tables in this section).

The following are the ``truth tables'' for \props{prop4}{prop5} for
$q=1$, $q=2$, $q\geq 4$ and $q=3$ (the most difficult case).

\begin{figure}[H]
$$
\begin{matrix}
 r &\vdots  &\vdots&\vdots&\vdots&\vdots&\updots
      \cr 3 &+&+&+&+&+&\cdots&
      \cr 2 &+&+&+&+&+&\cdots&
      \cr 1 &+&+&+&+&+&\cdots&
      \cr 0& &+&+&-&-&\cdots&
      \cr &0&1 &2&3&4&k&
\end{matrix}
$$
\caption{Monotone approx. ($q=1$), validity of $|f(x)-P_n(x)|\le
c(k,r)\rho_n^r(x)\omega_k\left(f^{(r)},\rho_n(x)\right)$  for $x\in
[-1,1]$ and $n\geq \NN$ with $P_n\in\Delta^{(1)}$}\label{fig6}
\end{figure}

These results appeared in the papers by  Lorentz and Zeller
\cite{LZ1968}, DeVore and Yu \cite{DY1985}, Shevchuk \cite{S1992}
(see also \cite{S-doklady}), Wu and Zhou \cite{WZ1992}, and Leviatan
and Shevchuk \cite{LS1998}.

\samefigure
\begin{figure}[H]
{\scriptsize
$$
\begin{matrix}
 r &\vdots  &\vdots&\vdots&\vdots&\vdots&\updots
      \cr 3 &\hboxx{\cite{S1992}}&\hboxx{\cite{S1992}}&\hboxx{\cite{S1992}}
      &\hboxx{\cite{S1992}}&\hboxx{\cite{S1992}}& \cdots&
      \cr 2 &\hboxx{\cite{DY1985}}&\hboxx{\cite{S1992}}&\hboxx{\cite{S1992}}
      &\hboxx{\cite{S1992}}&\hboxx{\cite{S1992}}&\cdots&
      \cr 1 &\hboxx{\cite{LZ1968}}&\hboxx{\cite{DY1985}}&\hboxx{\cite{S1992}}
      &\hboxx{\cite{S1992}}&\hboxx{\cite{S1992}}&\cdots&
      \cr 0& &\hboxx{\cite{LZ1968}}&\hboxx{\cite{DY1985}}&\hboxx{\cite{LS1998}}
      &\hboxx{\cite{WZ1992}}&\cdots&
      \cr &0&1 &2&3&4&k&
\end{matrix}
$$
}
\caption{References for \newfig{fig6}} \label{fig6a}
\end{figure}
\esamefigure

\begin{figure}[H]
$$
\begin{matrix}
 r &\vdots  &\vdots&\vdots&\vdots&\vdots&\vdots&\updots
      \cr 3 &+&+   &+  &+ &+ & +& \cdots&
      \cr 2 &+&+ &+  &+  &+  & +&\cdots&
      \cr 1 &+& +&+ &-  & - &  -&\cdots&
      \cr 0& &+&+&+ &-& - &\cdots&
      \cr &0&1 &2&3&4&5&k&
\end{matrix}
$$
\caption{Convex approx. ($q=2$), validity of $|f(x)-P_n(x)|\le
c(k,r)\rho_n^r(x)\omega_k\left(f^{(r)},\rho_n(x)\right)$ for $x\in
[-1,1]$ and $n\geq \NN$ with $P_n\in\Delta^{(2)}$}\label{fig7}
\end{figure}

The results appear in the papers by Beatson \cite{Be1980}, Leviatan
\cite{Le1986}, Mania and Shevchuk (see in \cite{M1992}*{Theorem 17.2}), Kopotun \cite{K1994}, Wu and Zhou \cite{WZ1992}, and
Yushchenko \cite{Y2002}.

\samefigure
\begin{figure}[H]
{\scriptsize
$$
\begin{matrix}
 r &\vdots  &\vdots&\vdots&\vdots&\vdots&\vdots&\updots
      \cr 3 &\hboxx{\cite{M1992}}&\hboxx{\cite{M1992}}&\hboxx{\cite{M1992}}
      &\hboxx{\cite{M1992}}&\hboxx{\cite{M1992}}&\hboxx{\cite{M1992}}& \cdots&
      \cr 2 &\hboxx{\cite{Le1986}}&\hboxx{\cite{M1992}}&\hboxx{\cite{M1992}}
      &\hboxx{\cite{M1992}}&\hboxx{\cite{M1992}}&\hboxx{\cite{M1992}}&\cdots&
      \cr 1 &\hboxx{\cite{Be1980}}&\hboxx{\cite{Le1986}}&\hboxx{\cite{K1994}}
      &\hboxx{\cite{Y2002}}&\hboxx{\cite{Y2002}}&\hboxx{\cite{Y2002}}&\cdots&
      \cr 0& &\hboxx{\cite{Be1980}}&\hboxx{\cite{Le1986}}&\hboxx{\cite{K1994}}
      &\hboxx{\cite{Y2002}}&\hboxx{\cite{WZ1992}}&\cdots&
      \cr &0&1 &2&3&4&5&k&
\end{matrix}
$$
}
\caption{References for \newfig{fig7}} \label{fig7a}
\end{figure}
\esamefigure

\begin{figure}[H]
$$
\begin{matrix}
 r &\vdots  &\vdots&\vdots&\vdots&\vdots&\updots
      \cr 3 &- &-   &-  &- &- & \cdots&
      \cr 2 &+&- &-  &-  &- &\cdots&
      \cr 1 &+& +&-&- &-  & \cdots&
      \cr 0& &+&+&- &- &\cdots&
      \cr &0&1 &2&3&4&k&
\end{matrix}
$$
\caption{$q$-monotone approx. ($q\geq 4$), validity of
$|f(x)-P_n(x)|\le
c(k,r)\rho_n^r(x)\omega_k\left(f^{(r)},\rho_n(x)\right)$ for $x\in
[-1,1]$  and $n\geq \NN$ with $P_n\in\Delta^{(q)}$}\label{fig8}
\end{figure}

These results appeared in the papers by Beatson \cite{Be1980}, Cao
and Gonska \cite{CG1994}, and Bondarenko and Prymak \cite{BP2004}.

\samefigure
\begin{figure}[H]
{\scriptsize
$$
\begin{matrix}
 r &\vdots  &\vdots&\vdots&\vdots&\vdots&\updots
      \cr 3 &\hboxx{\cite{BP2004}}&\hboxx{\cite{BP2004}}&\hboxx{\cite{BP2004}}
      &\hboxx{\cite{BP2004}}&\hboxx{\cite{BP2004}}& \cdots&
      \cr 2 &\hboxx{\cite{CG1994}}&\hboxx{\cite{BP2004}}&\hboxx{\cite{BP2004}}
      &\hboxx{\cite{BP2004}}&\hboxx{\cite{BP2004}}&\cdots&
      \cr 1 &\hboxx{\cite{Be1980}}&\hboxx{\cite{CG1994}}&\hboxx{\cite{BP2004}}
      &\hboxx{\cite{BP2004}}&\hboxx{\cite{BP2004}}&\cdots&
      \cr 0& &\hboxx{\cite{Be1980}}&\hboxx{\cite{CG1994}}&\hboxx{\cite{BP2004}}
      &\hboxx{\cite{BP2004}}&\cdots&
      \cr &0&1 &2&3&4&k&
\end{matrix}
$$
}
\caption{References for \newfig{fig8}} \label{fig8a}
\end{figure}
\esamefigure

\begin{figure}[H]
$$
\begin{matrix}
 r &\vdots  &\vdots&\vdots&\vdots&\vdots&\vdots&\vdots&\updots
      \cr 5 &?^*&?^*   &?^*  &?^* &?^* &?^*& ?^*& \cdots&
      \cr 4 &?&?^*   &?^*  &?^*  &?^* &?^* &?^* & \cdots&
      \cr 3 &?& ? &?^*  &?^* &?^* &?^* & ?^* & \cdots&
      \cr 2 &+& ? &? &-&-&-& -& \cdots&
      \cr 1 &+& +&? & ?& -& - &  -& \cdots&
      \cr 0& &+&+&? & ?&-& - &\cdots&
      \cr &0&1 &2&3&4&5&6&k&
\end{matrix}
$$
\caption{$3$-monotone approx. ($q=3$), validity of $|f(x)-P_n(x)|\le
c(k,r)\rho_n^r(x)\omega_k\left(f^{(r)},\rho_n(x)\right)$ for $x\in
[-1,1]$  and $n\geq \NN$ with $P_n\in\Delta^{(3)}$}\label{fig9}
\end{figure}

\begin{remark}
 It follows from  Bondarenko and Gilewicz \cite{BG2002}
  that ``$?^*$'' in the cases
$(k\ge 0,r\geq 5)$, $(k\ge 1,r=4)$ and $(k\ge 2,r=3)$  in \newfig{fig9} cannot be replaced by ``$+$''.
\end{remark}

These results appeared in the papers by Beatson \cite{Be1980}, Cao
and Gonska \cite{CG1994}, Wu and Zhou \cite{WZ1992}, and Yushchenko \cite{Y2002}.

\samefigure
\begin{figure}[H]
{\scriptsize
$$
\begin{matrix}
 r &\vdots  &\vdots&\vdots&\vdots&\vdots&\vdots&\vdots&\updots
      \cr 5 &\hboxx{\cite{BG2002}}&\hboxx{\cite{BG2002}}   &\hboxx{\cite{BG2002}}  &\hboxx{\cite{BG2002}}
      &\hboxx{\cite{BG2002}} &\hboxx{\cite{BG2002}}& \hboxx{\cite{BG2002}}& \cdots&
      \cr 4 &?&\hboxx{\cite{BG2002}}    &\hboxx{\cite{BG2002}}  &\hboxx{\cite{BG2002}} &\hboxx{\cite{BG2002}}  &\hboxx{\cite{BG2002}} & \hboxx{\cite{BG2002}} & \cdots&
      \cr 3 &?& ? &\hboxx{\cite{BG2002}} &\hboxx{\cite{BG2002}} &\hboxx{\cite{BG2002}} &\hboxx{\cite{BG2002}} & \hboxx{\cite{BG2002}} &\cdots&
      \cr 2 &\hboxx{\cite{CG1994}}&?&? &\hboxx{\cite{Y2002}}
      &\hboxx{\cite{Y2002}}&\hboxx{\cite{Y2002}}&\hboxx{\cite{Y2002}}&\cdots&
      \cr 1 &\hboxx{\cite{Be1980}}&\hboxx{\cite{CG1994}}&?&?
      &\hboxx{\cite{Y2002}}&\hboxx{\cite{Y2002}}&\hboxx{\cite{Y2002}}&\cdots&
      \cr 0& &\hboxx{\cite{Be1980}}&\hboxx{\cite{CG1994}}&?&?
      &\hboxx{\cite{Y2002}}&\hboxx{\cite{WZ1992}}&\cdots&
      \cr &0&1 &2&3&4&5&6&k&
\end{matrix}
$$
}
\caption{References for \newfig{fig9}} \label{fig9a}
\end{figure}
\esamefigure

Again, we would like to replace the question marks with definitive
answers.

\begin{problem}
Determine which symbol $($among ``$+$'', ``$\ominus$'' and
``$-$''$)$ replaces ``$?$'' or ``$?^*$'' in {\rm\newfig{fig9}} in any
of the places.
\end{problem}

\begin{remark}
For the triples $(k\le2, r\le 2-k, q\ge 1)$, DeVore and
Yu\hboxx{\cite{DY1985}} (for $q=1$), Leviatan\hboxx{\cite{Le1986}}
(for $q=2$) and Cao and Gonska\hboxx{\cite{CG1994}} (for $q\in\N$),
proved that \ineq{eq7} and \ineq{eq9} are satisfied with $\rho_n(x)$
replaced by the smaller quantity $n^{-1} \varphi(x)$, so that, in
particular, the polynomials interpolate the function at the
endpoints.
One can also achieve interpolation at the endpoints if $k\geq 3$
(see \cite{K1994}*{(8)}, for example), but it is known that, in
general, the quantity $\rho_n(x)$ cannot be replaced by  $n^{-1}
\varphi(x)$ (see, \eg  \cites{yu85, li, dah, kop-sim, glsw} for
details).
\end{remark}

\section{Ditzian-Totik type estimates for $q$-monotone approximation, $q\geq 1$
} \label{S8}

In this section, we discuss the validity of the following
statements.

\begin{proposition} \label{prop6}
If $k\in\N$, $r\in\N_0$, $q\in\N$, $\NN\in\N$, and
  $f\in\Delta^{(q)}\cap \C^r_\varphi$, then
$$
E_n^{(q)}(f)\le
\frac{c(k,r,q)}{n^r}\omega^\varphi_{k,r}\bigl(f^{(r)},n^{-1}
\bigr),\quad n\ge \NN.
$$
\end{proposition}

If $k=0$, \prop{prop6} is modified as follows.

\begin{proposition} \label{prop7}
If   $r\in\N$, $q\in\N$, $\NN\in\N$, and $f\in\Delta^{(q)}\cap
\B^r$, then
$$
E_n^{(q)}(f)\le \frac{c(r,q)}{n^r}\|\varphi^r
f^{(r)}\|_{\LL_\infty[-1,1]},\quad n\ge \NN.
$$
\end{proposition}

The following are the ``truth tables'' for \props{prop6}{prop7}.

\begin{figure}[H]
$$
\begin{matrix}
 r &\vdots  &\vdots&\vdots&\vdots&\vdots&\vdots&\updots
      \cr 4 &+&+ &+  &+  &+ &+&\cdots&
      \cr 3 &+&+ &+  &+  &+ &+&\cdots&
      \cr 2 &+&\ominus &-  &-  &- &-&\cdots&
      \cr 1 &+& +&\ominus&- &- &-&\cdots&
      \cr 0& &+&+&\ominus&-&-&\cdots&
      \cr &0&1 &2&3&4&5&k&
\end{matrix}
$$
\caption{Monotone approx. ($q=1$), validity of $E_n^{(1)}(f)\leq
 c(k,r) n^{-r}\omega^\varphi_{k,r}\left(f^{(r)},n^{-1} \right)$,   $n\geq \NN$}\label{fig12}
\end{figure}

The results appear in the papers by Shvedov\hboxx{\cite{S1979}},
Leviatan\hboxx{\cite{Le1986}}, Wu and Zhou\hboxx{\cite{WZ1992}},
Dzyubenko, Listopad and Shevchuk\hboxx{\cite{DLS1993}}, Kopotun and
Listopad\hboxx{\cite{KL1994}}, Kopotun\hboxx{\cite{K1995}}, Leviatan
and Shevchuk\hboxx{\cites{LS1998, LS2000}},  and Nesterenko and
Petrova\hboxx{\cite{NP2005}}.

Recalling that, in the case ``$\ominus$'', a reference to the
negative result is followed by a reference to the positive one, we
summarize all  references for \newfig{fig12} in the following table.

\samefigure
\begin{figure}[H]
{\scriptsize
$$
\begin{matrix}
 r &\vdots  &\vdots&\vdots&\vdots&\vdots&\vdots&\updots
      \cr 4 &\hboxx{\cite{DLS1993}}&\hboxx{\cite{K1995}}&\hboxx{\cite{K1995}}&\hboxx{\cite{K1995}}
      &\hboxx{\cite{K1995}}&\hboxx{\cite{K1995}}&\cdots&
      \cr 3 &\hboxx{\cite{DLS1993}}&\hboxx{\cite{K1995}}&\hboxx{\cite{K1995}}&\hboxx{\cite{K1995}}
      &\hboxx{\cite{K1995}}&\hboxx{\cite{K1995}}&\cdots&
      \cr 2 &\hboxx{\cite{Le1986}}&\hboxx{\cite{KL1994}},\hboxx{\cite{LS1998}}&\hboxx{\cite{NP2005}}
      &\hboxx{\cite{LS2000}}&\hboxx{\cite{LS2000}}&\hboxx{\cite{LS2000}}&\cdots&
      \cr 1 &\hboxx{\cite{Le1986}}& \hboxx{\cite{Le1986}}&\hboxx{\cite{KL1994}},\hboxx{\cite{LS1998}}
      &\hboxx{\cite{NP2005}} &\hboxx{\cite{LS2000}}&\hboxx{\cite{LS2000}}&\cdots&
      \cr 0& &\hboxx{\cite{Le1986}}&\hboxx{\cite{Le1986}}&\hboxx{\cite{S1979}},\hboxx{\cite{LS1998}}
      &\hboxx{\cite{WZ1992}}&\hboxx{\cite{WZ1992}}&\cdots&
      \cr &0&1 &2&3&4&5&k&
\end{matrix}
$$
}
\caption{References for \newfig{fig12}}\label{fig12a}
\end{figure}
\esamefigure

\begin{figure}[H]
$$
\begin{matrix}
 r &\vdots  &\vdots&\vdots&\vdots&\vdots&\vdots&\vdots&\updots
      \cr 6 &+&+&+&+&+&+&+&\cdots&
      \cr 5 &+&+&+&+&+&+&+&\cdots&
      \cr 4 &\ominus&\ominus&-&-&-&-&-&\cdots&
      \cr 3 &+&\ominus&\ominus&-&-&-&-&\cdots&
      \cr 2 &+&+&\ominus&\ominus&-&-&-&\cdots&
      \cr 1 &+&+&+&\ominus&-&-&-&\cdots&
      \cr 0& &+&+&+&\ominus&-&-&\cdots&
      \cr &0&1 &2&3&4&5&6&k&
\end{matrix}
$$
\caption{Convex approx. ($q=2$), validity of $E_n^{(2)}(f)\leq
 c(k,r) n^{-r}\omega^\varphi_{k,r}\left(f^{(r)},n^{-1} \right)$,   $n\geq \NN$}\label{fig13}
\end{figure}

We emphasize that \newfig{fig13} shows that we now have not only cases
when \prop{prop6} is invalid, but also the case $(r=4, k=0)$ when
\prop{prop7} is valid only with $\NN$ that depends on the function
$f$. This is in contrast with \newfig{fig12} for monotone approximation
that did not  contain ``$\ominus$'' in the column corresponding to
$k=0$.

These results appeared in the papers by Shvedov\hboxx{\cite{S1979}},
Leviatan\hboxx{\cite{Le1986}},  Wu and Zhou\hboxx{\cite{WZ1992}},
Mania (see \cite{M1992}), Kopotun\hboxx{\cite{K1992}}, \cites{K1994,
K1995}, Leviatan and Shevchuk\hboxx{\cite{LS2003}}, Kopotun,
Leviatan and Shevchuk\hboxx{\cite{KLS2005}}, and Nissim and
Yushchenko\hboxx{\cite{NY2003}}.

\samefigure
\begin{figure}[H]
{\scriptsize
$$
\begin{matrix}
 r &\vdots  &\vdots&\vdots&\vdots&\vdots&\vdots&\vdots&\updots
      \cr 6 &\hboxx{\cite{K1992}}&\hboxx{\cite{K1995}}&\hboxx{\cite{K1995}}&\hboxx{\cite{K1995}}
      &\hboxx{\cite{K1995}}&\hboxx{\cite{K1995}}&\hboxx{\cite{K1995}}&\cdots&
      \cr 5 &\hboxx{\cite{K1992}}&\hboxx{\cite{K1995}}&\hboxx{\cite{K1995}}&\hboxx{\cite{K1995}}
      &\hboxx{\cite{K1995}}&\hboxx{\cite{K1995}}&\hboxx{\cite{K1995}}&\cdots&
      \cr 4 &\hboxx{\cite{K1992}},\hboxx{\cite{LS2003}}&\hboxx{\cite{K1992}},\hboxx{\cite{KLS2005}}
      &\hboxx{\cite{KLS2005}}&\hboxx{\cite{KLS2005}}&\hboxx{\cite{KLS2005}}
      &\hboxx{\cite{KLS2005}}&\hboxx{\cite{KLS2005}}&\cdots&
      \cr 3 &\hboxx{\cite{K1992}}&\hboxx{\cite{K1992}},\hboxx{\cite{LS2003}}&\hboxx{\cite{K1992}},
      \hboxx{\cite{KLS2005}}&\hboxx{\cite{KLS2005}}&\hboxx{\cite{KLS2005}}&
      \hboxx{\cite{KLS2005}}&\hboxx{\cite{KLS2005}}&\cdots&
      \cr 2 &\hboxx{\cite{Le1986}}&\hboxx{\cite{K1994}}&\hboxx{\cite{K1992}},\hboxx{\cite{LS2003}}
      &\hboxx{\cite{K1992}},\hboxx{\cite{KLS2005}}&\hboxx{\cite{KLS2005}}&
      \hboxx{\cite{KLS2005}}&\hboxx{\cite{KLS2005}}&\cdots&
      \cr 1 &\hboxx{\cite{Le1986}}&\hboxx{\cite{Le1986}}&\hboxx{\cite{K1994}}&\hboxx{\cite{M1992}},
      \hboxx{\cite{LS2003}}& \hboxx{\cite{NY2003}}&\hboxx{\cite{NY2003}}&\hboxx{\cite{NY2003}}&\cdots&
      \cr 0& &\hboxx{\cite{Le1986}}&\hboxx{\cite{Le1986}}&\hboxx{\cite{K1994}}&\hboxx{\cite{S1979}},
      \hboxx{\cite{LS2003}}&\hboxx{\cite{WZ1992}}&\hboxx{\cite{WZ1992}}&\cdots&
      \cr &0&1 &2&3&4&5&6&k&
\end{matrix}
$$
}
\caption{References for \newfig{fig13}}\label{fig13a}
\end{figure}
\esamefigure

\begin{figure}[H]
$$
\begin{matrix}
r &\vdots  &\vdots&\vdots&\vdots&\vdots&\vdots&\vdots&\updots
      \cr 5 &?&?   &?  &? &? &?& ?& \cdots&
      \cr 4 &?&?   &?  &? &? &?& ?& \cdots&
      \cr 3 &+& ? &? &?&?&?& ?& \cdots&
      \cr 2 &+& + &? &?^*&-&-& -& \cdots&
      \cr 1 &+& +&+ & ?& ?^* & - &  -& \cdots&
      \cr 0& &+&+&+ & ?&?^*& - &\cdots&
      \cr &0&1 &2&3&4&5&6&k&
\end{matrix}
$$
\caption{ $3$-monotone approx. ($q=3$), validity of $E_n^{(3)}(f)\le
c(k,r) n^{-r}\omega^\varphi_{k,r}\left(f^{(r)},n^{-1} \right)$,
$n\ge \NN$ }\label{fig5new}
\end{figure}
Since
\begin{equation} \label{newJuly9}
\omega^\varphi_{k,r}\left(f^{(r)},n^{-1} \right) \leq
\omega_{k}\left(f^{(r)},n^{-1} \right) ,
\end{equation}
all positive estimates in \newfig{fig5new} imply positive estimates in
corresponding places in \newfig{fig5} (and all negative results in
\newfig{fig5} imply negative results in corresponding places in
\newfig{fig5new}). At this time, \newfig{fig5new} is identical with
\newfig{fig5}
(all cases ``$+$'' follow from the article by
Bondarenko\hboxx{\cite{Bo2002}}, and all other cases are determined
by corresponding entries in \newfig{fig5}). However, we emphasize that
there is no guarantee that once all question marks  are replaced by
definitive symbols these tables will remain identical.

\begin{problem}
Determine which symbol $($among ``$+$'', ``$\ominus$'' and $``-")$
replaces $``?"$ or $``?^*"$ in {\rm\newfig{fig5new}} in  any of the
places.
\end{problem}

\begin{figure}[H]
$$
\begin{matrix}
 r &\vdots  &\vdots&\vdots&\vdots&\vdots&\updots
      \cr 4 &-&-  &-  &- &- & \cdots&
      \cr 3 &- &-   &-  &- &- & \cdots&
      \cr 2 &+&- &-  &-  &- &\cdots&
      \cr 1 &+& +&-&- &-  & \cdots&
      \cr 0& &+&+&- &- &\cdots&
      \cr &0&1 &2&3&4&k&
\end{matrix}
$$
\caption{ $q$-monotone approx. ($q\geq 4$), validity of
$E_n^{(q)}(f)\le c(k,r)
n^{-r}\omega^\varphi_{k,r}\left(f^{(r)},n^{-1} \right)$, $n\ge \NN$
}\label{fig4new}
\end{figure}

Because of \ineq{newJuly9}, all cases ``$-$'' in \newfig{fig4new}
follow from corresponding   ``$-$'''s in \newfig{fig4}. The cases
``$+$'' can be derived from results in the article by Gavrea,
Gonska, P{\u a}lt{\u a}nea and Tachev\hboxx{\cite{GGP2003}},
combined with the $q$-monotonicity preservation properties of the
Gavrea operators (see Gavrea\hboxx{\cite{Ga1996}}), appearing in the
paper of Cottin, Gavrea, Gonska, Kacs\'o and
Zhou\hboxx{\cite{CGG1999}}.

\section{Relations between   degrees of best unconstrained and $q$-monotone approximation}\label{S6}

Clearly, for each $f\in \C[-1,1]$ and $q\geq 1$,
$$
E_n(f)\le E_n^{(q)}(f).
$$
Moreover, Lorentz and Zeller\hboxx{\cite{LZ1968}},  proved that
there is a function $f\in\Delta^{(q)} \cap \C^q[-1,1]$, such that
\begin{equation}
\limsup_{n\to\infty}\frac{ E_n^{(q)}(f)}{ E_n(f)}=\infty. \label{eq11}
\end{equation}
It is well known that in the unconstrained case, for all $f\in
\C^r[-1,1]$, we have
$$
E_n(f)\le\frac{c(r)}{n^r}E_{n-r}(f^{(r)}),\quad n>r.
$$

At the same time, it was shown in \hboxx{\cite{LS1995}}  and
\hboxx{\cite{S1996}} that,  for each $n>q$, there is a function
$f_n\in\Delta^{(q)}\cap \C^q[-1,1]$, such that
\begin{equation}
E_n^{(q)}(f_n) >  c(q)E_{n-q}(f_n^{(q)}),\quad c(q)>0,\label{eq12}
\end{equation}
and so  the (almost obvious) estimate \ineq{eq2} may not in general be improved.

\begin{remark}
It was proved in \cite{BP2004} that, if $q\geq 4$ and  $r\geq q-1$,
then for any nonnegative sequence  $\{\alpha_n\}$ such that
$\lim_{n\to\infty} \alpha_n = \infty$, there exists a function $f =
f_{r,q} \in\Delta^{(q)} \cap \C^r[-1,1]$ for which
 \begin{equation} \label{inbp}
 \limsup_{n\to\infty} \alpha_n E_n^{(q)}(f)n^{r-q+3} = \infty .
\end{equation}

\end{remark}

\begin{problem} \label{prob-July9}
Given $q, r\in\N$ such that $r\geq q+1$, is it true that, for each
$n>r$, there is a function $f_n\in\Delta^{(q)}\cap \C^r[-1,1]$, such
that
\[
E_n^{(q)}(f_n) > c(q,r)E_{n-r}(f_n^{(r)}),\quad c(q,r)>0 \, ?
\]
\end{problem}

\begin{problem} \label{prob-July9-2}
Given $q, r\in\N$ such that $r\geq q+1$ and $1\leq q\leq 3$, is it
true that, for each $f\in\Delta^{(q)}\cap \C^r[-1,1]$,
\[
E_n^{(q)}(f)\leq c(q,r) n^{-r} E_{n-r}(f^{(r)}),\quad n>r \, ?
\]
\end{problem}
Note that, in the case $q\geq 4$,   estimate \ineq{inbp} implies
that the answer  to the question posed in \op{prob-July9-2}  is
``no''.

\begin{problem} \label{prob3}
Given $q, r\in\N$, does there exist $f\in\Delta^{(q)}\cap
\C^r[-1,1]$ such that
$$
\limsup_{n\to\infty}\frac{ E_n^{(q)}(f)}{ E_{n-r}(f^{(r)})}>0\,?
$$
\end{problem}

The following  is a weaker version of \op{prob3}.

\begin{problem} \label{newopr}
Given $1\leq q\leq 3$ and $r\geq q$, does there exist
$f\in\Delta^{(q)}\cap \C^r[-1,1]$  such that
$$
\limsup_{n\to\infty}\frac{ n^r
E_n^{(q)}(f)}{E_{n-r}(f^{(r)})}=\infty\,?
$$
\end{problem}
Note that, for $q\geq 4$,
``$-$'''s in Table \ref{fig4} imply that the answer to the question
posed in \op{newopr}  is ``yes''. The answer is also ``yes'' if
$0\leq r<q\leq 3$ (this follows from the cases ``$-$'' in
Tables~\ref{fig2}, \ref{fig3} and \ref{fig5}).

\section{$\alpha$-relations for $q$-monotone approximation}\label{S7}

Notwithstanding \ineq{eq11}, it was shown in \cite{KLS2009} for
$q=2$ (the case $q=1$ is similar) that, for  each $\alpha>0$ and
$f\in\Delta^{(q)}$, $1\leq q \leq 2$,
$$
n^\alpha E_n(f)\le1,\quad n\ge1 \quad  \Longrightarrow \quad
n^\alpha E_n^{(q)}(f)\le c(\alpha),\quad n\ge1.
$$
What happens if we only have the information on the left-hand side
beginning from some fixed $\NN\ge2$? In other words, what can be
said about the validity
 of the
following statement?

\begin{proposition} \label{prop55}
Let $\alpha>0$ and $f\in\Delta^{(q)}$, $1\le q\le2$. Then
\begin{equation}
n^\alpha E_n(f)\le1,  \quad  n\ge \NN\quad\Longrightarrow\quad
n^\alpha E_n^{(q)}(f)\le
c(\alpha,\NN),\quad n\ge \NN^*. \label{eq13}
\end{equation}
\end{proposition}

We are interested in determining all cases for which \prop{prop55}
is valid and the exact dependence of
 $\NN^*$ on the various parameters involved
in \ineq{eq13}, namely, $\alpha$, $\NN$, and perhaps $f$ itself.

It turns out that the parameter that has an influence on the
behavior of $\NN^*$ is $\lceil{\alpha/2}\rceil$, where
$\lceil{\mu}\rceil$  is the ceiling function (\ie the smallest
integer not smaller than $\mu$).

\begin{figure}[H]
$$
\begin{matrix}
 \lceil{\alpha/2}\rceil &\vdots  &\vdots&\vdots&\vdots&\updots
      \cr 3 &+  &+  &+ &+ & \cdots&
      \cr 2  &+   &+  &+ &+ & \cdots&
      \cr 1 &+ &+ &\ominus  &\ominus &\cdots&
             \cr &1 &2&3&4&\NN&
\end{matrix}
$$
\caption{$\alpha$-relations for monotone approx. ($q=1$), validity
of ``$n^\alpha E_n(f)\le1,  \;  n\ge \NN   \;  \Longrightarrow \;
n^\alpha E_n^{(1)}(f)\le c(\alpha,\NN),\; n\ge \NN^*$''}
\label{fig10}
\end{figure}

\begin{figure}[H]
$$
\begin{matrix}
 \lceil{\alpha/2}\rceil &\vdots  &\vdots&\vdots&\vdots&\vdots&\updots
      \cr 4 &+&+  &+  &+ &+ & \cdots&
      \cr 3 &+ &+   &+  &+ &+ & \cdots&
      \cr 2 &+&+ &+ &\ominus  &\ominus &\cdots&
      \cr 1 &+& +&+&\ominus &\ominus  & \cdots&
       \cr &1 &2&3&4&5&\NN&
\end{matrix}
$$
\caption{$\alpha$-relations for convex approx. ($q=2$), validity of
``$n^\alpha E_n(f)\le1,\;n\ge \NN\;\Longrightarrow \; n^\alpha
E_n^{(2)}(f)\le c(\alpha,\NN),\;n\ge \NN^*$''} \label{fig11}
\end{figure}

We recall that ``+" in  \newfigs{fig10}{fig11} means that $\NN^*=\NN$,
and ``$\ominus$'' means that $\NN^*$ depends on $\alpha$, $\NN$ and,
in addition, must depend on $f$, that is, it cannot be independent
of $f$. Note that \newfigs{fig10}{fig11} show that \ineq{eq13} is
always true, that is, \prop{prop55} is always valid (perhaps with
$\NN^*$ depending on $f$).

For $q=2$, \newfig{fig11} was completed in \cites{KLS2009, KLS2010}.
For $q=1$, one can use the same arguments to complete \newfig{fig10}.

\begin{problem}
Determine $\alpha$-relations for $3$-monotone approximation, \ie
construct a table analogous to Tables  {\rm \ref{fig10}} and
 {\rm \ref{fig11}} in the case $q=3$.
\end{problem}

In the case $q\geq 4$ and $\alpha > 2$ it is known that the
implication \ineq{eq13} is invalid for any $\NN\in\N$, since
  there exists a function $f =
f_{\alpha, q} \in \Delta^{(q)}$ such that $E_n(f)\le n^{-\alpha}$,
$n\in\N$, but for some $\beta < \alpha $,  $E_n^{(q)}(f) \geq
c(\alpha, \beta, q) n^{-\beta}$, for infinitely many $n$.

Indeed, if $\alpha > q-2\geq 2$ then \ineq{inbp} (with $r:=\lceil
\alpha \rceil$ and $\alpha_n := n^\epsilon$) implies that  there
exists $f  \in\Delta^{(q)} \cap \C^r[-1,1]$ such that
$\limsup_{n\to\infty}   E_n^{(q)}(f)n^{r-q+3+\epsilon} = \infty$,
and so (if  $\epsilon >0$ is sufficiently small)
\[
E_n^{(q)}(f) \geq n^{-r+q-3-\epsilon} \geq n^{-\alpha+\epsilon} ,
\]
for infinitely many $n$. At the same time, $E_n (f) \leq c n^{-r}
\leq c n^{-\alpha}$.

It was also shown in \cite{BP2004} that $E_n^{(q)}(x_+^{q-1}) \geq
c(q)n^{-2}$ if $q\geq 4$. Therefore, for $2< \alpha \leq  q-1$, $E_n
(x_+^{q-1}) \leq c n^{-q+1} \leq c n^{-\alpha}$ and
$E_n^{(q)}(x_+^{q-1}) \geq c(q) n^{-\alpha+\epsilon}$.

\begin{problem}
Let $\alpha > 0$, $q\geq 4$, $f\in \Delta^{(q)}$ and $\NN\in\N$.
Determine the largest $\beta > 0$ such that
\[
n^\alpha E_n(f)\le1,  \quad  n\ge \NN \quad\Longrightarrow\quad
n^\beta E_n^{(q)}(f)\le c(\alpha, \beta, q, \NN),\quad n\ge \NN^*,
\]
and investigate the dependence of $\NN^*$ on $\alpha$, $\beta$, $q$,
$\NN$ and $f$.
\end{problem}

\section{Comonotone and coconvex approximation: introducing the case ``$\oplus$''}  \label{S9}

For $s\in\N$, let $Y_s=\{y_i\}_{i=1}^{s}$ be a collection of $s$
points $ -1<y_s<\dots <y_1<1$, and denote $y_0:=1$ and
$y_{s+1}:=-1$.

We say that $f\in\Delta^{(1)}(Y_s)$, if $f$ is continuous on
$[-1,1]$, changes monotonicity at the points $Y_s$, and is
nondecreasing on $[y_1, 1]$.   In other words, $f\in\C[-1,1]$ is in
$\Delta^{(1)}(Y_s)$ if it is
nondecreasing on $[y_{2i+1},y_{2i}]$
and
nonincreasing on $[y_{2i},y_{2i-1}]$.

We say that $f\in\Delta^{(2)}(Y_s)$, if $f\in \C[-1,1]$ and
 $Y_s$ is its set of inflection points  with $f$ being convex
in $[y_1,1]$ (\ie   $f$ is convex on $[y_{2i+1},y_{2i}]$
and concave on $[y_{2i},y_{2i-1}]$).

Note that a function $f$ may belong to more than one class
$\Delta^{(q)}(Y_s)$, for $q=1$ or $q=2$, both for different sets
$Y_s$ (with the same number of changes $s$), and for different sets
$Y_s$ with a different number of changes. For example, this may
happen when $f$ is constant on a subinterval of $[-1,1]$. Also $f$
may belong to $\Delta^{(1)}(Y_{s_1})$ and to
$\Delta^{(2)}(Y_{s_2})$.

In particular, if $f\in \C^q(-1,1)\cap\C[-1,1]$, then
$f\in\Delta^{(q)}(Y_s)$ if and only if
$$
f^{(q)}(x)\prod_{i=1}^s(x-y_i)\ge 0,\quad x\in(-1,1).
$$
For $f\in \Delta^{(q)}(Y_s)$, $q=1,2$, we define
$$
E_n^{(q)}(f,Y_s):=\inf_{P_n\in{\Delta^{(q)}(Y_s)}\cap\P_n}\norm{f-P_n}.
$$
Newman \cite{N1979} proved the first strongly positive result (the
case ``$+$'')  showing that, if $f\in\Delta^{(1)}(Y_s)$, then
$$
E_n^{(1)}(f,Y_s)\le c(s)\omega \left(f, n^{-1} \right),\quad n\ge1.
$$
 Furthermore, Shvedov \cite{S1981-co}   proved that, if
$f\in\Delta^{(1)}(Y_s)$, then
$$
E_n^{(1)}(f,Y_s)\le c(s)\omega_2\left(f, n^{-1} \right),\quad n\ge
\NN,
$$
where $\NN=\NN(Y_s)$, and that this estimate is no longer valid with
$\NN$ independent of $Y_s$. Namely, for each $A>0$, $n\in\N$ and
$s\ge 1$, there are a collection $Y_s$ and a function
$f\in\Delta^{(1)}(Y_s)$ such that, for every polynomial
$P_n\in\Delta^{(1)}(Y_s)\cap\P_n$,
$$
\|f-P_n\|_{\C[-1,1]}>A\omega_2\left(f,n^{-1}\right).
$$

Thus, we arrive at a phenomenon we have not seen before: the
constant $\NN$ cannot depend only on $k$, $r$, $s$ and $q$, and at
the same time, it does not have to (fully) depend on the function
$f$ itself. Rather, it  depends only on where $f$ changes its
monotonicity (the set $Y_s$). We refer to the cases where $\NN$
depends on the location of the changes in monotonicity or convexity
as ``weakly positive'' cases and denote them by ``$\oplus$''.

\section{Comonotone approximation: uniform and pointwise estimates} \label{S10}

\subsection{Jackson-Stechkin type   estimates for comonotone approximation} \label{subS10-1}

We begin by discussing the validity of the following statements.

\begin{proposition}
If $k\in\N$, $r\in\N_0$, $s\in\N$, $\NN\in\N$, and
$f\in\Delta^{(1)}(Y_s)\cap \C^r[-1,1]$, then
\begin{equation}
E^{(1)}_n(f,Y_s)\le\frac{c(k,r,s)}{n^r}\omega_k\bigl(f^{(r)}, n^{-1}
\bigr),\quad
n\ge \NN. \label{eq14}
\end{equation}
\end{proposition}

\begin{proposition}
If  $r\in\N$, $s\in\N$, $\NN\in\N$, and
 $f\in\Delta^{(1)}(Y_s)\cap \W^r$, then
$$
E^{(1)}_n(f,Y_s)\le\frac{c(r,s)}{n^r}\|f^{(r)}\|_{\LL_\infty[-1,1]},\quad
n\ge \NN.
$$
\end{proposition}

We have the following truth tables.

\begin{figure}[H]
$$
\begin{matrix}
 r &\vdots  &\vdots&\vdots&\vdots&\vdots&\vdots&\updots

      \cr 3 &+& + &+ &+&+&+&  \cdots&
      \cr 2 &+& + &+&+&+&+&  \cdots&
      \cr 1 &+& +&+ & \oplus& -& - &   \cdots&
      \cr 0& &+&\oplus&-& -&-& \cdots&
      \cr &0&1 &2&3&4&5&k&
\end{matrix}
$$
\caption{Comonotone approx. with $s=1$, validity of
$E^{(1)}_n(f,Y_1)\le c(k,r) n^{-r}\omega_k\left(f^{(r)}, n^{-1}
\right)$, $n\ge \NN$} \label{fig15}
\end{figure}

\begin{figure}[H]
$$
\begin{matrix}
r &\vdots  &\vdots&\vdots&\vdots&\vdots&\vdots&\updots&
      \cr 4 &+& + &+ &+&+&+&  \cdots
      \cr 3 &+& + &+ &+&+&+&  \cdots
      \cr 2 &+& + &+&\oplus&\oplus&\oplus&  \cdots
      \cr 1 &+& +&\oplus & \oplus& -& - &   \cdots
      \cr 0& &+&\oplus&-& -&-& \cdots
      \cr
      &0&1 &2&3&4&5&k&
\end{matrix}
$$
\caption{Comonotone approx. with $s=2$, validity of
$E^{(1)}_n(f,Y_2)\le c(k,r) n^{-r}\omega_k\left(f^{(r)}, n^{-1}
\right)$, $n\ge \NN$} \label{fig16}
\end{figure}

\begin{figure}[H]
$$
\begin{matrix}
 r   &\vdots  &\vdots&\vdots&\vdots&\vdots&\vdots& \vdots & \vdots & \updots  \cr
 s+2 &+& + &+ &+&+&+&  \cdots & + & \cdots   \cr
 s+1 &+& + &+ &+&+&+&  \cdots & + & \cdots  \cr
 s   &+& + & +&\oplus&\oplus&\oplus&  \cdots & \oplus  & \cdots  \cr
 s-1 &+& + &\oplus&\oplus&\oplus&\oplus& \cdots & \oplus  &  \cdots   \cr
 \vdots &\vdots& \vdots &\vdots&\vdots&\vdots&\vdots &  \cdots & \vdots  & \cdots  \cr
 2   &+& + &\oplus&\oplus&\oplus&\oplus& \cdots & \oplus  & \cdots \cr
 1   &+& +&\oplus & \oplus& -& - & \cdots & - &  \cdots  \cr
 0   & &+&\oplus&-& -&-& \cdots & - & \cdots  \cr
     &0&1 &2&3&4&5&  \cdots &  \cdots &   k \cr
\end{matrix}
$$
\caption{Comonotone approx. with $s\geq 3$, validity of
$E^{(1)}_n(f,Y_s)\le c(k,r,s) n^{-r}\omega_k\left(f^{(r)}, n^{-1}
\right)$, $n\ge \NN$} \label{fig17}
\end{figure}

The results for $r=0$ in   Tables~\ref{fig15}, \ref{fig16} and
\ref{fig17}  are by Newman\hboxx{\cite{N1979}} ($k=1$),
Shvedov\hboxx{\cite{S1981}}
 ($k=2$), and  Zhou \cite{Z1993} ($k\ge3$).
The ``$+$'' result for $(k=1, r=1)$ (and so for $(k=0, r=2)$) is due
to Beatson and Leviatan\hboxx{\cite{BL1983}}. All other cases are
due to Gilewicz and Shevchuk \cite{GS1996}, and Dzyubenko, Gilewicz
and Shevchuk\hboxx{\cite{DGS1998}}. We will not go into details.
Instead, we advise the interested reader to consult these two papers
for their exact results.

We now summarize all references for the three previous tables in one
table. We refer to \hboxx{\cite{GS1996}} and \hboxx{\cite{DGS1998}}
together, so we present them in the table as
\hboxx{\cite{GS1996}\&\cite{DGS1998}}.

\samefigure
\begin{figure}[H]
{\scriptsize
$$
\begin{matrix}
 r &\vdots  &\vdots&\vdots&\vdots&\vdots&\vdots&\updots
       \cr 3 &\hboxx{\cite{GS1996}\&\cite{DGS1998}}&\hboxx{\cite{GS1996}\&\cite{DGS1998}}
       &\hboxx{\cite{GS1996}\&\cite{DGS1998}}&\hboxx{\cite{GS1996}\&\cite{DGS1998}}
       &\hboxx{\cite{GS1996}\&\cite{DGS1998}}&\hboxx{\cite{GS1996}\&\cite{DGS1998}}&\cdots&
       \cr 2 &\hboxx{\cite{BL1983}}&\hboxx{\cite{GS1996}\&\cite{DGS1998}}
       &\hboxx{\cite{GS1996}\&\cite{DGS1998}}&\hboxx{\cite{GS1996}\&\cite{DGS1998}}
       &\hboxx{\cite{GS1996}\&\cite{DGS1998}}&\hboxx{\cite{GS1996}\&\cite{DGS1998}}&\cdots&
       \cr 1 &\hboxx{\cite{N1979}}&\hboxx{\cite{BL1983}}
       &\hboxx{\cite{GS1996}\&\cite{DGS1998}}&\hboxx{\cite{GS1996}\&\cite{DGS1998}}
       &\hboxx{\cite{GS1996}\&\cite{DGS1998}}&\hboxx{\cite{GS1996}\&\cite{DGS1998}}&\cdots&
       \cr 0 & &\hboxx{\cite{N1979}}&\hboxx{\cite{S1981}}&\hboxx{\cite{Z1993}}
       &\hboxx{\cite{Z1993}} &\hboxx{\cite{Z1993}}
       &\cdots&
      \cr &0&1 &2&3&4&5&k&
\end{matrix}
$$
}
\caption{References for Tables~\ref{fig15}, \ref{fig16} and
\ref{fig17}} \label{fig17a}
\end{figure}
\esamefigure

\begin{remark}
 Estimates for comonotone approximation were
first discussed, independently, by Iliev\hboxx{\cite{I1978}} and
Newman\hboxx{\cite{N1979}}. Iliev claimed the same estimates as
Newman's   for  $(k=1,r=0)$ in Tables~\ref{fig15}, \ref{fig16} and
\ref{fig17}.
 However, his proof was somewhat incomplete except in the case $s=1$.
\end{remark}

\subsection{Ditzian-Totik type   estimates for comonotone approximation} \label{subS10-2}

In this section, we investigate the validity of the following
statements.

\begin{proposition}
If $k\in\N$, $r\in\N_0$, $s\in\N$, $\NN\in\N$, and
 $f\in\Delta^{(1)}(Y_s)\cap \C^r_\varphi$, then
$$
E^{(1)}_n(f,Y_s)\le\frac{c(k,r,s)}{n^r}\omega^\varphi_{k,r}\bigl(f^{(r)},
n^{-1}\bigr),\quad n\ge \NN.
$$
\end{proposition}

\begin{proposition}
If  $r\in\N$, $s\in\N$, $\NN\in\N$, and
 $f\in\Delta^{(1)}(Y_s)\cap \B^r$, then
$$
E^{(1)}_n(f,Y_s)\le\frac{c(r,s)}{n^r}\|\varphi^rf^{(r)}\|_{\LL_\infty[-1,1]},\quad
n\ge \NN.
$$
\end{proposition}

We have the following truth tables.

\begin{figure}[H]
$$
\begin{matrix}
 r &\vdots  &\vdots&\vdots&\vdots&\vdots&\updots
      \cr 6 &+& + &+ &+&+&\cdots&
      \cr 5 &+& + &+ &+&+&\cdots&
      \cr 4 &\oplus&\oplus&\oplus &\oplus&\oplus&\cdots&
      \cr 3 &+&\oplus&\oplus&\oplus&\oplus&\cdots&
      \cr 2 &\oplus&\oplus&-&-&-&\cdots&
      \cr 1 &+&\oplus&\oplus&-&-&\cdots&
      \cr 0& &+&\oplus&-&-&\cdots&
      \cr &0&1 &2&3&4&k&
\end{matrix}
$$
\caption{Comonotone approx. with $s=1$, validity of
$E^{(1)}_n(f,Y_1)\le c(k,r) n^{-r}\omega^\varphi_{k,r}\left(f^{(r)},
n^{-1} \right)$, $n\ge \NN$} \label{fig18}
\end{figure}

\begin{figure}[H]
$$
\begin{matrix}
 r &\vdots  &\vdots&\vdots&\vdots&\vdots&\updots
      \cr 2s+4 &+&+&+&+&+&\cdots&
      \cr 2s+3 &+&+&+&+&+&\cdots&
      \cr 2s+2 &\oplus&\oplus&\oplus&\oplus&\oplus&\cdots&
      \cr 2s+1 &\oplus&\oplus&\oplus&\oplus&\oplus&\cdots&
      \cr\vdots&\vdots&\vdots&\vdots&\vdots&\vdots&\cdots&
      \cr 3 &\oplus&\oplus&\oplus&\oplus&\oplus&\cdots&
     \cr 2 &\oplus&\oplus&-&-&-&\cdots&
      \cr 1 &+&\oplus&\oplus&-&-&\cdots&
      \cr 0& &+&\oplus&-&-&\cdots&
      \cr &0&1 &2&3&4&k&
\end{matrix}
$$
\caption{Comonotone approx. with $s\geq 2$, validity of
$E^{(1)}_n(f,Y_s)\le c(k,r,s)
n^{-r}\omega^\varphi_{k,r}\left(f^{(r)}, n^{-1} \right)$, $n\ge
\NN$} \label{fig19}
\end{figure}

\begin{remark}
 Note the unique phenomenon of the case ``$+$''
for $(k=0, r=3)$  when the function has just one point of inflection
($s=1$); see Leviatan and Shevchuk \cite{LS1999a}.
\end{remark}

The results for $r=0,1,2$ in \newfigs{fig18}{fig19} are by
Shvedov\hboxx{\cite{S1981}},  Zhou\hboxx{\cite{Z1993}}, Leviatan and
Shevchuk \cites{LS1997, LS1999a, LS2000}, Kopotun and Leviatan
\cite{KL1997}, and Nesterenko and Petrova \cite{NP2005}.

Again, we   summarize all references for both tables in one table of
references.

\samefigure
\begin{figure}[H]
{\scriptsize
$$
\begin{matrix}
 r &\vdots  &\vdots&\vdots&\vdots&\vdots&\updots
      \cr 2s+4
      &\hboxx{\cite{LS1999a}}&\hboxx{\cite{LS1999a}}&\hboxx{\cite{LS1999a}}
      &\hboxx{\cite{LS1999a}}&\hboxx{\cite{LS1999a}}&\cdots&
      \cr 2s+3 &\hboxx{\cite{LS1999a}}&\hboxx{\cite{LS1999a}}&\hboxx{\cite{LS1999a}}
      &\hboxx{\cite{LS1999a}}&\hboxx{\cite{LS1999a}}&\cdots&
      \cr 2s+2 &\hboxx{\cite{LS1997},\cite{LS1999a}}&\hboxx{\cite{LS1997},\cite{LS1999a}}
      &\hboxx{\cite{LS1997},\cite{LS1999a}}
      &\hboxx{\cite{LS1997},\cite{LS1999a}}&\hboxx{\cite{LS1997},\cite{LS1999a}}&\cdots&
      \cr\vdots&\vdots&\vdots&\vdots&\vdots&\vdots&\cdots&
      \cr 4 &\hboxx{\cite{LS1997},\cite{LS1999a}}&\hboxx{\cite{LS1997},\cite{LS1999a}}
      &\hboxx{\cite{LS1997},\cite{LS1999a}}
      &\hboxx{\cite{LS1997},\cite{LS1999a}}&\hboxx{\cite{LS1997},\cite{LS1999a}}&\cdots&
      \cr 3 &\hboxx{\cite{LS1999a}}&\hboxx{\cite{LS1997},\cite{LS1999a}}&\hboxx{\cite{LS1997},\cite{LS1999a}}
      &\hboxx{\cite{LS1997},\cite{LS1999a}}&\hboxx{\cite{LS1997},\cite{LS1999a}}&\cdots&
      \cr 2 &\hboxx{\cite{LS1997}},\hboxx{\cite{KL1997}}&\hboxx{\cite{LS1997}},\hboxx{\cite{LS2000}}
      &\hboxx{\cite{NP2005}}&\hboxx{\cite{LS2000}} &\hboxx{\cite{LS2000}}&\cdots&
      \cr 1 &\hboxx{\cite{LS1997}}&\hboxx{\cite{LS1997}},\hboxx{\cite{KL1997}}
      &\hboxx{\cite{LS1997}},\hboxx{\cite{LS2000}} &\hboxx{\cite{NP2005}}
      &\hboxx{\cite{LS2000}}&\cdots&
      \cr 0 & &\hboxx{\cite{LS1997}}&\hboxx{\cite{S1981}},\hboxx{\cite{KL1997}}
      &\hboxx{\cite{Z1993}} &\hboxx{\cite{Z1993}}  & \cdots&
      \cr &0&1 &2&3&4&k&
\end{matrix}
$$
}
\caption{References for \newfigs{fig18}{fig19}} \label{fig19a}
\end{figure}
\esamefigure

\subsection{Pointwise estimates for comonotone approximation} \label{subS10-3}

We now discuss the validity of analogous pointwise estimates.

\begin{proposition} \label{prop12}
If $k\in\N$, $r\in\N_0$, $s\in\N$, $\NN\in\N$, and
 $f\in\Delta^{(1)}(Y_s)\cap \C^r[-1,1]$, then
there is a sequence $\{P_n\}_{n=1}^\infty$ of polynomials
$P_n\in\Delta^{(1)}(Y_s) \cap\P_n$, such that for every $n\ge \NN$
and each $x\in[-1,1]$, we have
\begin{equation}\label{eq15}
|f(x)-P_n(x)|\le c(k,r,s)\rho_n^r(x)\omega_k\bigl(f^{(r)},\rho_n(x)\bigr).  %
\end{equation}
\end{proposition}

\begin{proposition} \label{prop13}
If  $r\in\N$, $s\in\N$, $\NN\in\N$, and
 $f\in\Delta^{(1)}(Y_s)\cap \W^{r}$, then
there is a sequence $\{P_n\}_{n=1}^\infty$ of polynomials
$P_n\in\Delta^{(1)}(Y_s)\cap \P_n$, such that every $n\ge \NN$, and
each $x\in[-1,1]$, we have
\begin{equation}
|f(x)-P_n(x)|\le
c(r,s)\rho_n^r(x)\|f^{(r)}\|_{\LL_\infty[-1,1]}. \label{eq16}
\end{equation}
\end{proposition}

The truth tables for \props{prop12}{prop13}   are the following.

\begin{figure}[H]
$$
\begin{matrix}
 r &\vdots  &\vdots&\vdots&\vdots&\vdots&\vdots&\updots

      \cr 3 &+& + &+ &+&+&+&  \cdots&
      \cr 2 &+& + &+&+&+&+&  \cdots&
      \cr 1 &+& +&+ & \oplus& -& - &   \cdots&
      \cr 0& &+&\oplus&-& -&-& \cdots&
      \cr &0&1 &2&3&4&5&k&
\end{matrix}
$$
\caption{Comonotone approx. with $s=1$, validity of
$|f(x)-P_n(x)|\le
c(k,r)\rho_n^r(x)\omega_k\left(f^{(r)},\rho_n(x)\right)$ for $x\in
[-1,1]$ and $n\geq \NN$ with $P_n\in\Delta^{(1)}(Y_1)$ }
\label{fig20}
\end{figure}

\begin{figure}[H]
$$
\begin{matrix}
 r &\vdots  &\vdots&\vdots&\vdots&\vdots&\vdots&\updots

      \cr 3 &\oplus& \oplus&\oplus &\oplus&\oplus&\oplus&  \cdots&
      \cr 2 &\oplus& \oplus &\oplus&\oplus&\oplus&\oplus&  \cdots&
      \cr 1 &\oplus& \oplus&\oplus & \oplus& -& - &   \cdots&
      \cr 0& &\oplus&\oplus&-& -&-& \cdots&
      \cr &0&1 &2&3&4&5&k&
\end{matrix}
$$
\caption{Comonotone approx. with $s\geq 2$, validity of
$|f(x)-P_n(x)|\le
c(k,r,s)\rho_n^r(x)\omega_k\left(f^{(r)},\rho_n(x)\right)$ for $x\in
[-1,1]$ and $n\geq \NN$ with $P_n\in\Delta^{(1)}(Y_s)$ }
\label{fig21}
\end{figure}

The results for $r=0$ are by Shvedov \cite{S1981},
 Zhou \cite{Z1993}, and
Dzyubenko \cite{D1994}. All other results are due to Dzyubenko,
Gilewicz and Shevchuk \cite{DGS1998}.

We  summarize all references for both tables in one table of
references.

\samefigure
\begin{figure}[H]
{\scriptsize
$$
\begin{matrix}
 r &\vdots  &\vdots&\vdots&\vdots&\vdots&\vdots&\updots
      \cr 2 &\hboxx{\cite{DGS1998}}&\hboxx{\cite{DGS1998}}&\hboxx{\cite{DGS1998}}
      &\hboxx{\cite{DGS1998}}&\hboxx{\cite{DGS1998}}&\hboxx{\cite{DGS1998}}&\cdots&
      \cr 1 &\hboxx{\cite{DGS1998}}&\hboxx{\cite{DGS1998}}&\hboxx{\cite{DGS1998}}&\hboxx{\cite{DGS1998}}
      &\hboxx{\cite{DGS1998}}&\hboxx{\cite{DGS1998}}&\cdots&
      \cr 0& &\hboxx{\cite{DGS1998}}&\hboxx{\cite{S1981}},\hboxx{\cite{D1994}}&\hboxx{\cite{Z1993}}&\hboxx{\cite{Z1993}}
      &\hboxx{\cite{Z1993}}& \cdots&
      \cr &0&1 &2&3&4&5&k&
\end{matrix}
$$
}
\caption{References for \newfigs{fig20}{fig21}} \label{fig21a}
\end{figure}
\esamefigure

\section{Coconvex approximation: uniform and pointwise estimates} \label{S11}

\subsection{Jackson-Stechkin type estimates for coconvex approximation} \label{S11-1}

We use the structure  of section~\ref{S10} and begin with the
following two statements.

\begin{proposition}
If $k\in\N$, $r\in\N_0$, $s\in\N$, $\NN\in\N$, and
$f\in\Delta^{(2)}(Y_s)\cap \C^r[-1,1]$, then
$$
E^{(2)}_n(f,Y_s)\le\frac{c(k,r,s)}{n^r}\omega_k\bigl(f^{(r)},n^{-1}\bigr),\quad
n\ge \NN.
$$
\end{proposition}

\begin{proposition}
If   $r\in\N$, $s\in\N$, $\NN\in\N$, and
 $f\in\Delta^{(2)}(Y_s)\cap \W^r$, then
$$
E^{(2)}_n(f,Y_s)\le\frac{c(r,s)}{n^r}\|f^{(r)}\|_{\LL_\infty[-1,1]},\quad
n\ge \NN.
$$
\end{proposition}

 We have the following truth tables.

\begin{figure}[H]
$$
\begin{matrix}
 r &\vdots  &\vdots&\vdots&\vdots&\vdots&\vdots&\updots
      \cr 4 &+&+&+&+&+&+&\cdots&
      \cr 3 &+&+&+&+&+&+&\cdots&
      \cr 2 &+&+&+&\oplus&-&-&\cdots&
      \cr 1 &+&+&\oplus&-&-&-&\cdots&
      \cr 0 & &+&+&\oplus&-&-&\cdots&
      \cr &0&1&2&3&4&5&k&
\end{matrix}
$$
\caption{Coconvex approx. with $s=1$, validity of
$E^{(2)}_n(f,Y_1)\le c(k,r) n^{-r}\omega_k\left(f^{(r)}, n^{-1}
\right)$, $n\ge \NN$} \label{fig22}
\end{figure}

The results are due to Zhou\hboxx{\cite{Z1993}}, Kopotun, Leviatan
and Shevchuk\hboxx{\cite{KLS1999}}, Pleshakov and
Shatalina\hboxx{\cite{PS2000}}, Gilewicz and Yushchenko
\cite{GY2002}, and Leviatan and Shevchuk \cites{LS2002, LS2003}.

\samefigure
\begin{figure}[H]
{\scriptsize
$$
\begin{matrix}
 r &\vdots  &\vdots&\vdots&\vdots&\vdots&\vdots&\updots
       \cr 3 &\hboxx{\cite{LS2003}}&\hboxx{\cite{LS2003}}&\hboxx{\cite{LS2003}}
       &\hboxx{\cite{LS2003}}&\hboxx{\cite{LS2003}}&\hboxx{\cite{LS2003}}&\cdots&
      \cr 2 &\hboxx{\cite{LS2002}}&\hboxx{\cite{LS2003}}&\hboxx{\cite{LS2003}}&
      \hboxx{\cite{PS2000}},\hboxx{\cite{LS2003}}&\hboxx{\cite{GY2002}}&\hboxx{\cite{GY2002}}&  \cdots&
      \cr 1 &\hboxx{\cite{LS2002}}& \hboxx{\cite{LS2002}}&\hboxx{\cite{PS2000}},\hboxx{\cite{KLS1999}}&
      \hboxx{\cite{Z1993}}&\hboxx{\cite{Z1993}}&\hboxx{\cite{Z1993}}&\cdots&
      \cr 0& &\hboxx{\cite{LS2002}}&\hboxx{\cite{LS2002}}&\hboxx{\cite{PS2000}},\hboxx{\cite{KLS1999}}
      &\hboxx{\cite{Z1993}}&\hboxx{\cite{Z1993}}& \cdots&
      \cr &0&1 &2&3&4&5&k&
\end{matrix}
$$
}
\caption{References for \newfig{fig22}} \label{fig22a}
\end{figure}
\esamefigure

\begin{figure}[H]
$$
\begin{matrix}
 r &\vdots  &\vdots&\vdots&\vdots&\vdots&\vdots&\updots
      \cr 4 &\oplus& \oplus&\oplus &\oplus&\oplus&\oplus&  \cdots&
      \cr 3 &\oplus& \oplus&\oplus &\oplus&\oplus&\oplus&  \cdots&
      \cr 2 &\oplus& \oplus &\oplus&\oplus&-&-&  \cdots&
      \cr 1 &\oplus& \oplus&\oplus &-& -& - &   \cdots&
      \cr 0& &\oplus&\oplus&\oplus& -&-& \cdots&
      \cr &0&1 &2&3&4&5&k&
\end{matrix}
$$
\caption{Coconvex approx. with $s \ge 2$, validity of
$E^{(2)}_n(f,Y_s)\le c(k,r,s) n^{-r}\omega_k\left(f^{(r)}, n^{-1}
\right)$, $n\ge \NN$} \label{fig23}
\end{figure}

The results are due to Zhou\hboxx{\cite{Z1993}}, Kopotun, Leviatan
and Shevchuk\hboxx{\cite{KLS1999}},
 Gilewicz and
Yushchenko \cite{GY2002},  Leviatan and Shevchuk \cites{LS2002,
LS2003}, and Kopotun, Leviatan and Shevchuk \cite{KLS2006}  (some
earlier negative results are due to Pleshakov and Shatalina
\cite{PS2000}).

\samefigure
\begin{figure}[H]
{\scriptsize
$$
\begin{matrix}
 r &\vdots  &\vdots&\vdots&\vdots&\vdots&\vdots&\updots \cr
4
&\hboxx{\cite{KLS2006}},\hboxx{\cite{LS2003}}&\hboxx{\cite{KLS2006}},
\hboxx{\cite{LS2003}}&\hboxx{\cite{KLS2006}},\hboxx{\cite{LS2003}}
       &\hboxx{\cite{KLS2006}},\hboxx{\cite{LS2003}}&\hboxx{\cite{KLS2006}},
       \hboxx{\cite{LS2003}}&\hboxx{\cite{KLS2006}},\hboxx{\cite{LS2003}}&\cdots&\cr

 3 &\hboxx{\cite{LS2002}},\hboxx{\cite{KLS1999}}&\hboxx{\cite{LS2002}},\hboxx{\cite{LS2003}}
 &\hboxx{\cite{LS2002}},\hboxx{\cite{LS2003}}
       &\hboxx{\cite{LS2002}},\hboxx{\cite{LS2003}}&\hboxx{\cite{LS2002}},\hboxx{\cite{LS2003}}
       &\hboxx{\cite{LS2002}},
       \hboxx{\cite{LS2003}}&\cdots&\cr
2
&\hboxx{\cite{LS2002}}&\hboxx{\cite{LS2002}},\hboxx{\cite{KLS1999}}&\hboxx{\cite{LS2002}},\hboxx{\cite{LS2003}}&
      \hboxx{\cite{LS2002}},\hboxx{\cite{LS2003}}&\hboxx{\cite{GY2002}}&\hboxx{\cite{GY2002}}&  \cdots&\cr
 1 &\hboxx{\cite{LS2002}}& \hboxx{\cite{LS2002}}&\hboxx{\cite{LS2002}},\hboxx{\cite{KLS1999}}&
      \hboxx{\cite{Z1993}}&\hboxx{\cite{Z1993}}&\hboxx{\cite{Z1993}}&\cdots& \cr
 0& &\hboxx{\cite{LS2002}}&\hboxx{\cite{LS2002}}&\hboxx{\cite{LS2002}},\hboxx{\cite{KLS1999}}
      &\hboxx{\cite{Z1993}}&\hboxx{\cite{Z1993}}& \cdots& \cr
 &0&1 &2&3&4&5&k&
\end{matrix}
$$
}
\caption{References for \newfig{fig23}} \label{fig23a}
\end{figure}
\esamefigure

\subsection{Ditzian-Totik type   estimates for coconvex approximation} \label{subS11-2}

For the generalized D-T moduli, we have the following two
statements.

\begin{proposition}
If $k\in\N$, $r\in\N_0$, $s\in\N$, $\NN\in\N$, and
 $f\in\Delta^{(2)}(Y_s)\cap\C^\varphi_r$, then
$$
E^{(2)}_n(f,Y_s)\le\frac{c(k,r,s)}{n^r}\omega^\varphi_{k,r}\bigl(f^{(r)},n^{-1}\bigr),\quad
n\ge \NN.
$$
\end{proposition}

\begin{proposition}
If   $r\in\N$, $s\in\N$, $\NN\in\N$, and
 $f\in\Delta^{(2)}(Y_s)\cap \B^r$, then
$$
E^{(2)}_n(f,Y_s)\le\frac{c(r,s)}{n^r}\|\varphi^rf^{(r)}\|_{\LL_\infty[-1,1]},\quad
n\ge \NN.
$$
\end{proposition}

The truth tables in this case are the following.

\begin{figure}[H]
$$
\begin{matrix}
 r &\vdots  &\vdots&\vdots&\vdots&\vdots&\vdots&\updots

      \cr 8 &+& + &+ &+&+&+&  \cdots&
      \cr 7 &+& + &+ &+&+&+&  \cdots&
      \cr 6 &\oplus& \oplus &\oplus &\oplus&\oplus&\oplus&  \cdots&
      \cr 5 &\oplus& \oplus &\oplus &\oplus&\oplus&\oplus&  \cdots&
      \cr 4 &\oplus&\ominus &- &-&-&-&  \cdots&
      \cr 3 &\oplus& \oplus &\ominus &-&-&-&  \cdots&
      \cr 2 &+& \oplus &\oplus& \ominus&-&-&  \cdots&
      \cr 1 &+& +& \oplus&-& -& - &   \cdots&
      \cr 0& &+&+& \oplus& -&-& \cdots&
      \cr &0&1 &2&3&4&5&k&
\end{matrix}
$$
\caption{Coconvex approx. with $s=1$, validity of
$E^{(2)}_n(f,Y_1)\le c(k,r) n^{-r}\omega^\varphi_{k,r}\left(f^{(r)},
n^{-1} \right)$, $n\ge \NN$} \label{fig24}
\end{figure}

\begin{figure}[H]
$$
\begin{matrix}
 r &\vdots  &\vdots&\vdots&\vdots&\vdots&\vdots&\updots

      \cr 6 &\oplus& \oplus&\oplus &\oplus&\oplus&\oplus&  \cdots&
      \cr 5 &\oplus& \oplus&\oplus &\oplus&\oplus&\oplus&  \cdots&
      \cr 4 &\oplus& \oplus&- &-&-&-&  \cdots&
      \cr 3 &\oplus& \oplus&\oplus &-&-&-&  \cdots&
      \cr 2 &\oplus& \oplus &\oplus&\oplus&-&-&  \cdots&
      \cr 1 &\oplus& \oplus&\oplus &-& -& - &   \cdots&
      \cr 0& &\oplus&\oplus&\oplus& -&-& \cdots&
      \cr &0&1 &2&3&4&5&k&
\end{matrix}
$$

\caption{Coconvex approx. with $s\geq 2$, validity of
$E^{(2)}_n(f,Y_s)\le c(k,r,s)
n^{-r}\omega^\varphi_{k,r}\left(f^{(r)}, n^{-1} \right)$, $n\ge
\NN$} \label{fig25}
\end{figure}

\begin{remark}
It is interesting to note that, in some cases, estimates for
coconvex approximation with $s=1$ are ``better'' than those with
$s\geq 2$ (see all cases ``$+$'' in Table~\ref{fig24}), and, in some
other cases,  they are ``worse'' (see cases ``$\ominus$'').

\end{remark}

The results are due to Zhou \cite{Z1993}, Kopotun, Leviatan and
Shevchuk \cites{KLS1999,KLS2006}, Pleshakov and Shatalina
\cite{PS2000}, Gilewicz and Yuschenko \cite{GY2002},  and Leviatan
and Shevchuk \cite{LS2002}.

We summarize the references for both tables in one table.

\samefigure
\begin{figure}[H]
{\scriptsize
$$
\begin{matrix}
 r &\vdots  &\vdots&\vdots&\vdots&\vdots&\vdots&\updots
\cr 5
&\hboxx{\cite{KLS2006}}&\hboxx{\cite{KLS2006}}&\hboxx{\cite{KLS2006}}&\hboxx{\cite{KLS2006}}
&\hboxx{\cite{KLS2006}}&\hboxx{\cite{KLS2006}}&\cdots \cr 4
&\hboxx{\cite{KLS2006}}&\hboxx{\cite{KLS2006}}&\hboxx{\cite{KLS2006}}&\hboxx{\cite{KLS2006}}
&\hboxx{\cite{KLS2006}}&\hboxx{\cite{KLS2006}}&\cdots \cr 3
&\hboxx{\cite{LS2002}},\hboxx{\cite{KLS1999}}&\hboxx{\cite{KLS2006}}&\hboxx{\cite{KLS2006}}
&\hboxx{\cite{KLS2006}}&\hboxx{\cite{KLS2006}}&\hboxx{\cite{KLS2006}}&\cdots
\cr 2
&\hboxx{\cite{LS2002}},\hboxx{\cite{KLS1999}}&\hboxx{\cite{LS2002}},\hboxx{\cite{KLS1999}}
&\hboxx{\cite{LS2002}},\hboxx{\cite{KLS2006}}&
\hboxx{\cite{LS2002}},\hboxx{\cite{KLS2006}}
&\hboxx{\cite{GY2002}}&\hboxx{\cite{GY2002}}&\cdots \cr 1
&\hboxx{\cite{LS2002}},\hboxx{\cite{KLS1999}}&\hboxx{\cite{LS2002}},\hboxx{\cite{KLS1999}}
&\hboxx{\cite{PS2000}},\hboxx{\cite{KLS1999}}
&\hboxx{\cite{Z1993}}&\hboxx{\cite{Z1993}}&\hboxx{\cite{Z1993}}&\cdots
\cr 0 &
&\hboxx{\cite{LS2002}},\hboxx{\cite{KLS1999}}&\hboxx{\cite{LS2002}},\hboxx{\cite{KLS1999}}
&\hboxx{\cite{PS2000}},\hboxx{\cite{KLS1999}}
&\hboxx{\cite{Z1993}}&\hboxx{\cite{Z1993}}&\cdots \cr &0&1
&2&3&4&5&k&
\end{matrix}
$$
} \caption{References for \newfigs{fig24}{fig25}} \label{fig24a}
\end{figure}
\esamefigure

\subsection{Pointwise  estimates for coconvex approximation} \label{subS11-3}

Finally, we discuss the pointwise estimates for coconvex
approximation.

\begin{proposition}
If $k\in\N$, $r\in\N_0$, $s\in\N$, $\NN\in\N$, and
 $f\in\Delta^{(2)}(Y_s)\cap \C^r[-1,1]$, then
there is a sequence $\{P_n\}_{n=1}^\infty$ of polynomials
$P_n\in\Delta^{(2)}(Y_s)\cap\P_n$  such that  \ineq{eq15}  holds for
every $n\ge \NN$ and each $x\in[-1,1]$.
\end{proposition}

\begin{proposition}
If  $r\in\N$, $s\in\N$, $\NN\in\N$, and
 $f\in\Delta^{(2)}(Y_s)\cap \W^{r}$, then
there is a sequence $\{P_n\}_{n=1}^\infty$ of polynomials
$P_n\in\Delta^{(2)}(Y_s)\cap\P_n$  such that  \ineq{eq16}  holds for
every $n\ge \NN$  and each $x\in[-1,1]$.
\end{proposition}

The cases $s=1$ and $s\ge2$ are significantly different. We begin
with $s=1$,  where we have the following truth table.

\begin{figure}[H]
$$
\begin{matrix}
 r &\vdots  &\vdots&\vdots&\vdots&\vdots&\vdots&\updots
\cr 4 &\ominus&\ominus   &\ominus  &\ominus&\ominus &\ominus&\cdots
\cr 3 &\ominus&\ominus &\ominus &\ominus&\ominus &\ominus&\cdots \cr
2 &+&\ominus &\ominus &\ominus &- &- &\cdots \cr 1 &+& +&\ominus &-
&- &- &\cdots \cr 0& &+&+&\ominus &- &-&\cdots\cr &0&1 &2&3&4&5&k
\end{matrix}
$$
\caption{Coconvex approx. with $s=1$, validity of $|f(x)-P_n(x)|\le
c(k,r)\rho_n^r(x)\omega_k\left(f^{(r)},\rho_n(x)\right)$ for $x\in
[-1,1]$ and $n\geq \NN$ with $P_n\in\Delta^{(2)}(Y_1)$ }
\label{fig26}
\end{figure}

The results are due to  Zhou \cite{Z1993}, Gilewicz and Yushchenko
\cite{GY2002}, Dzyubenko, Gilewicz and Shevchuk \cite{DGS2002}, and
Dzyubenko, Leviatan and Shevchuk \cite{DLS2010a},
 \cite{DLS2010b}. We summarize the references in the following table.

\samefigure
\begin{figure}[H]
{\scriptsize
$$
\begin{matrix}
 r &\vdots  &\vdots&\vdots&\vdots&\vdots&\updots& \cr
 4 &\hboxx{\cite{DGS2002}, \cite{DLS2010a}}&\hboxx{\cite{DGS2002}, \cite{DLS2010a}}
      &\hboxx{\cite{DGS2002}, \cite{DLS2010a}}&\hboxx{\cite{DGS2002}, \cite{DLS2010a}}
      &\hboxx{\cite{DGS2002}, \cite{DLS2010a}}&\cdots&
      \cr
 3 &\hboxx{\cite{DGS2002}, \cite{DLS2010a}}&\hboxx{\cite{DGS2002}, \cite{DLS2010a}}
      &\hboxx{\cite{DGS2002}, \cite{DLS2010a}}&\hboxx{\cite{DGS2002}, \cite{DLS2010a}}
      &\hboxx{\cite{DGS2002}, \cite{DLS2010a}}&\cdots&
      \cr
2 &\hboxx{\cite{DGS2002}}&\hboxx{\cite{DGS2002}, \cite{DLS2010a}}
      &\hboxx{\cite{DGS2002}, \cite{DLS2010a}}&\hboxx{\cite{DGS2002}, \cite{DLS2010a}}&\hboxx{\cite{GY2002}}&\cdots&
      \cr 1 &\hboxx{\cite{DGS2002}}&\hboxx{\cite{DGS2002}}&\hboxx{\cite{DGS2002}, \cite{DLS2010b}}
      &\hboxx{\cite{Z1993}}&\hboxx{\cite{Z1993}}&\cdots&
      \cr 0 & &\hboxx{\cite{DGS2002}}&\hboxx{\cite{DGS2002}}&\hboxx{\cite{DGS2002}, \cite{DLS2010b}}
      &\hboxx{\cite{Z1993}}&\cdots&
      \cr &0&1 &2&3&4&k&
\end{matrix}
$$
}
\caption{References for \newfig{fig26}} \label{fig26a}
\end{figure}
\esamefigure

Now we will discuss the case $s\ge2$.

Recall that the case ``$\oplus$'' means that we have an estimate of
the form
\begin{equation} \label{sh1}
\left| f(x) - P_n(x)\right| \leq
 c(k,r,s,q) \rho_n^r(x) \omega_k\bigl(f^{(r)},\rho_n(x)\bigr), \quad n\ge \NN=\NN(k,r,q,Y_s),
\end{equation}
which in general is not true if the constant $\NN$ is independent of
$Y_s$. For comonotone approximation, a Whitney type estimate has
been proved in \cite{S1995} (in the case ``$\oplus$''):
\[
E_{k+r}^{(1)}(f,Y_s)\le c(k,r,Y_s) \omega_k(f^{(r)},1).
\]
This inequality, combined with \eqref{sh1}, implies for all
comonotone cases ``$\oplus$'' above that, together with the
estimate~\eqref{sh1}, we also have the estimate
\begin{equation} \label{sh2}
\left| f(x) - P_n(x)\right| \leq
 c(k,r,q,Y_s) \rho_n^r(x)\omega_k\bigl(f^{(r)},\rho_n(x)\bigr),\quad n\ge\NN= k+r.
\end{equation}
For the case of coconvex pointwise approximation with $s\ge2$, a new
case  appears, which is somewhere between ``$\oplus$" and ``$\ominus$".
Namely, inequality~\eqref{sh2} is valid (see \cite{PS2000} for the
corresponding Whitney inequality),  but the inequality~\eqref{sh1}
is not valid.

We refer to this case ``$\oslash$'' and formally define it as
follows.

\begin{description}
\item[Case ``$\oslash$'':]  \label{slashref}
Inequality
\[
\left| f(x) - P_n(x)\right| \leq
 c \rho_n^r(x) \omega_k\bigl(f^{(r)},\rho_n(x)\bigr), \quad n\ge \NN,
\]
holds with $c=c(k,r,Y_s)$ and $\NN=k+r$, as well as  with
$c=c(k,r,s)$ and $\NN=\NN(k,r,Y_s,f)$, but this inequality does not
hold with $c=c(k,r,s)$ and $\NN$ which may depend on $k$, $r$, and
$Y_s$, but is independent of $f$.
\end{description}

We have the following truth table.

\begin{figure}[H]
$$
\begin{matrix}
 r &\vdots  &\vdots&\vdots&\vdots&\vdots&\vdots&\updots
\cr 4 &\oslash&\oslash   &\oslash  &\oslash&\oslash
&\oslash&\cdots\cr 3 &\oslash&\oslash &\oslash &\oslash&\oslash
&\oslash&\cdots\cr 2 &\oplus&\oslash &\oslash &\oslash &- &-
&\cdots\cr 1 &\oplus& \oplus&\oslash &- &- &- &\cdots\cr 0&
&\oplus&\oplus&\oslash &- &-&\cdots \cr &0&1&2&3&4&5&k
\end{matrix}
$$
\caption{Coconvex approx. with $s\geq 2$, validity of
$|f(x)-P_n(x)|\le
c\rho_n^r(x)\omega_k\left(f^{(r)},\rho_n(x)\right)$ for $x\in
[-1,1]$ and $n\geq \NN$ with $P_n\in\Delta^{(2)}(Y_s)$ }
\label{fig27}
\end{figure}

The results are due to Zhou \cite{Z1993}, Gilewicz and Yushchenko
\cite{GY2002}, Leviatan and Shevchuk \cite{LS2002}, Dzyubenko,
Gilewicz and Shevchuk \cite{DGS2002}, \cite{DGS2006}, Dzyubenko and
Zalizko \cite{DZ2004}, \cite{DZ2005}, and Dzyubenko, Leviatan and
Shevchuk \cite{DLS2010a}, \cite{DLS2010b}. We summarize the
references in the following table.

\newcommand{\shrink}{}

\samefigure
\begin{figure}[H]

\resizebox{1.02\textwidth}{!} { {\scriptsize $
\begin{matrix}
 r &\vdots  &\vdots&\vdots&\vdots&\vdots & \updots& \cr
 3 &\hboxx{\shrink \cite{DGS2006}, \cite{DZ2004}, \cite{DLS2010a}}
 &\hboxx{\shrink \cite{DGS2006}, \cite{DZ2005}, \cite{DLS2010a}}
      &\hboxx{\shrink \cite{DGS2006}, \cite{DZ2005}, \cite{DLS2010a}}
      &\hboxx{\shrink \cite{DGS2006}, \cite{DZ2005}, \cite{DLS2010a}}
      &\hboxx{\shrink \cite{DGS2006},  \cite{DZ2005},\cite{DLS2010a}}&\cdots&\cr
2 &\hboxx{\cite{LS2002}, \cite{DGS2002}}&\hboxx{\shrink
\cite{DGS2006}, \cite{DZ2004}, \cite{DLS2010a}}
      &\hboxx{\shrink \cite{DGS2006}, \cite{DZ2005}, \cite{DLS2010a}}&\hboxx{\shrink
      \cite{DGS2006}, \cite{DZ2005}, \cite{DLS2010a}}&\hboxx{\cite{GY2002}}&\cdots&\cr
1 &\hboxx{\cite{LS2002}, \cite{DGS2002}}&\hboxx{\cite{LS2002},
\cite{DGS2002}} &\hboxx{\shrink \cite{DGS2006}, \cite{DZ2004},
\cite{DLS2010b}}
      &\hboxx{\cite{Z1993}}&\hboxx{\cite{Z1993}}&\cdots&\cr
0 & &\hboxx{\cite{LS2002}, \cite{DGS2002}}&\hboxx{\cite{LS2002},
\cite{DGS2002}} &\hboxx{\shrink \cite{DGS2006}, \cite{DZ2004},
\cite{DLS2010b}}
      &\hboxx{\cite{Z1993}}&\cdots&\cr
&0&1 &2&3&4&k&
\end{matrix}
$ } }

\caption{References for \newfig{fig27}} \label{fig27a}
\end{figure}
\esamefigure

\section{$\alpha$-relations for comonotone and coconvex approximation} \label{S12}

It follows from Leviatan and Shevchuk \cites{LS1997, LS1999} and
from Kopotun, Leviatan and Shevchuk \cite{KLS1999}  that,  if $Y_1$
and $\alpha>0$, $\alpha\ne 2$, are given, and if a function
$f\in\Delta^{(1)}(Y_1)$ satisfies
\begin{equation}
n^\alpha E_n(f)\le1,\quad n\ge1, \label{eq17}
\end{equation}
then
$$
n^\alpha E_n^{(1)}(f,Y_1)\le c(\alpha),\quad n\ge1.
$$
For coconvex approximation, Kopotun, Leviatan and Shevchuk
\cite{KLS2009} recently proved that if $Y_1$ and $\alpha>0$,
$\alpha\ne 4$, are given, and if a function $f\in\Delta^{(2)}(Y_1)$
satisfies \ineq{eq17}, then
$$
n^\alpha E_n^{(2)}(f,Y_1)\le c(\alpha),\quad n\ge1.
$$
For $\alpha=4$ this statement is false (see \cite{KLS2009} or
\cite{KLS2010}) since, for every $Y_1$, there is a function
$f\in\Delta^{(2)}(Y_1)$, satisfying  \ineq{eq17} with $\alpha=4$,
and
$$
\sup_{n\ge1}n^4 E_n^{(2)}(f,Y_1)> c \left|\ln\varphi(y_1)\right|.
$$
Similar arguments yield that, for every $Y_1$, there is a function
$f\in\Delta^{(1)}(Y_1)$, satisfying  \ineq{eq17}  with $\alpha=2$,
and
$$
\sup_{n\ge1}n^2 E_n^{(1)}(f,Y_1)> c \left|\ln\varphi(y_1)\right|.
$$
However, it is still true that, if a function
$f\in\Delta^{(2)}(Y_1)$ satisfies \ineq{eq17}  with $\alpha=4$, or
if a function $f\in\Delta^{(1)}(Y_1)$ satisfies \ineq{eq17} with
$\alpha=2$, then
$$
n^4E_n^{(2)}(f,Y_1)\le c,\quad n\ge \NN^*(Y_1)
$$
or
$$
n^2E_n^{(1)}(f,Y_1)\le c,\quad n\ge \NN^*(Y_1),
$$
respectively.

Also, for $s\ge2$, it has been shown in  \cite{KLS2009}  that,
 if $\,Y_s$ and $\alpha>0$, are given, and if a function
$f\in\Delta^{(2)}(Y_s)$ satisfies  \ineq{eq17}, then
$$
n^\alpha E_n^{(2)}(f,Y_s)\le c(\alpha,s),\quad n\ge
\NN^*(\alpha,Y_s),
$$
but we cannot have the above with $\NN^*$ that is independent of
$Y_s$.

Surprisingly, the case $s\ge2$ in comonotone approximation is more
complicated. Leviatan, Radchenko and Shevchuk \cite{LRS2011},
recently obtained the following.

\begin{theorem}
If $Y_s$, $s\ge1$, and $\alpha>0$ are given, and if a function
$f\in\Delta^{(1)}(Y_s)$ satisfies  $\eqref{eq17}$, then
$$
n^\alpha E_n^{(1)}(f,Y_s)\le c(\alpha,s),\quad n\ge \NN^*,
$$
where $\NN^*=\NN^*(\alpha,Y_s)$, if $\alpha=j,\,j=1,\dots,2 \lfloor\frac s2\rfloor$, or $\alpha=2j,\,j=1,\dots,s$, and $\NN^*=1$ in all other
cases.

Moreover, this statement cannot be improved since, for $s\ge1$ and
$\alpha=j,\,j=1,\dots,2\lfloor\frac s2\rfloor$, or $\alpha=2j$, $j=1,\dots,s$,
for every $m$ there are a collection $Y_s$ and a function
$f\in\Delta^{(1)}(Y_s)$ satisfying $\eqref{eq17}$ and
$$
m^\alpha E_m^{(1)}(f,Y_s)\ge c(s)\ln m,
$$
where $c(s)>0$, depends only on $s$.
\end{theorem}

For coconvex approximation, we discuss the validity of the following
statement.

\begin{proposition} \label{prop20}
If  $f\in\Delta^{(2)}(Y_s)$ and
\begin{equation}
n^\alpha E_n(f)\le1,\quad n\ge \NN, \label{eq18}
\end{equation}
then
$$
n^\alpha E_n^{(2)}(f,Y_s)\le c(s,\alpha),\quad n\ge \NN^*.
$$
\end{proposition}

Obviously $\NN^*\ge \NN$.  Is $\NN^*=\NN^*(\NN,\alpha,s)$ (strongly
positive), or $\NN^*=\NN^*(\NN,\alpha,Y_s)$ (weakly positive), or
$\NN^*=\NN^*(\NN,\alpha,Y_s,f)$ (weakly negative), or is it possible
that the above is invalid (strongly negative)?

We have the following truth tables which show, in particular, that
the strongly negative case is impossible. That is, \prop{prop20} is
always valid and there always exists an $\NN^*$. We observe that the
validity of \prop{prop20} depends on $\lceil{\frac\alpha2}\rceil$
rather than on $\alpha$ itself with only one exception. Namely,
 the symbol ``$\circplusnew$" for $\lceil{\frac\alpha2}\rceil=2$ and
$\NN=1,2$ in Table~\ref{fig28} is meant to indicate that $\alpha=4$ behaves
differently than the other $\alpha$'s with
$\lceil{\frac\alpha2}\rceil=2$.  Namely, we have ``$\circplusnew$''$:=$
``$\oplus$'' for $\alpha=4$, and ``$\circplusnew$''$:=$``$+$''
otherwise.

\begin{figure}[H]
$$
\begin{matrix}
 \lceil{\frac\alpha2}\rceil    &\vdots  &\vdots&\vdots&\vdots&\vdots&\vdots&\updots
      \cr 5 &+   &+  &+ &+ &+& +& \cdots&
      \cr 4 &+   &+  &+ &+ &+& +& \cdots&
      \cr 3 & + &+ &+&+&\oplus&\oplus& \cdots&
      \cr 2 &\circplusnew   &\circplusnew &\oplus&\oplus&\ominus& \ominus& \cdots&
      \cr 1 & +&+ & \oplus& \ominus & \ominus &  \ominus& \cdots&
      \cr &1 &2&3&4&5&6&\NN&
\end{matrix}
$$
\caption{$\alpha$-relations for coconvex approx. with $s=1$,
validity of ``$n^\alpha E_n(f)\le1,  \;  n\ge \NN   \;
\Longrightarrow \;   n^\alpha E_n^{(2)}(f,Y_1)\le c(\alpha),\; n\ge
\NN^*$''} \label{fig28}
\end{figure}

We  emphasize once again that,  in all cases ``$+$'' in
Table~\ref{fig28},  \prop{prop20} holds with $\NN^*=\NN$.

\begin{figure}[H]
$$
\begin{matrix}
  \lceil{\frac\alpha2}\rceil   &\vdots  &\vdots&\vdots&\vdots&\vdots&\vdots&\vdots&\updots
      \cr 4 &\oplus&  \oplus  &\oplus  &\oplus &\oplus &\oplus& \oplus& \cdots&
      \cr 3 &\oplus& \oplus &\oplus &\oplus&\oplus&\oplus& \oplus& \cdots&
      \cr 2 &\oplus& \oplus &\oplus &\oplus&\oplus&\ominus&\ominus& \cdots&
      \cr 1 &\oplus& \oplus&\oplus & \oplus& \ominus & \ominus &  \ominus& \cdots&
      \cr &1 &2&3&4&5&6&7&\NN&
\end{matrix}
$$
\caption{$\alpha$-relations for coconvex approx. with $s=2$,
validity of ``$n^\alpha E_n(f)\le1,  \;  n\ge \NN   \;
\Longrightarrow \;   n^\alpha E_n^{(2)}(f,Y_2)\le c(\alpha),\; n\ge
\NN^*$''}\label{fig29}
\end{figure}

\begin{figure}[H]
$$
\begin{matrix}
  \lceil{\frac\alpha2}\rceil   &\vdots  &\vdots&\vdots&\vdots&\vdots&\vdots&\vdots&\updots
      \cr 4 &\oplus   &\cdots  &\oplus &\oplus &\oplus& \oplus&\oplus& \cdots&
      \cr 3 &\oplus &\cdots &\oplus&\oplus&\oplus& \oplus&\oplus& \cdots&
      \cr 2 &\oplus &\cdots &\oplus&\oplus& ?^* &\ominus& \ominus&\cdots&
      \cr 1 &\oplus&\cdots & \oplus& \oplus & \ominus &  \ominus& \ominus&\cdots&
      \cr &1 &\cdots&s+1&s+2&s+3&s+4&s+5&\NN&
\end{matrix}
$$
\caption{$\alpha$-relations for coconvex approx. with $s\ge3$,
validity of ``$n^\alpha E_n(f)\le1,  \;  n\ge \NN   \;
\Longrightarrow \;   n^\alpha E_n^{(2)}(f,Y_s)\le c(s,\alpha),\;
n\ge \NN^*$''} \label{fig30}
\end{figure}

\begin{remark}
Note that ``$?^*$''  in \newfig{fig30}  can  be replaced neither by ``$+$'' nor by ``$-$''.
\end{remark}

All the results in Tables \ref{fig28}, \ref{fig29} and \ref{fig30},
have appeared in \cite{KLS2009} and \cite{KLS2010}.

\renewcommand{\thesection}{\hskip2ex}
\section[References]{}

{
\bigskip
% first author's name, affiliation, and address, including email and web page
\hskip1.4 em\vbox{\noindent K. A. Kopotun\\Department of Mathematics\\University of Manitoba \\
Winnipeg, Manitoba R3T 2N2\\Canada\\
 {\tt kopotunk@cc.umanitoba.ca }\\
 }

\hskip1.4 em\vbox{\noindent D. Leviatan\\Raymond and Beverly Sackler School of Mathematics\\
Tel Aviv University\\Tel Aviv 69978 \\
Israel\\
 {\tt leviatan@post.tau.ac.il }\\
}

% first author's name, affiliation, and address, including email and web page
\hskip1.4 em\vbox{\noindent A. Prymak\\Department of Mathematics\\University of Manitoba \\
Winnipeg, Manitoba R3T 2N2\\Canada\\
 {\tt prymak@cc.umanitoba.ca }\\
 }

\hskip1.4 em\vbox{\noindent I. A. Shevchuk\\National
Taras Shevchenko University of Kyiv \\ Kyiv, Ukraine\\
 {\tt shevchukh@ukr.net }\\
} }

\endddoc
\nocite{*}

\begin{bibsection}
\begin{biblist} %[\normalsize]

\bib{Be1978}{article}{
   author={Beatson, R. K.},
   title={The degree of monotone approximation},
   journal={Pacific J. Math.},
   volume={74},
   date={1978},
   number={1},
   pages={5--14},
   %%issn={0030-8730},
   %%review={\MR{0487165 (58 \#6825)}},
}


\bib{Be1980}{article}{
   author={Beatson, R. K.},
   title={Joint approximation of a function and its derivatives},
   conference={
      title={Approximation theory, III (Proc. Conf., Univ. Texas, Austin,
      Tex., 1980)},
   },
   book={
      publisher={Academic Press},
      place={New York},
   },
   date={1980},
   pages={199--206},
   %%review={\MR{602716 (82c:41019)}},
}


\bib{BL1983}{article}{
   author={Beatson, R. K.},
   author={Leviatan, D.},
   title={On comonotone approximation},
   journal={Canad. Math. Bull.},
   volume={26},
   date={1983},
   number={2},
   pages={220--224},
   %%issn={0008-4395},
   %%review={\MR{697804 (84f:41022)}},
   %%doi={10.4153/CMB-1983-034-0},
}


\bib{Bern1912}{article}{
author={Bernstein, S. N.}, title={D\'{e}monstration du
th\'{e}or\`{e}me de Weierstrass fond\'{e}e sur la calcul des
probabilit\'{e}s}, journal={Commun. Soc. Math.  Kharkov},
 volume={13},
date={1912},
 pages={1--2},
 }

\bib{B1927}{article}{
   author={Bernstein, S. N.},
   title={Sur les polynomes multiplement monotone},
 %  language={French},
   journal={Commun. Soc. Math.  Kharkov}, % taken from Lorentz "Bernstein polynomials"
   volume={4},
   date={1927},
   number={1},
   pages={1--11},
}

\bib{ber1937}{book}{
   author={Bernstein, S. N.},
   title={Ekstremal'nye svojstva polinomov.
 (Extremal properties of polynomials.)},
   language={Russian},
   publisher={Glavnaja Redakcija Obschetehnicheskoj Literatury (GROL)},
   place={Leningrad-Moscow},
   date={1937},
}


\bib{collected1}{book}{
   author={Bernstein, S. N.},
   title={Sobranie so\v cinen\u\i. Tom I. Konstruktivnaya teoriya funkci\u\i\
   [1905--1930] (Collected works. Vol. I. The constructive theory of functions [1905--1930].)},
   language={Russian},
   publisher={Akad. Nauk SSSR},
   place={Moscow},
   date={1952},
   pages={581 pp. (1 plate)},
   %%review={\MR{0048360 (14,2c)}},
}


\bib{collected2}{book}{
   author={Bernstein, S. N.},
   title={Sobranie so\v cineni\u\i. Tom II. Konstruktivnaya teoriya funkci\u\i\
   [1931--1953] (Collected works. Vol. II. The constructive theory of functions [1931--1953].)},
   language={Russian},
   publisher={Akad. Nauk SSSR},
   place={Moscow},
   date={1954},
   pages={627},
 %  review={\MR{0065452 (16,433o)}},
}





\bib{Bo2002}{article}{
   author={Bondarenko, A. V.},
   title={Jackson type inequality in 3-convex approximation},
   journal={East J. Approx.},
   volume={8},
   date={2002},
   number={3},
   pages={291--302},
   %%issn={1310-6236},
   %%review={\MR{1932335 (2003h:41014)}},
}




\bib{BG2002}{article}{
   author={Bondarenko, A. V.},
   author={Gilewicz, J.  },
   title={A negative result in pointwise 3-convex approximation by
   polynomials},
   language={Ukrainian, with English and Russian summaries},
   journal={Ukra\"\i n. Mat. Zh.},
   volume={61},
   date={2009},
   number={4},
   pages={563--567},
   %%issn={1027-3190},
   %%review={\MR{2588680 (2011a:41012)}},
}

\bib{BP2004}{article}{
   author={Bondarenko, A. V.},
   author={Primak, A. V.},
   title={Negative results in shape-preserving higher-order approximations},
   language={Russian, with Russian summary},
   journal={Mat. Zametki},
   volume={76},
   date={2004},
   number={6},
   pages={812--823},
   %%issn={0025-567X},
   translation={
      journal={Math. Notes},
      volume={76},
      date={2004},
      number={5-6},
      pages={758--769},
      %%issn={0001-4346},
   },
   %%review={\MR{2127492 (2006d:41010)}},
   %%doi={10.1023/B:MATN.0000049675.24438.f2},
}




\bib{B1959}{article}{
   author={Brudnyi, Ju. A.},
   title={Approximation by integral functions on the exterior of a segment
   or on a semi-axis},
   language={Russian},
   journal={Dokl. Akad. Nauk SSSR},
   volume={124},
   date={1959},
   pages={739--742},
%   issn={0002-3264},
%   review={\MR{0101444 (21 \#255)}},
}



\bib{Br1963}{article}{
   author={Brudnyi, Ju. A.},
   title={Generalization of a theorem of A. F. Timan},
   language={Russian},
   journal={Dokl. Akad. Nauk SSSR},
   volume={148},
   date={1963},
   pages={1237--1240},
   %%issn={0002-3264},
   %%review={\MR{0146574 (26 \#4096)}},
}


\bib{Ch1873}{article}{
   label={Che1873},
   author={Chebyshev, P. L.},
   title={On functions having least deviation from zero},
   journal={Prilozh. k XXII Zapis. Akad., N 1. Journal de M. Liouville, II serie, t, XIX, 1874.},
   date={1873},
   language={Russian},
}

\bib{C1955}{article}{
   author={Chebyshev, P. L.},
   title={Selected works},
   language={Russian},
   book={
   publisher={Izdatel'stvo Akademii Nauk SSSR},
     place={Moscow},
   },
   date={1955},
   pages={579--608},
}


\bib{CG1994}{article}{
   author={Cao, Jia Ding},
   author={Gonska, Heinz H.},
   title={Pointwise estimates for higher order convexity preserving
   polynomial approximation},
   journal={J. Austral. Math. Soc. Ser. B},
   volume={36},
   date={1994},
   number={2},
   pages={213--233},
   %%issn={0334-2700},
   %%review={\MR{1312232 (95k:41011)}},
   %%doi={10.1017/S0334270000010365},
}


\bib{CGG1999}{article}{
   author={Cottin, Claudia},
   author={Gavrea, Ioan},
   author={Gonska, Heinz H.},
   author={Kacs\'o, Daniela P.},
   author={Zhou, Ding-Xuan},
   title={Global smoothness preservation and the
   variation-diminishing property},
   journal={J. Inequ. \&Appl.},
   volume={4},
   date={1999},
   pages={91--114},
}



\bib{dah}{article}{
   author={Dahlhaus, R.},
   title={Pointwise approximation by algebraic polynomials},
   journal={J. Approx. Theory},
   volume={57},
   date={1989},
   number={3},
   pages={274--277},
%   issn={0021-9045},
%   review={\MR{999862 (90h:41012)}},
%   doi={10.1016/0021-9045(89)90042-7},
}

\bib{DV1974}{article}{
   author={DeVore, R. A.},
   title={Degree of monotone approximation},
   conference={
      title={Linear operators and approximation, II (Proc. Conf.,
      Oberwolfach Math. Res. Inst., Oberwolfach, 1974)},
   },
   book={
      publisher={Birkh\"auser},
      place={Basel},
   },
   date={1974},
   pages={337--351. Internat. Ser. Numer. Math., Vol. 25},
   %%review={\MR{0382930 (52 \#3812)}},
}


\bib{DV1976}{article}{
   author={DeVore, R. A.},
   title={Degree of approximation},
   conference={
      title={Approximation theory, II (Proc. Internat. Sympos., Univ. Texas,
      Austin, Tex., 1976)},
   },
   book={
      publisher={Academic Press},
      place={New York},
   },
   date={1976},
   pages={117--161},
   %%review={\MR{0440865 (55 \#13733)}},
}



\bib{DV1977}{article}{
   author={DeVore, R. A.},
   title={Monotone approximation by polynomials},
   journal={SIAM J. Math. Anal.},
   volume={8},
   date={1977},
   number={5},
   pages={906--921},
   %%issn={0036-1410},
   %%review={\MR{0510582 (58 \#23252)}},
}




\bib{DY1985}{article}{
   author={DeVore, R. A.},
   author={Yu, X. M.},
   title={Pointwise estimates for monotone polynomial approximation},
   journal={Constr. Approx.},
   volume={1},
   date={1985},
   number={4},
   pages={323--331},
   %%issn={0176-4276},
   %%review={\MR{891762 (88h:41010)}},
   %%doi={10.1007/BF01890039},
}

\bib{Dit2007}{article}{
   author={Ditzian, Z.},
   title={Polynomial approximation and $\omega^r_\varphi(f,t)$ twenty
   years later},
   journal={Surv. Approx. Theory},
   volume={3},
   date={2007},
   pages={106--151},
%   review={\MR{2342231 (2008f:41010)}},
}


\bib{DT}{book}{
   author={Ditzian, Z.},
   author={Totik, V.},
   title={Moduli of smoothness},
   series={Springer Series in Computational Mathematics},
   volume={9},
   publisher={Springer-Verlag},
   place={New York},
   date={1987},
   pages={x+227},
%   isbn={0-387-96536-X},
%   review={\MR{914149 (89h:41002)}},
}


\bib{D1956}{article}{
   author={Dzyadyk, V. K.},
   title={Constructive characterization of functions satisfying the
   condition ${\rm Lip}\,\alpha(0<\alpha<1)$ on a finite segment of
   the real axis},
   language={Russian},
   journal={Izv. Akad. Nauk SSSR. Ser. Mat.},
   volume={20},
   date={1956},
   pages={623--642},
%   issn={0373-2436},
%   review={\MR{0081376 (18,392c)}},
}


\bib{Dz1958}{article}{
   author={Dzyadyk, V. K.},
   title={A further strengthening of Jackson's theorem on the approximation
   of continuous functions by ordinary polynomials},
   language={Russian},
   journal={Dokl. Akad. Nauk SSSR},
   volume={121},
   date={1958},
   pages={403--406},
   %%issn={0002-3264},
   %%review={\MR{0101438 (21 \#249)}},
}


\bib{s-book}{book}{
   author={Dzyadyk, V. K.},
  author={Shevchuk, I. A.},
  title={Theory of Uniform Approximation of Functions by Polynomials},
  publisher={Walter de Gruyter},
  place={Berlin},
  date={2008},
}




\bib{D1994}{article}{
   author={Dzyubenko, G. A.},
   title={A pointwise estimate of a comonotone approximation},
   language={Russian, with English and Ukrainian summaries},
   journal={Ukra\"\i n. Mat. Zh.},
   volume={46},
   date={1994},
   number={11},
   pages={1467--1472},
   %%issn={0041-6053},
   translation={
      journal={Ukrainian Math. J.},
      volume={46},
      date={1994},
      number={11},
      pages={1620--1626 (1996)},
      %%issn={0041-5995},
   },
   %%review={\MR{1356197 (96g:41021)}},
   %%doi={10.1007/BF01058880},
}




\bib{DGS1998}{article}{
   author={Dzyubenko, G. A.},
   author={Gilewicz, J.},
   author={Shevchuk, I. A.},
   title={Piecewise monotone pointwise approximation},
   journal={Constr. Approx.},
   volume={14},
   date={1998},
   number={3},
   pages={311--348},
   %%issn={0176-4276},
   %%review={\MR{1626710 (99f:41008)}},
   %%doi={10.1007/s003659900077},
}

\bib{DGS2002}{article}{
   author={Dzyubenko, G. A.},
   author={Gilewicz, J.},
   author={Shevchuk, I. A.},
   title={Conconvex pointwise approximation},
   language={English, with English and Ukrainian summaries},
   journal={Ukra\"\i n. Mat. Zh.},
   volume={54},
   date={2002},
   number={9},
   pages={1200--1212},
   %%issn={0041-6053},
   translation={
      journal={Ukrainian Math. J.},
      volume={54},
      date={2002},
      number={9},
      pages={1445--1461},
      %%issn={0041-5995},
   },
   %%review={\MR{2016140 (2004i:41011)}},
   %%doi={10.1023/A:1023411817844},
}



\bib{DGS2006}{article}{
   author={Dzyubenko, G. A.},
   author={Gilewicz, J.},
   author={Shevchuk, I. A.},
   title={New phenomena in coconvex approximation},
   language={English, with English and Russian summaries},
   journal={Anal. Math.},
   volume={32},
   date={2006},
   number={2},
   pages={113--121},
   %%issn={0133-3852},
   %%review={\MR{2248067 (2007b:41031)}},
   %%doi={10.1007/s10476-006-0005-x},
}



\bib{DLS2010a}{article}{
   author={Dzyubenko, G. A.},
   author={Leviatan, D.},
   author={Shevchuk, I. A.},
   title={Nikolskii-type estimates for coconvex approximation of functions with one inflection point},
   journal={Jaen J. Approx.},
   volume={2},
   date={2010},
   number={1},
   pages={51--64},
}

\bib{DLS2010b}{article}{
   author={Dzyubenko, G. A.},
   author={Leviatan, D.},
   author={Shevchuk, I. A.},
   title={Coconvex pointwise approximation},
   journal={Supplemento ai Rendiconti del Circolo Matematico di Palermo, Serie II,},
   volume={82},
   date={2010},
   pages={359--374},
}

\bib{DLS1993}{article}{
   author={Dzyubenko, G. A.},
   author={Listopad, V. V.},
   author={Shevchuk, I. A.},
   title={Uniform estimates for monotonic polynomial approximation},
   journal={Ukrainian Math. J.},
   volume={45},
   date={1993},
   number={1},
   pages={40--47},
   %%issn={0041-5995},
   %%review={\MR{1426568 (97k:41008)}},
}


\bib{DZ2004}{article}{
   author={Dzyubenko, G. A.},
   author={Zal{\={\i}}zko, V. D.},
   title={Co-convex approximation of functions that have more than one
   inflection point},
   language={Ukrainian, with English and Ukrainian summaries},
   journal={Ukra\"\i n. Mat. Zh.},
   volume={56},
   date={2004},
   number={3},
   pages={352--365},
   %%issn={0041-6053},
   translation={
      journal={Ukrainian Math. J.},
      volume={56},
      date={2004},
      number={3},
      pages={427--445},
      %%issn={0041-5995},
   },
   %%review={\MR{2105891 (2005h:41043)}},
   %%doi={10.1023/B:UKMA.0000045688.71949.44},
}

\bib{DZ2005}{article}{
   author={Dzyubenko, G. A.},
   author={Zal{\={\i}}zko, V. D.},
   title={Pointwise estimates for the coconvex approximation of
   differentiable functions},
   language={Ukrainian, with English and Ukrainian summaries},
   journal={Ukra\"\i n. Mat. Zh.},
   volume={57},
   date={2005},
   number={1},
   pages={47--59},
   %%issn={0041-6053},
   translation={
      journal={Ukrainian Math. J.},
      volume={57},
      date={2005},
      number={1},
      pages={52--69},
      %%issn={0041-5995},
   },
   %%review={\MR{2190954}},
   %%doi={10.1007/s11253-005-0171-1},
}

\bib{F1959}{article}{
   author={Freud, G.},
   title={\"Uber die Approximation reeller stetigen Funktionen durch
   gew\"ohnliche Polynome},
   language={German},
   journal={Math. Ann.},
   volume={137},
   date={1959},
   pages={17--25},
   %%issn={0025-5831},
   %%review={\MR{0101439 (21 \#250)}},
}

\bib{Ga1996}{article}{
   author={Gavrea, I.},
   title={The approximation of the continuous functions by means of
   some linear positive operators},
   journal={Results in Math.},
   volume={30},
   date={1996},
   pages={55--66},
}

\bib{GGP2003}{article}{
   author={Gavrea, I.},
   author={Gonska, H.},
   author={P{\u a}lt{\u a}nea, R.},
   author={Tachev, G.},
   title={General estimates for the Ditzian-Totik modulus},
   journal={East J. Approx.},
   volume={9},
   date={2003},
   pages={175--194},
}

\bib{gehner}{article}{
   author={Gehner, K. R.},
   title={Characterization theorems for constrained approximation problems
   via optimization theory},
   journal={J. Approximation Theory},
   volume={14},
   date={1975},
   pages={51--76},
%   issn={0021-9045},
%   review={\MR{0372496 (51 \#8703)}},
}



\bib{gel}{article}{
   author={Gelfond, A. O.},
   title={On uniform approximations by polynomials with integral rational
   coefficients},
   language={Russian},
   journal={Uspehi Mat. Nauk (N.S).},
   volume={10},
   date={1955},
   number={1(63)},
   pages={41--65},
%   issn={0042-1316},
%   review={\MR{0070750 (17,30e)}},
}



\bib{GS1996}{article}{
   author={Gilevich, Ya.},
   author={Shevchuk, I. A.},
   title={Comonotone approximation},
   language={Russian, with English and Russian summaries},
   journal={Fundam. Prikl. Mat.},
   volume={2},
   date={1996},
   number={2},
   pages={319--363},
   %%issn={1560-5159},
   %%review={\MR{1793407}},
}

\bib{GY2002}{article}{
   author={Gilewicz, J.},
   author={Yushchenko, L. P.},
   title={A counterexample in coconvex and $q$-coconvex approximations},
   journal={East J. Approx.},
   volume={8},
   date={2002},
   number={2},
   pages={131--144},
   %%issn={1310-6236},
   %%review={\MR{1983845 (2004g:41026)}},
}

\bib{glsw}{article}{
   author={Gonska, H. H.},
   author={Leviatan, D.},
   author={Shevchuk, I. A.},
   author={Wenz, H.-J.},
   title={Interpolatory pointwise estimates for polynomial approximation},
   journal={Constr. Approx.},
   volume={16},
   date={2000},
   number={4},
   pages={603--629},
%   issn={0176-4276},
%   review={\MR{1771698 (2001h:41007)}},
%   doi={10.1007/s003650010008},
}





\bib{gop}{article}{
   author={Gopengauz, I. E.},
   title={A question concerning the approximation of functions on a segment
   and in a region with corners},
   language={Russian},
   journal={Teor. Funkci\u\i\ Funkcional. Anal. i Prilo\v zen. Vyp.},
   volume={4},
   date={1967},
   pages={204--210},
%   issn={0321-4427},
%   review={\MR{0223790 (36 \#6838)}},
}




\bib{HLY1994}{article}{
   author={Hu, Y.},
   author={Leviatan, D.},
   author={Yu, X. M.},
   title={Convex polynomial and spline approximation in $C[-1,1]$},
   journal={Constr. Approx.},
   volume={10},
   date={1994},
   number={1},
   pages={31--64},
   %%issn={0176-4276},
   %%review={\MR{1260358 (95a:41018)}},
   %%doi={10.1007/BF01205165},
}


\bib{ibr1946}{article}{
   author={Ibraghimoff, I.},
   title={Sur la valeur asymptotique de la meilleure approximation d'une
   fonction ayant un point singulier r\'eel},
   language={Russian, with French summary},
   journal={Bull. Acad. Sci. URSS. S\'er. Math. [Izvestia Akad. Nauk SSSR]},
   volume={10},
   date={1946},
%   pages={429--460},
%   review={\MR{0019772 (8,459c)}},
}



\bib{ibr}{article}{
   author={Ibragimov, I. I.},
   title={On the best approzimation by polynomials of the functions
   $[ax+b\vert x\vert ]\vert x\vert ^s$ on the interval $[-1,+1]$},
   language={Russian},
   journal={Izvestiya Akad. Nauk SSSR. Ser. Mat.},
   volume={14},
   date={1950},
   pages={405--412},
%   issn={0373-2436},
%   review={\MR{0037932 (12,331e)}},
}


\bib{I1978}{article}{
   author={Iliev, G. L.},
   title={Exact estimates for partially monotone approximation},
   language={English, with Russian summary},
   journal={Anal. Math.},
   volume={4},
   date={1978},
   number={3},
   pages={181--197},
   %%issn={0133-3852},
   %%review={\MR{514758 (80h:41007)}},
   %%doi={10.1007/BF01908988},
}

\bib{I1978-serdica}{article}{
   author={Iliev, G. L.},
   title={Partially monotone interpolation},
   journal={Serdica},
   volume={4},
   date={1978},
   number={2-3},
   pages={267--276},
%   issn={0204-4110},
%   review={\MR{541821 (81a:41014)}},
}




\bib{I1977}{article}{
   author={Ishisaki, K.},
   title={Jackson-type estimates for monotone approximation},
   journal={Proc. Japan Acad. Ser. A Math. Sci.},
   volume={53},
   date={1977},
   number={5},
   pages={171--173},
%   issn={0386-2194},
%   review={\MR{0460972 (57 \#961)}},
}

\bib{KL1976}{article}{
   author={Kimchi, E.},
   author={Leviatan, D.},
   title={On restricted best approximation to functions with restricted
   derivatives},
   journal={SIAM J. Numer. Anal.},
   volume={13},
   date={1976},
   number={1},
   pages={51--53},
%   issn={0036-1429},
%   review={\MR{0407489 (53 \#11264)}},
}


\bib{KL2003}{article}{
   author={Konovalov, V. N.},
   author={Leviatan, D.},
   title={Shape preserving widths of Sobolev-type classes of $s$-monotone
   functions on a finite interval},
   journal={Israel J. Math.},
   volume={133},
   date={2003},
   pages={239--268},
   %%issn={0021-2172},
   %%review={\MR{1968430 (2003m:41045)}},
   %%doi={10.1007/BF02773069},
}

\bib{K1992}{article}{
   author={Kopotun, K. A.},
   title={Uniform estimates for coconvex approximation of functions by
   polynomials},
   language={Russian},
   journal={Mat. Zametki},
   volume={51},
   date={1992},
   number={3},
   pages={35--46, 143},
   %%issn={0025-567X},
   translation={
      journal={Math. Notes},
      volume={51},
      date={1992},
      number={3-4},
      pages={245--254},
      %%issn={0001-4346},
   },
   %%review={\MR{1172223 (93f:41014)}},
   %%doi={10.1007/BF01206386},
}

\bib{K1994}{article}{
   author={Kopotun, K. A.},
   title={Pointwise and uniform estimates for convex approximation of
   functions by algebraic polynomials},
   journal={Constr. Approx.},
   volume={10},
   date={1994},
   number={2},
   pages={153--178},
   %%issn={0176-4276},
   %%review={\MR{1305916 (95k:41014)}},
   %%doi={10.1007/BF01263061},
}

 \bib{K1995}{article}{
   author={Kopotun, K. A.},
   title={Uniform estimates of monotone and convex approximation of smooth
   functions},
   journal={J. Approx. Theory},
   volume={80},
   date={1995},
   number={1},
   pages={76--107},
   %%issn={0021-9045},
   %%review={\MR{1308595 (95j:41006)}},
   %%doi={10.1006/jath.1995.1005},
}

\bib{kop-sim}{article}{
   author={Kopotun, K. A.},
   title={Simultaneous approximation by algebraic polynomials},
   journal={Constr. Approx.},
   volume={12},
   date={1996},
   number={1},
   pages={67--94},
%   issn={0176-4276},
%   review={\MR{1389920 (97c:41021)}},
%   doi={10.1007/s003659900003},
}




\bib{KL1997}{article}{
   author={Kopotun, K. A.},
   author={Leviatan, D.},
   title={Comonotone polynomial approximation in ${\bf L}_p[-1,1]$,
   $0<p\leq\infty$},
   journal={Acta Math. Hungar.},
   volume={77},
   date={1997},
   number={4},
   pages={301--310},
   %%issn={0236-5294},
   %%review={\MR{1488227 (98m:41011)}},
   %%doi={10.1023/A:1006589610577},
}



\bib{KLS1999}{article}{
   author={Kopotun, K. A.},
   author={Leviatan, D.},
   author={Shevchuk, I. A.},
   title={The degree of coconvex polynomial approximation},
   journal={Proc. Amer. Math. Soc.},
   volume={127},
   date={1999},
   number={2},
   pages={409--415},
   %%issn={0002-9939},
   %%review={\MR{1459130 (99c:41010)}},
   %%doi={10.1090/S0002-9939-99-04452-4},
}

\bib{KLS2005}{article}{
   author={Kopotun, K. A.},
   author={Leviatan, D.},
   author={Shevchuk, I. A.},
   title={Convex polynomial approximation in the uniform norm: conclusion},
   journal={Canad. J. Math.},
   volume={57},
   date={2005},
   number={6},
   pages={1224--1248},
   %%issn={0008-414X},
   %%review={\MR{2178560 (2006f:41018)}},
   %%doi={10.4153/CJM-2005-049-6},
}

\bib{KLS2006}{article}{
   author={Kopotun, K. A.},
   author={Leviatan, D.},
   author={Shevchuk, I. A.},
   title={Coconvex approximation in the uniform norm: the final frontier},
   journal={Acta Math. Hungar.},
   volume={110},
   date={2006},
   number={1-2},
   pages={117--151},
   %%issn={0236-5294},
   %%review={\MR{2198418 (2006i:41014)}},
   %%doi={10.1007/s10474-006-0010-3},
}



\bib{KLS2009}{article}{
   author={Kopotun, K. A.},
   author={Leviatan, D.},
   author={Shevchuk, I. A.},
   title={Are the degrees of best (co)convex and unconstrained polynomial
   approximation the same?},
   journal={Acta Math. Hungar.},
   volume={123},
   date={2009},
   number={3},
   pages={273--290},
   %%issn={0236-5294},
   %%review={\MR{2500917 (2010c:41003)}},
   %%doi={10.1007/s10474-009-8111-4},
}



\bib{KLS2010}{article}{
   author={Kopotun, K. A.},
   author={Leviatan, D.},
     author={Shevchuk, I. A.},
    title={Are the degrees of best (co)convex and unconstrained polynomial
   approximation the same? II},
   language={English, with English and Ukrainian summaries},
   journal={Ukra\"\i n. Mat. Zh.},
   volume={62},
   date={2010},
   number={3},
   pages={369--386},
   translation={
      journal={Ukrainian Math. J.},
      volume={62},
      date={2010},
      number={3},
      pages={420--440},
   },
}




\bib{KL1994}{article}{
   author={Kopotun, K. A.},
   author={Listopad, V. V.},
   title={Remarks on monotone and convex approximation by algebraic
   polynomials},
   language={English, with English and Ukrainian summaries},
   journal={Ukra\"\i n. Mat. Zh.},
   volume={46},
   date={1994},
   number={9},
   pages={1266--1270},
   %%issn={0041-6053},
   translation={
      journal={Ukrainian Math. J.},
      volume={46},
      date={1994},
      number={9},
      pages={1393--1398 (1996)},
      %%issn={0041-5995},
   },
   %%review={\MR{1319576 (95k:41027)}},
   %%doi={10.1007/BF01059430},
}

\bib{L1957}{article}{
   author={Lebed', G. K.},
   title={Inequalities for polynomials and their derivatives},
   language={Russian},
   journal={Dokl. Akad. Nauk SSSR (N.S.)},
   volume={117},
   date={1957},
   pages={570--572},
%   issn={0002-3264},
%   review={\MR{0094652 (20 \#1165)}},
}



\bib{Le1986}{article}{
   author={Leviatan, D.},
   title={Pointwise estimates for convex polynomial approximation},
   journal={Proc. Amer. Math. Soc.},
   volume={98},
   date={1986},
   number={3},
   pages={471--474},
   %%issn={0002-9939},
   %%review={\MR{857944 (88i:41010)}},
   %%doi={10.2307/2046205},
}

\bib{LRS2011}{article}{
   author={Leviatan, D.},
   author={Radchenko, D.},
   author={Shevchuk, I. A.},
   title={},
   journal={preprint},
}

\bib{LS1995}{article}{
   author={Leviatan, D.},
   author={Shevchuk, I. A.},
   title={Counterexamples in convex and higher order constrained
   approximation},
   journal={East J. Approx.},
   volume={1},
   date={1995},
   number={3},
   pages={391--398},
   %%issn={1310-6236},
   %%review={\MR{1404356 (97e:41044)}},
}

\bib{LS1997}{article}{
   author={Leviatan, D.},
   author={Shevchuk, I. A.},
   title={Some positive results and counterexamples in comonotone
   approximation},
   journal={J. Approx. Theory},
   volume={89},
   date={1997},
   number={2},
   pages={195--206},
   %%issn={0021-9045},
   %%review={\MR{1447838 (99a:41023)}},
   %%doi={10.1006/jath.1997.3038},
}

\bib{LS1998}{article}{
   author={Leviatan, D.},
   author={Shevchuk, I. A.},
   title={Monotone approximation estimates involving the third modulus of
   smoothness},
   conference={
      title={Approximation theory IX, Vol. I.},
      address={Nashville, TN},
      date={1998},
   },
   book={
      series={Innov. Appl. Math.},
      publisher={Vanderbilt Univ. Press},
      place={Nashville, TN},
   },
   date={1998},
   pages={223--230},
   %%review={\MR{1743008 (2001b:41023)}},
}

\bib{LS1999a}{article}{
   author={Leviatan, D.},
   author={Shevchuk, I. A.},
   title={Some positive results and counterexamples in comonotone
   approximation. II},
   journal={J. Approx. Theory},
   volume={100},
   date={1999},
   number={1},
   pages={113--143},
   %%issn={0021-9045},
   %%review={\MR{1710556 (2000f:41026)}},
   %%doi={10.1006/jath.1999.3334},
}

 \bib{LS1999}{article}{
   author={Leviatan, D.},
   author={Shevchuk, I. A.},
   title={Constants in comonotone polynomial approximation---a survey},
   conference={
      title={New developments in approximation theory},
      address={Dortmund},
      date={1998},
   },
   book={
      series={Internat. Ser. Numer. Math.},
      volume={132},
      publisher={Birkh\"auser},
      place={Basel},
   },
   date={1999},
   pages={145--158},
   %%review={\MR{1724917 (2000h:41024)}},
}

 \bib{LS2000}{article}{
   author={Leviatan, D.},
   author={Shevchuk, I. A.},
   title={More on comonotone polynomial approximation},
   journal={Constr. Approx.},
   volume={16},
   date={2000},
   number={4},
   pages={475--486},
   %%issn={0176-4276},
   %%review={\MR{1771692 (2001h:41029)}},
   %%doi={10.1007/s003650010003},
}

 \bib{LS2002}{article}{
   author={Leviatan, D.},
   author={Shevchuk, I. A.},
   title={Coconvex approximation},
   journal={J. Approx. Theory},
   volume={118},
   date={2002},
   number={1},
   pages={20--65},
   %%issn={0021-9045},
   %%review={\MR{1928255 (2003f:41027)}},
   %%doi={10.1006/jath.2002.3695},
}

\bib{LS2003}{article}{
   author={Leviatan, D.},
   author={Shevchuk, I. A.},
   title={Coconvex polynomial approximation},
   journal={J. Approx. Theory},
   volume={121},
   date={2003},
   number={1},
   pages={100--118},
   %%issn={0021-9045},
   %%review={\MR{1962998 (2004b:41018)}},
   %%doi={10.1016/S0021-9045(02)00045-X},
}

% \bib{L1953}{book}{
%   author={Lorentz, G. G.},
%   title={Bernstein polynomials},
%   series={Mathematical Expositions, no. 8},
%   publisher={University of Toronto Press},
%   place={Toronto},
%   date={1953},
%   pages={x+130},
%   %%review={\MR{0057370 (15,217a)}},
%}

\bib{li}{article}{
   author={Li, W.},
   title={On Timan type theorems in algebraic polynomial approximation},
   language={Chinese},
   journal={Acta Math. Sinica},
   volume={29},
   date={1986},
   number={4},
   pages={544--549},
%   issn={0583-1431},
%   review={\MR{867708 (88e:41022)}},
}


\bib{lim}{article}{
   author={Lim, K. P.},
   title={Note on monotone approximation},
   journal={Bull. London Math. Soc.},
   volume={3},
   date={1971},
   pages={366--368},
%   issn={0024-6093},
%   review={\MR{0293303 (45 \#2380)}},
}



 \bib{L1972}{article}{
   author={Lorentz, G. G.},
   title={Monotone approximation},
   conference={
      title={Inequalities, III (Proc. Third Sympos., Univ. California, Los
      Angeles, Calif., 1969; dedicated to the memory of Theodore S.
      Motzkin)},
   },
   book={
      publisher={Academic Press},
      place={New York},
   },
   date={1972},
   pages={201--215},
   %%review={\MR{0346375 (49 \#11100)}},
}


\bib{LZ1968}{article}{
   author={Lorentz, G. G.},
   author={Zeller, K. L.},
   title={Degree of approximation by monotone polynomials. I},
   journal={J. Approximation Theory},
   volume={1},
   date={1968},
   pages={501--504},
   %%issn={0021-9045},
   %%review={\MR{0239342 (39 \#699)}},
}

\bib{LZ1969}{article}{
   author={Lorentz, G. G.},
   author={Zeller, K. L.},
   title={Degree of approximation by monotone polynomials. II},
   journal={J. Approximation Theory},
   volume={2},
   date={1969},
   pages={265--269},
   %%issn={0021-9045},
   %%review={\MR{0244677 (39 \#5991)}},
}

\bib{LZ1970}{article}{
   author={Lorentz, G. G.},
   author={Zeller, K. L.},
   title={Monotone approximation by algebraic polynomials},
   journal={Trans. Amer. Math. Soc.},
   volume={149},
   date={1970},
   pages={1--18},
%   issn={0002-9947},
%   review={\MR{0285843 (44 \#3060)}},
}

\bib{LR1971}{article}{
   author={Lorentz, R. A.},
   title={Uniqueness of best approximation by monotone polynomials},
   journal={J. Approximation Theory},
   volume={4},
   date={1971},
   pages={401--418},
%   issn={0021-9045},
%   review={\MR{0291688 (45 \#779)}},
}



\bib{MR1978}{article}{
   author={Myers, D.-C.},
   author={Raymon, L.},
   title={Exact comonotone approximation},
   journal={J. Approx. Theory},
   volume={24},
   date={1978},
   number={1},
   pages={35--50},
   issn={0021-9045},
%   review={\MR{510920 (80a:41005)}},
%   doi={10.1016/0021-9045(78)90035-7},
}




\bib{NP2005}{article}{
   author={Nesterenko, O. N.},
   author={Petrova, T. O.},
   title={On a problem for co-monotone approximation},
   language={Ukrainian, with English and Ukrainian summaries},
   journal={Ukra\"\i n. Mat. Zh.},
   volume={57},
   date={2005},
   number={10},
   pages={1424--1429},
   %%issn={0041-6053},
   translation={
      journal={Ukrainian Math. J.},
      volume={57},
      date={2005},
      number={10},
      pages={1667--1673},
      %%issn={0041-5995},
   },
   %%review={\MR{2219774 (2007a:41029)}},
   %%doi={10.1007/s11253-006-0021-9},
}



\bib{N1979}{article}{
   author={Newman, D. J.},
   title={Efficient co-monotone approximation},
   journal={J. Approx. Theory},
   volume={25},
   date={1979},
   number={3},
   pages={189--192},
   %%issn={0021-9045},
   %%review={\MR{531408 (80j:41011)}},
   %%doi={10.1016/0021-9045(79)90009-1},
}


\bib{NPR1972}{article}{
   author={Newman, D. J.},
   author={Passow, Eli},
   author={Raymon, Louis},
   title={Piecewise monotone polynomial approximation},
   journal={Trans. Amer. Math. Soc.},
   volume={172},
   date={1972},
   pages={465--472},
%   issn={0002-9947},
%   review={\MR{0310506 (46 \#9604)}},
}


\bib{N1946}{article}{
   author={Nikolskii, S. M.},
   title={On the best approximation of functions satisfying Lipschitz's
   conditions by polynomials},
   language={Russian, with English summary},
   journal={Bull. Acad. Sci. URSS. S\'er. Math. [Izvestia Akad. Nauk SSSR]},
   volume={10},
   date={1946},
%   pages={295--322},
%   review={\MR{0017439 (8,153i)}},
}





\bib{NY2003}{article}{
   author={Nissim, R.},
   author={Yushchenko, L. P.},
   title={Negative result for nearly $q$-convex approximation},
   journal={East J. Approx.},
   volume={9},
   date={2003},
   number={2},
   pages={209--213},
   %%issn={1310-6236},
   %%review={\MR{1991851 (2004d:41029)}},
}

\bib{PR1974}{article}{
   author={Passow, E.},
   author={Raymon, L.},
   title={Monotone and comonotone approximation},
   journal={Proc. Amer. Math. Soc.},
   volume={42},
   date={1974},
   pages={390--394},
%   issn={0002-9939},
%   review={\MR{0336176 (49 \#952)}},
}



\bib{PRR1974}{article}{
   author={Passow, E.},
   author={Raymon, L.},
   author={Roulier, J. A.},
   title={Comonotone polynomial approximation},
   journal={J. Approximation Theory},
   volume={11},
   date={1974},
   pages={221--224},
%   issn={0021-9045},
%   review={\MR{0352807 (50 \#5293)}},
}

\bib{PRo1976}{article}{
   author={Passow, E.},
   author={Roulier, J. A.},
   title={Negative theorems on generalized convex approximation},
   journal={Pacific J. Math.},
   volume={65},
   date={1976},
   number={2},
   pages={437--447},
%   issn={0030-8730},
%   review={\MR{0422966 (54 \#10950)}},
}




\bib{PS2000}{article}{
   author={Pleshakov, M. G.},
   author={Shatalina, A. V.},
   title={Piecewise coapproximation and the Whitney inequality},
   journal={J. Approx. Theory},
   volume={105},
   date={2000},
   number={2},
   pages={189--210},
   %%issn={0021-9045},
   %%review={\MR{1775144 (2002d:41014)}},
   %%doi={10.1006/jath.2000.3458},
}

\bib{PS1974}{article}{
   author={Popov, V.},
   author={Sendov, Bl.},
   title={Approximation of monotone functions by monotone polynomials in
   Hausdorff metric},
   journal={Rev. Anal. Num\'er. Th\'eorie Approximation},
   volume={3},
   date={1974},
   number={1},
   pages={79--88},
%   issn={0025-5505},
%   review={\MR{0380199 (52 \#1099)}},
}






\bib{p1934}{article}{
  author={Popoviciu, T.},
  title={Sur l'approximation des fonctions convexes d'ordre sup\'erieur},
  language={French},
  journal={Mathematica},
  volume={10},
  date={1934},
  pages={49--54},
}

\bib{R1968}{article}{
   author={Roulier, J. A.},
   title={Monotone approximation of certain classes of functions},
   journal={J. Approximation Theory},
   volume={1},
   date={1968},
   pages={319--324},
%   issn={0021-9045},
%   review={\MR{0236580 (38 \#4875)}},
}

\bib{R1971}{article}{
   author={Roulier, J. A.},
   title={Linear operators invariant on nonnegative monotone functions},
   journal={SIAM J. Numer. Anal.},
   volume={8},
   date={1971},
   pages={30--35},
%   issn={0036-1429},
%   review={\MR{0283465 (44 \#696)}},
}

\bib{R1973}{article}{
   author={Roulier, J. A.},
   title={Polynomials of best approximation which are monotone},
   journal={J. Approximation Theory},
   volume={9},
   date={1973},
   pages={212--217},
%   issn={0021-9045},
%   review={\MR{0352793 (50 \#5279)}},
}

\bib{R1975}{article}{
   author={Roulier, J. A.},
   title={Some remarks on the degree of monotone approximation},
   journal={J. Approximation Theory},
   volume={14},
   date={1975},
   number={3},
   pages={225--229},
%   issn={0021-9045},
%   review={\MR{0377375 (51 \#13547)}},
}


\bib{R1976}{article}{
   author={Roulier, J. A.},
   title={Negative theorems on monotone approximation},
   journal={Proc. Amer. Math. Soc.},
   volume={55},
   date={1976},
   number={1},
   pages={37--43},
%   issn={0002-9939},
%   review={\MR{0393969 (52 \#14776)}},
}



\bib{S-doklady}{article}{
   author={Shevchuk, I. A.},
   title={On co-approximation of monotone functions},
   language={Russian},
   journal={Dokl. Akad. Nauk SSSR},
   volume={308},
   date={1989},
   number={3},
   pages={537--541},
   %%issn={0002-3264},
   translation={
      journal={Soviet Math. Dokl.},
      volume={40},
      date={1990},
      number={2},
      pages={349--354},
      %%issn={0197-6788},
   },
   %%review={\MR{1021110 (91c:41050)}},
}

 \bib{S1992}{article}{
   author={Shevchuk, I. A.},
   title={Approximation of monotone functions by monotone polynomials},
   language={Russian, with Russian summary},
   journal={Mat. Sb.},
   volume={183},
   date={1992},
   number={5},
   pages={63--78},
   %%issn={0368-8666},
   translation={
      journal={Russian Acad. Sci. Sb. Math.},
      volume={76},
      date={1993},
      number={1},
      pages={51--64},
      %%issn={1064-5616},
   },
   %%review={\MR{1184310 (94i:41022)}},
   %%doi={10.1070/SM1993v076n01ABEH003401},
}

 \bib{M1992}{book}{
   author={Shevchuk, I. A.},
   title={Polynomial approximation and traces of functions continuous on a segment},
   publisher={Naukova Dumka, Kiev},
    language={Russian},
   date={1992},
}


\bib{S1995}{article}{
   author={Shevchuk, I. A.},
   title={Whitney's inequality and coapproximation},
   booktitle={Proceedings of the XIX Workshop on Function Theory (Beloretsk,
   1994)},
   journal={East J. Approx.},
   volume={1},
   date={1995},
   number={4},
   pages={479--500},
%   issn={1310-6236},
%   review={\MR{1407979 (97e:41046)}},
}

\bib{S1996}{article}{
   author={Shevchuk, I. A.},
   title={One example in monotone approximation},
   journal={J. Approx. Theory},
   volume={86},
   date={1996},
   number={3},
   pages={270--277},
   %%issn={0021-9045},
   %%review={\MR{1405980 (97e:41047)}},
   %%doi={10.1006/jath.1996.0068},
}

\bib{S1965}{article}{
   author={Shisha, O.},
   title={Monotone approximation},
   journal={Pacific J. Math.},
   volume={15},
   date={1965},
   pages={667--671},
   %%issn={0030-8730},
   %%review={\MR{0185334 (32 \#2802)}},
}



\bib{S1979}{article}{
   author={Shvedov, A. S.},
   title={Jackson's theorem in $L^{p}$, $0<p<1$, for algebraic
   polynomials and orders of comonotone approximations},
   language={Russian},
   journal={Mat. Zametki},
   volume={25},
   date={1979},
   number={1},
   pages={107--117, 159},
   %%issn={0025-567X},
   %%review={\MR{527004 (81c:41017)}},
}

\bib{S1980}{article}{
   author={Shvedov, A. S.},
   title={Comonotone approximation of functions by polynomials},
   language={Russian},
   journal={Dokl. Akad. Nauk SSSR},
   volume={250},
   date={1980},
   number={1},
   pages={39--42},
   %%issn={0002-3264},
   %%review={\MR{556116 (81d:41026)}},
}




\bib{S1981}{article}{
   author={Shvedov, A. S.},
   title={Orders of coapproximations of functions by algebraic polynomials},
   language={Russian},
   journal={Mat. Zametki},
   volume={29},
   date={1981},
   number={1},
   pages={117--130, 156},
   %%issn={0025-567X},
   %%review={\MR{604156 (82c:41009)}},
}

\bib{S1981-co}{article}{
   author={Shvedov, A. S.},
   title={Co-approximation of piecewise monotone functions by polynomials},
   language={Russian},
   journal={Mat. Zametki},
   volume={30},
   date={1981},
   number={6},
   pages={839--846, 958},
%   issn={0025-567X},
%   review={\MR{641657 (83c:41018)}},
}




\bib{T1951}{article}{
   author={Timan, A. F.},
   title={A strengthening of Jackson's theorem on the best approximation of
   continuous functions by polynomials on a finite segment of the real axis},
   language={Russian},
   journal={Doklady Akad. Nauk SSSR (N.S.)},
   volume={78},
   date={1951},
   pages={17--20},
   %%review={\MR{0041276 (12,823g)}},
}

\bib{T1957}{article}{
   author={Timan, A. F.},
   title={Converse theorems in the constructive theory of functions given on
   a finite segment of the real axis},
   language={Russian},
   journal={Dokl. Akad. Nauk SSSR (N.S.)},
   volume={116},
   date={1957},
   pages={762--765},
%   issn={0002-3264},
%   review={\MR{0092891 (19,1175e)}},
}




 \bib{T1994}{book}{
   author={Timan, A. F.},
   title={Theory of approximation of functions of a real variable},
   note={Translated from the Russian by J. Berry;
   Translation edited and with a preface by J. Cossar;
   Reprint of the 1963 English translation},
   publisher={Dover Publications Inc.},
   place={New York},
   date={1994},
   pages={viii+631},
%   isbn={0-486-67830-X},
%   review={\MR{1262128 (94j:41001)}},
}

\bib{Trigub62}{article}{
   author={Trigub, R. M.},
   title={Approximation of functions by polynomials with integer
   coefficients},
   language={Russian},
   journal={Izv. Akad. Nauk SSSR Ser. Mat.},
   volume={26},
   date={1962},
   pages={261--280},
%   issn={0373-2436},
%   review={\MR{0136912 (25 \#373)}},
}




\bib{WZ1992}{article}{
   author={Wu, X.},
   author={Zhou, S. P.},
   title={On a counterexample in monotone approximation},
   journal={J. Approx. Theory},
   volume={69},
   date={1992},
   number={2},
   pages={205--211},
   %%issn={0021-9045},
   %%review={\MR{1160255 (93a:41034)}},
   %%doi={10.1016/0021-9045(92)90143-C},
}



\bib{yu85}{article}{
   author={Yu, X. M.},
   title={Pointwise estimate for algebraic polynomial approximation},
   journal={Approx. Theory Appl.},
   volume={1},
   date={1985},
   number={3},
   pages={109--114},
%   issn={1000-9221},
%   review={\MR{816613 (87g:41042)}},
}




\bib{Y2002}{article}{
   author={Yushchenko, L. P.},
   title={A counterexample in convex approximation},
   language={Ukrainian, with English and Ukrainian summaries},
   journal={Ukra\"\i n. Mat. Zh.},
   volume={52},
   date={2000},
   number={12},
   pages={1715--1721},
   %%issn={0041-6053},
   translation={
      journal={Ukrainian Math. J.},
      volume={52},
      date={2000},
      number={12},
      pages={1956--1963 (2001)},
      %%issn={0041-5995},
   },
   %%review={\MR{1834633}},
   %%doi={10.1023/A:1010420313286},
}

\bib{Z1973}{article}{
   author={Zeller, K. L.},
   title={Monotone approximation},
   conference={
      title={Approximation theory (Proc. Internat. Sympos., Univ. Texas,
      Austin, Tex., 1973)},
   },
   book={
      publisher={Academic Press},
      place={New York},
   },
   date={1973},
   pages={523--525},
%   review={\MR{0372482 (51 \#8689)}},
}




\bib{Z1993}{article}{
   author={Zhou, S. P.},
   title={On comonotone approximation by polynomials in $L^p$ space},
   journal={Analysis},
   volume={13},
   date={1993},
   number={4},
   pages={363--376},
%   issn={0174-4747},
%   review={\MR{1256520 (95d:41021)}},
}


\end{biblist}
\end{bibsection}
